\documentclass[a4paper,11pt]{article}
\usepackage{makeidx}
\usepackage[mathscr]{eucal}
\usepackage[dvips]{color}
\usepackage{amssymb, amsfonts, amsmath, amsthm, graphics} 
\usepackage{here}
\usepackage{tasks}
\DeclareMathOperator{\sech}{sech}

\newtheorem{definition}{Definition}[section]
\newtheorem{proposition}{Proposition}[section]
\newtheorem{theorem}{Theorem}[section]

\newtheorem{lemma}[proposition]{Lemma}
\newtheorem{remark}{Remark}[section]

\newtheorem*{thank}{Acknowledgments}

\newcommand{\N}{\mathbb{N}}

\newcommand{\R}{\mathbb{R}}

\newcommand{\T}{\mathbb{T}}

\newcommand{\Z}{\mathbb{Z}}

\newcommand{\boC}{\mathcal{C}}

\newcommand{\boF}{\mathcal{F}}

\numberwithin{equation}{section}

\author{}
\date{}
\begin{document}
\title{Transverse stability of line soliton and 
characterization of ground state for wave guide Schr\"{o}dinger equations}
\author{\renewcommand{\thefootnote}{\arabic{footnote}} Yakine Bahri\footnotemark[12], \ Slim Ibrahim\footnotemark[12] \ and Hiroaki Kikuchi\footnotemark[3]}

\footnotetext[1]{Department of Mathematics and Statistics, University of Victoria 3800 Finnerty Road, Victoria, B.C., Canada V8P 5C2 . E-mail: {\tt ybahri@uvic.ca}}

\footnotetext[2]{Pacific Institute for Mathematical Sciences, 4176-2207 Main Mall,
Vancouver, British Columbia
Canada, V6T 1Z4. E-mail: {\tt ibrahims@uvic.ca}}

\footnotetext[3]{Department of Mathematics, Tsuda University 2-1-1 Tsuda-machi, Kodaira-shi, Tokyo 187-8577, JAPAN. E-mail: {\tt hiroaki@tsuda.ac.jp}}

\maketitle
\begin{center}
\textit{In the memory of Walter Craig}
\end{center}
\begin{abstract}
In this paper, we study the transverse stability of the 
line Schr\"{o}dinger soliton under a 
full wave guide Schr\"{o}dinger flow on a cylindrical domain $\mathbb R\times\mathbb T$. 
When the nonlinearity is of power type 
$|\psi|^{p-1}\psi$ with $p>1$, 
we show that there exists a critical 
frequency $\omega_{p} >0$ such that 
the line standing wave is stable for $0<\omega < \omega_{p}$ 
and unstable for $\omega > \omega_{p}$. Furthermore, we characterize the ground state of the wave guide Schr\"{o}dinger equation. More precisely, we prove that there exists $\omega_{*} \in (0, \omega_{p}]$ such that the ground states coincide with the line standing waves for $\omega \in (0, \omega_{*}]$ and are different from the line standing waves 
for $\omega \in (\omega_{*}, \infty)$. 
\end{abstract}

\section{Introduction}
We consider the following wave guide Schr\"{o}dinger 
equation: 
\begin{equation} 
\label{WS}
i \partial_{t} \psi + \partial_{xx} \psi - |D_{y}| \psi 
+ |\psi|^{p-1} \psi = 0, \qquad 
\mbox{in $\mathbb{R}_t \times \mathbb{R}_x \times \mathbb{T}_y$}, 
\end{equation}
where $1 < p< 5,$ $|D_{y}| := \sqrt{- \partial_{yy}}$ and $\mathbb{T} = \mathbb{R}/2\pi \Z $. The equation is Hamiltonian and its Hamiltonian, the energy, is given by
\begin{align}
& \mathcal{H}(u) 
:= \frac{1}{2} \int_{\mathbb{R} \times \mathbb{T}}
\left(|\partial_{x} u(x, y)|^{2} + 
|D_{y}| u(x, y) \overline{u(x, y)}\right) dxdy 
%\\ & \quad
- \frac{1}{p+1}\int_{\mathbb{R} \times \mathbb{T}} |u(x, y)|^{p+1} dxdy. 
\label{hamiltonian}
\end{align} 
The energy as well as the mass i.e.
\begin{equation}
\mathcal{M}(u) := \frac{1}{2} \|u\|_{L_{x, y}^{2}(\mathbb{R} \times \mathbb{T})}^{2}, \label{mass}
\end{equation}
are conservation laws which means 
\begin{equation}
\label{conserv}
\mathcal{H}(\psi(t)) = \mathcal{H}(\psi_{0}), \qquad \mathcal{M}(\psi(t)) = \mathcal{M}(\psi_{0}),
\end{equation}
where $\psi_{0} $ is the initial condition.

The fractional Laplacian $|D_{y}|$ 
can be defined as follow. 
For a function $f(x, y)$ defined on $\mathbb{R} \times \mathbb{T}$, we write its Fourier series 
expansion with respect to $y$ as 
\begin{equation*}%\label{coefficient}
f(x, y) = \sum_{n \in \mathbb{Z}} f_{n}(x) e^{i n y}, 
\end{equation*} 
where each Fourier coefficient $f_n(x)$ is defined by 
\begin{equation*} %\label{fourier-coefficient}
f_{n}(x) = \frac{1}{2\pi} \int_{-\pi}^{\pi} f(x, y) e^{- i n y} dy. 
\end{equation*}
Then, the operator $|D_{y}|$ is defined by 
\begin{equation*}
|D_{y}| f(x, y) = \sum_{n \in \mathbb{Z}} |n| f_{n}(x) e^{i n y}. 
\end{equation*}
We define a Hilbert space $X$ by 
\begin{equation} \label{energyspace}
X = H^{1}_{x}L^{2}_{y} \cap L^{2}_{x}H^{\frac{1}{2}}_{y}(\mathbb{R} \times \mathbb{T}). 
\end{equation}
endowed with the norm 
\begin{equation*}
\|u\|_{X} := \left\{ 
\int_{\mathbb{R} \times \mathbb{T}} 
\left(|\partial_{x} u(x, y)|^{2} + |D_{y}|u(x, y) \overline{u(x, y)} 
+ |u(x, y)|^{2}\right) dxdy
\right\}^{\frac{1}{2}}. 
\end{equation*}

%Concerning the Cauchy problem, 
%we have the following:
%\begin{theorem} \label{assum-1}
%Let $s > 1/2$. 
%For any $\psi_{0} \in H^{1}_{x}L^{2}_{y} \cap L^{2}_{x}H^{s}_{y}(\mathbb{R} \times \mathbb{T})$, there exist $T >0$ and 
%a unique local solution $\psi \in C([0, T), 
%H^{1}_{x}L^{2}_{y} \cap L^{2}_{x}H^{s}_{y}(\mathbb{R} \times \mathbb{T}))$ to \eqref{WS} 
%with $\psi|_{t=0} = \psi_{0}$. 
%Moreover, as long as the solution exists, the following conservation laws hold; 
%
%where $I_{\max}$ is the maximal existence time and 
%\end{theorem}
Equation \eqref{WS} was first considered by Xu~\cite{Xu}, where the author studied the large time behavior of small solutions to the cubic defocusing wave guide Schr\"{o}dinger equation on spatially cylindrical domain $\mathbb{R}_{x} \times \mathbb{T}_{y}$. The author obtained modified scattering result for small smooth initial data. In \cite{BIK}, we studied the Cauchy problem associated to \eqref{WS} in the plane $\R ^2$ with initial data in $L^{2}_{x}H^{s}_{y}(\mathbb{R}^{2})$ and $s> 1/2$. We also showed both the existence and orbital stability of the ground states. In addition, we proved that traveling waves exist and converge to zero as the velocity tends to $1$. 

\begin{remark}
As pointed out by Xu in ~\cite{Xu} and by the authors in \cite{BIK}, even the local well-posedness of
the Cauchy problem of \eqref{WS} in the energy space 
$X$ still remains an interesting unresolved problem. Standard arguments do not seem to work directly. We refer interested readers to Remark 1.6 in \cite{BIK} for more details about these technical difficulties. 
%In \cite{BIK}, the equation \eqref{WS} is locally 
%well-posed in $L_{x}^{2}H_{y}^{\frac{1}{2} + \varepsilon}(\mathbb{R}^{2})$ for any $\varepsilon %>0$. 
%Therefore, we may expect that the equation \eqref{WS} is locally well-posed in the energy %space $X$. 
\end{remark}
A standing wave to \eqref{WS} is a solution which has the form 
$\psi(t, x, y) = e^{i \omega t}\varphi(x, y)\; (\omega >0)$. 
Then, we see that $\varphi(x, y)$ satisfies the following elliptic equation: 
\begin{equation} \label{sp}
- \partial_{xx} \varphi + |D_{y}| \varphi + \omega \varphi 
- |\varphi|^{p-1}\varphi = 0 \qquad 
\mbox{in $\mathbb{R} \times \mathbb{T}$}. 
\end{equation}
Recall that 
\begin{equation} \label{line-sp}
- \partial_{xx} R_{\omega} + \omega R_{\omega} - 
|R_{\omega}|^{p-1} R_{\omega} = 0 \qquad \mbox{in $\mathbb{R}$}, 
\end{equation}
has a unique positive solution given by
\begin{equation}
\label{def-R-omega}
R_{\omega}(x) = \left(\frac{(p+1)\omega}{2}\right)^{\frac{1}{p-1}} 
\sech^{\frac{2}{p-1}}\left(\frac{(p-1)\sqrt{\omega}}{2} x\right).
\end{equation}
We refer to \cite[Theorem 8.1.6]{Cazenave} for the existence and uniqueness of $R_\omega$.
Observe that $e^{i \omega t}R_{\omega}(x)$ becomes a standing wave of not only \eqref{WS} but also of the following Schr\"{o}dinger equations: 
\begin{equation} \label{1d-nls}
i \partial_{t} \psi +\partial_{xx} \psi + |\psi|^{p-1} \psi = 0 \qquad 
\mbox{in $\mathbb{R} \times \mathbb{R}$}. 
\end{equation}
In addition, clearly $R_{\omega}\in X$ and satisfies \eqref{sp}.
% 
% 
%Putting $\widetilde{R}_{\omega}(x, y) = R_{\omega}(x)$ for all $y \in \mathbb{T}$, we see that $\widetilde{R}_{\omega}(x, y) \in X$ and satisfies \eqref{sp}. 
%If there is no fear of confusion, 
%we denote $\widetilde{R}_{\omega}$ by $R_{\omega}$, for simplicity. 
Here, we aim to study the transverse stability of the line standing wave $e^{i \omega t} R_{\omega}(x)$ under the flow of \eqref{WS}. It consists of studying two-dimensional perturbations around the one-dimensional nonlinear Schr\"odinger standing wave. The transverse stability defined as follows:

\begin{definition}
Let $s > 1/2$. 
We say that the standing wave $e^{i \omega t} \varphi_{\omega}(x, y)$ 
is transverse orbitally stable in 
$H^{1}_{x}L^{2}_{y} \cap L^{2}_{x}H^{s}_{y}(\mathbb{R} \times \mathbb{T})$ 
if for any $\varepsilon > 0$, there exists 
$\delta > 0$ such that for all $\psi_{0} \in 
H^{1}_{x}L^{2}_{y} \cap L^{2}_{x}H^{s}_{y}(\mathbb{R} \times \mathbb{T})$ 
satisfying $\|\psi_{0} - R_{\omega}\|_{X} < \delta$, 
the solution $\psi(t)$ to \eqref{WS} 
with the initial data $\psi(0) = \psi_{0}$ exists globally in time and satisfies 
\begin{equation*}
\sup_{t >0} \inf_{\theta \in \mathbb{R}, z \in \mathbb{R}} 
\|\psi(t, \cdot) - e^{i \theta} R_{\omega} (\cdot - z)\|_{X} < \varepsilon. 
\end{equation*}
\end{definition}
The orbital stability of the line standing wave $e^{i \omega t} R_{\omega}$ under the flow of \eqref{1d-nls} has been intensively studied in the literature. 
Results can be summarized as follows
\begin{itemize}
\item Cazenave and Lions proved in \cite{Caz-lions} that $e^{i \omega t} R_{\omega}$ is stable for $1<p<5$. 
%A second approach was introduced by Weinstein in \cite{Wein1, Wein2}.
\item Weinstein showed in \cite{Wein1} that $e^{i \omega t} R_{\omega}$ is unstable 
for $p=5$.
\item Berestycki and Cazenave proved in \cite{bc1} that 
$e^{i \omega t} R_{\omega}$ is unstable for $p > 5$.
\end{itemize}
where $p=5$ corresponds to the $L^2$-critical exponent. Note that the above stability and instability results can also be obtained 
by using the abstract theory provided by Grillakis, Shatah and Strauss 
in \cite{GSS, GSS2} except for $p=5$ where we refer to Weinstein~\cite{Wein1}.

The transverse stability question has also been widely studied for different models in PDEs. 
%It appears for example in surface gravity-capillary waves in a deep fluid. 
For example, in \cite{Pelinov} Pelinovsky showed instability band of a deep-water soliton of the hyperbolic nonlinear Schr\"{o}dinger equation. 
Later on, Deconinck, Pelinovsky and Carter in \cite{deco-peli} described a wave train in deep water by a one-dimensional soliton solution of a two dimensional
nonlinear Schr\"{o}dinger equation. 
Numerically, they examined the instability of such a
wave train with respect to transverse perturbations. 
Rousset and Tzvetkov~\cite{R-T-4} 
proved the linear and nonlinear instability of the line standing waves to water waves with respect to transverse perturbations (see also the references therein). Rousset and Tzvetkov studied in a series of papers \cite{R-T-1, R-T-2, R-T-3} the nonlinear long time instability of the KdV standing wave (with respect to periodic transverse perturbations) under a KP-I flow and the transverse nonlinear instability of standing waves for the cubic nonlinear Schr\"{o}dinger equation. 
They provided in \cite{R-T-5} a general result of transverse nonlinear instability of 1d standing waves for Hamiltonian PDE for both periodic and localized transverse perturbations.
In \cite{Yamazaki3}, Yamazaki investigated the transverse instability of the 
line standing waves under the flow of the two-dimensional nonlinear Schr\"{o}dinger equation:
\begin{equation} \label{NLS}
i \partial_{t} \psi + \Delta_{x, y} \psi + |\psi|^{p-1}\psi = 0 \qquad 
\mbox{in $\R \times \R \times \T$}. 
\end{equation}
The result of \cite{Yamazaki3} states that 
there exists $\omega_{c}>0$ such that 
the line standing wave $e^{i \omega t}R_{\omega}$ is 
stable for $0 < \omega < \omega_{c}$ and 
unstable for $\omega > \omega_{c}$ under some assumption on $p$. 
Moreover in \cite{Yamazaki2}, Yamazaki also studied 
the case of the critical frequency $\omega_{c}$ 
and showed that 
there exists $2 < p_{1} < p_{2} < 3$ such that 
the standing waves is stable for $2 \leq p < p_{1}$ 
and unstable for $p_{2} < p < 5$. 

In this paper, following Yamazaki~\cite{Yamazaki3}, we will discuss the stability of the line standing wave $e^{i \omega t}R_{\omega}$ with respect to the frequency $\omega >0$. 
Our first result provides a classification for transverse stability and instability for all frequencies excluding the critical one.
\begin{theorem} 
\label{thm-1}
Let $1<p<5$, $\omega_{p} 
= \frac{4}{(p-1)(p+3)}$ and $s>\frac{1}{2}$. 
Then, in the space $H^{1}_{x}L^{2}_{y} \cap L^{2}_{x}H^{s}_{y}(\mathbb{R} \times \mathbb{T})$, the following holds:
\begin{itemize}
\item[(i)] for $0 < \omega < \omega_{p}$, the standing wave $e^{i \omega t}R_{\omega}$ 
is stable under the flow of \eqref{WS}, 
\item[(ii)] for $\omega > \omega_{p}$, the standing wave $e^{i \omega t}R_{\omega}$ 
is unstable under the flow of \eqref{WS}. 
\end{itemize}
\end{theorem}
\begin{remark}
The critical frequency case is more challenging because the linearized operator is degenerate. We will undertake this problem in a forthcoming paper.
%{\color{red} 
%We refer to our forthcoming paper \cite{BIK1} 
%for full investigation of this case.}
\end{remark}
The proof of Theorem \ref{thm-1} {\rm (i)} relies on
the argument \`a la Grillakis, Shatah 
and Strauss~\cite{GSS}. 
We obtain Theorem \ref{thm-1} by showing some positivity 
of the linearized operator around 
the line ground state when $\omega < \omega_{p}$. 
For the instability case, we start by 
proving the linear instability, using the implicit function theorem
i.e. the linearized Hamiltonian has a positive eigenvalue. 
Then, the nonlinear instability follows using the argument of Yamazaki~\cite{Yamazaki1}.

Next, we give a characterization of the ground state to \eqref{sp}. 
We note that if we define a functional $\mathcal{S}_{\omega}$ on $X$ 
by 
\begin{equation} \label{action}
\mathcal{S}_{\omega} (u) := 
\frac{1}{2} \|\nabla u\|_{L_{x, y}^{2}(\R \times 
\T)}^{2} 
+ \frac{1}{2} \||D_{y}|^{\frac{1}{2}} u\|_{L_{x, y}^{2}(\R \times \T)}^{2}
+ \frac{\omega}{2} \|u\|_{L_{x, y}^{2}(\R \times \T)}^{2} 
- \frac{1}{p+1} \|u\|_{L_{x, y}^{p+1}(\R \times \T)}^{p+1}, 
\end{equation}
we see that $u$ is a solution to the equation \eqref{sp} 
if and only if $u$ is a critical point of 
$\mathcal{S}_{\omega
}$. 
Then, by a ground state, we mean the least action solution among all nontrivial solutions to \eqref{sp} in $X$. 
We seek a ground state to \eqref{sp} by considering the following minimization 
problem: 
\begin{equation} \label{mini-pro}
m_{\omega} := \inf \left\{\mathcal{S}_{\omega}(u) \colon 
u \in X \setminus \{0\}, \; \mathcal{N}_{\omega}(u) = 0 \right\}, 
\end{equation} 
where 
\begin{equation} \label{nehari}
\begin{split}
\mathcal{N}_{\omega} (u) 
& = \langle \mathcal{S}_{\omega}^{\prime}(u), u \rangle \\
& = 
\|\nabla u\|_{L_{x, y}^{2}(\R \times \T)}^{2} 
+ \||D_{y}|^{\frac{1}{2}} u\|_{L_{x, y}^{2}(\R \times \T)}^{2}
+ \omega \|u\|_{L_{x, y}^{2}(\R \times \T)}^{2} 
- \|u\|_{L_{x, y}^{p+1}(\R \times \T)}^{p+1}. 
\end{split}
\end{equation}
Then, by a standard argument, we shall show the following: 
\begin{theorem} \label{thm-exist-ground}
Let $1 < p< 5$ and $\omega >0$. 
For any $\omega >0$, there exist ground states $Q_{\omega} \neq 0$ 
satisfying $\mathcal{S}_{\omega}(Q_{\omega}) = m_{\omega}$. 
\end{theorem}
For the sake of completeness, we provide the 
proof of Theorem \ref{thm-exist-ground}. 
Here, we shall investigate 
whether the ground states $R_{\omega}$ 
to \eqref{line-sp} are also ground states to \eqref{sp} or not. 
Following Berestycki and Wei~\cite{BW} and 
Terracini,Tzvetkov and Visciglia~\cite{TTV}, 
we show that the line soliton $R_{\omega}$ coincides with
$Q_{\omega}$ for low frequencies. 
This implies the uniqueness of ground state to \eqref{sp}. 
On the other hand, we show that for high frequencies, 
the line soliton $R_{\omega}$ is different from $Q_{\omega}$. 
The results are given by the following:
\begin{theorem} \label{thm-chara-ground}
Let $1 < p < 5$, $R_{\omega}$ be the line soliton of \eqref{sp} 
and $Q_{\omega}$ be the ground state of \eqref{sp}. 
Then, there exists $\omega_{*} \in (0, \omega_{p}]$ such that 
$Q_{\omega} = R_{\omega}$ for all $0 < \omega \leq \omega_{*}$ and 
$Q_{\omega} \neq R_{\omega}$ for all $\omega > \omega_{*}$.
\end{theorem}
Note that 
\begin{equation} \label{rescale-R}
R_{\omega}(x) = \omega^{ \frac{1}{p-1}} 
R_1(\omega^{\frac{1}{2}}x) \quad {\rm for} \ x \in \R.
\end{equation}
Furthermore, 
we put 
\begin{equation} \label{rescale-Q}
Q_{\omega}(x, y) = \omega^{\frac{1}{p-1}} 
\widetilde{Q}_{\omega}(\omega^{\frac{1}{2}}x, y) 
\quad {\rm for} \ (x, y) \in \R \times \T. 
\end{equation}
From \eqref{rescale-R} and 
\eqref{rescale-Q}, it is sufficient to investigate whether 
$\widetilde{Q}_{\omega} = R_{1}$ or not. 
\par
To this end, we first show that 
the equality does not hold 
for sufficiently large $\omega > 0$. 
We can prove this 
by showing that 
$\widetilde{\mathcal{S}}_{\omega}
(\widetilde{Q}_{\omega}) < 
\widetilde{\mathcal{S}}_{\omega}(R_{1})$ for 
sufficiently large $\omega>0$, 
where $\widetilde{\mathcal{S}}_{\omega}$ 
is the action functional corresponding to the equation 
\eqref{msp} (see \eqref{rescale-action} below
for the precise definition of 
$\widetilde{\mathcal{S}}_{\omega}$). 
\par
Secondly, we will show that 
$\widetilde{Q}_{\omega} = R_{1}$ if 
$\omega > 0$ is sufficiently small. 
One can verify that for any $\omega>0$, 
$\widetilde{Q}_{\omega}$ satisfies 
\begin{equation} \label{msp}
- \partial_{xx} \widetilde{Q}_{\omega} + \omega^{-1}
|D_{y}| \widetilde{Q}_{\omega} + \widetilde{Q}_{\omega} - 
|\widetilde{Q}_{\omega}|^{p-1}\widetilde{Q}_{\omega} = 0 \qquad 
\mbox{in $\mathbb{R} \times \mathbb{T}$}. 
\end{equation} 
Then, we find that $\widetilde{Q}_{\omega}$ 
converges to $R_1$, up to translation, 
in the energy space $X$ as $\omega\to 0$. 
Furthermore, using this together with 
the equation \eqref{msp} and an energy estimate, 
we deduce that $|D_y|\widetilde{Q}_{\omega}=0$. 
Namely, $\widetilde{Q}_{\omega}$ does not depend 
on $y$, so that it satisfies the equation \eqref{line-sp} 
with $\omega = 1$. 
From the uniqueness of the solution to \eqref{line-sp}, 
we obtain $\widetilde{Q}_{\omega} = R_{1}$. 
\par 
Once we find that $Q_{\omega}$ coincides with $R_{\omega}$ for 
sufficiently small $\omega > 0$ and not for 
sufficiently large $\omega > 0$, 
we would like to investigate the existence of a threshold frequency $\omega_*$ separating the above two properties. Define
%we study a general frequency $\omega > 0$. Putting 
\[
\omega_{*} := \left\{\omega > 0 \colon 
\mathcal{S}_{\omega} (Q_{\omega}) 
= \mathcal{S}_{\omega} (R_{\omega}) \right\}. 
\]
From the above arguments, we clearly see that $0 < \omega_{*} < \infty$. 
Using the continuity of the function $\omega \mapsto m_\omega$, we will find that for $0 < \omega \leq \omega_{*}$, 
$Q_{\omega} = R_{\omega}$ and for $\omega >\omega_{*}$, 
$Q_{\omega} \neq R_{\omega}$. 
In addition, by analyzing the second eigenvalue 
of the linearized operator around the ground state, 
we shall show $\omega_{*} \leq \omega_{p}$.

\color{black}
The rest of this paper is organized as follows. 
In Section 2, we study the stability of $R_\omega$ 
when $\omega<\omega_p$ using the result of 
Grillakis Shatah and Strauss~\cite{GSS}. In Section 3, we prove the instability result in case of $\omega > \omega_{p}$. 
We first show the linear instability using the implicit function 
theorem and find that the linear instability implies the 
non-linear instability from a classical argument.
% In Section 4, we treat the critical frequency $\omega = \omega_{p}$. 
%We show that both stability and instability hold depending 
%on the value of the exponent $p \in (2,5)$. 
In Section 4, we provide a characterization of the ground state of \eqref{WS}. 
In Section 5, we obtain a regularity result of the ground state of \eqref{WS}. In Appendix, we give some technical tools for the sake of completeness. 
%%%%%%%%%%%%%%%%%%%%%%%%%%%%%%%%%%%%%%%%%%%%%%%%%%%%%%%%%

\subsection{Notation}
\noindent 
\textbullet~For a function $u$, denote by $\Re u$ and $\Im u$ by
the real part of $u$ and the imaginary 
part of $u$, respectively.

\noindent 
\textbullet~
For $f, g \in L_{x,y}^{2}(\mathbb{R} \times \mathbb{T})$, 
We define the inner product by 
\begin{equation*}
(f,g)_{L_{x, y}^{2}(\mathbb{R} \times \mathbb{T})} 
= \Re 
\int_{\mathbb R\times\mathbb T}
f(x,y)\overline{g(x,y)}\; dxdy. 
\end{equation*}
For $u, v \in L_{x}^{2}(\mathbb R)$, 
We define the inner product by 
\begin{equation*}
(u, v)_{L_{x}^{2}(\mathbb{R})} 
= \Re 
\int_{\mathbb R}u(x)\overline{v(x)}\;dx.
\end{equation*}

\noindent 
\textbullet Let $X^{*}$ be the dual space of 
the function space $X$. 
We use the convention that if 
$f \in X^{*}$ and $g \in X$, then 
$\langle f, g \rangle_{X^{*}, X}$ denotes the duality pairing of $X^{*}$ and $X$, namely 
\begin{equation*}%\label{17/08/15/06:03}
\langle f, g \rangle_{X^{*}, X}
=
\Re \int_{\mathbb{R} \times \mathbb{T}} 
(- \partial_{xx} + |D_{y}| + 1)^{-\frac{1}{2}}f(x, y) 
\overline{(-\partial_{xx} + |D_{y}| + 1)^{\frac{1}{2}} 
g(x, y)}\,dx.
\end{equation*}
Similarly, 
we denote the duality pairing of 
$H_{x}^{-1}(\mathbb{R})$ and $H_{x}^{1}(\mathbb{R})$ 
by $\langle u, v \rangle_{H_{x}^{-1}(\mathbb{R}), H_{x}^{1}(\mathbb{R})}$
for $u \in H_{x}^{-1}(\mathbb{R})$ and 
$v \in H_{x}^{1}(\mathbb{R})$, namely 
\begin{equation*}
\langle u, v \rangle_{H_{x}^{-1}(\mathbb{R}), H_{x}^{1}(\mathbb{R})}
=
\Re \int_{\mathbb{R}} 
(- \partial_{xx} + 1)^{-\frac{1}{2}}u(x) 
\overline{(-\partial_{xx} + 1)^{\frac{1}{2}} 
v(x)}\,dx.
\end{equation*}

\noindent 
\textbullet
For a Banach space $Y$, we denote by $\mathcal{L}(Y)$ 
a set of the bounded linear operators from $Y$ to itself 
and by $\|\cdot \|_{\mathcal{L}(Y)}$ the operator norm. 

\noindent 
\textbullet~For given positive quantities $a$ and $b$, the notation $a \lesssim b$ means the inequality $a\le C b$ for some positive constant; the constant $C$ depends only on $d$ and $p$, unless otherwise noted.

%%%%%%%%%%%%%%%%%%%%%%%%%%%%%%%%%%%%%%%%%%%%%%%%%%%%%%%%%

\section{Stability in case of $0<\omega<\omega_p$}
To prove Theorem \ref{thm-1}, we use a well-known sufficient condition for the 
stability due to Grillakis, Shatah and Strauss~\cite{GSS}. 
\begin{proposition}[\cite{GSS}] \label{prop-2}
Assume that there exists a constant $\delta >0$ 
such that 
\begin{equation*}
\langle \mathcal{S}_{\omega}^{\prime \prime}(R_{\omega})u, u \rangle_{X^{*}, X} \geq \delta 
\|u\|_{X}^{2}
\end{equation*}
for all $u \in X$ satisfying 
\begin{equation} \label{orthogonal-cond1}
(u, R_{\omega})_{L_{x, y}^{2}(\mathbb{R} \times \mathbb{T})} = 
(u, i R_{\omega})_{L_{x, y}^{2}(\mathbb{R} \times \mathbb{T})} = 
(u, \partial_{x} R_{\omega})_{L_{x, y}^{2}(\mathbb{R} \times \mathbb{T})} = 0.
\end{equation} 
Then, the standing wave $e^{i \omega t}R_{\omega}$ is stable. 
\end{proposition}
Note that, the Fr\'{e}chet derivative $\mathcal{S}_{\omega}^{\prime \prime}(R_{\omega}) $ satisfies 
\begin{equation*}
\langle \mathcal{S}_{\omega}^{\prime \prime}(R_{\omega})u, u \rangle_{X^{*}, X} 
= \langle L_{\omega, +}
\Re u, \Re u \rangle_{X^{*}, X} 
+ \langle L_{\omega, -} 
\Im u, \Im u \rangle_{X^{*}, X}, 
\end{equation*}
where the operators $L_{\omega, +}$ and $L_{\omega, -}$ are defined by 
\begin{equation} \label{1d-ope}
L_{\omega, +} = - \partial_{xx} + \omega + |D_{y}| - p R_{\omega}^{p-1}, \qquad 
L_{\omega, -} = - \partial_{xx} + \omega + |D_{y}| - R_{\omega}^{p-1}. 
\end{equation}
Moreover, for each $u \in X$, we put 
\begin{equation*}
\Re u = \sum_{n \in \mathbb{Z}} u_{n}^{R} e^{i n y}, \qquad 
\Im u = \sum_{n \in \mathbb{Z}} u_{n}^{I} e^{i n y},
\end{equation*}
with the reality condition $\overline{u_{n}^{R}} 
=u_{-n}^{R}$ and $\overline{u_{n}^{I}}=u_{-n}^{I}.$ This yields that 
\begin{equation*}
- \partial_{xx} + |D_{y}| = \bigoplus_{n \in \mathbb{Z}}\left( 
- \partial_{xx} + |n| \right). 
\end{equation*}
Therefore, if we denote 
\begin{equation} \label{2d-ope}
L_{\omega, +, n} = - \partial_{xx} + \omega + |n| - p R_{\omega}^{p-1}, \qquad 
L_{\omega, -, n} = - \partial_{xx} + \omega + |n| - R_{\omega}^{p-1}. 
\end{equation}
we can write
\begin{equation*}
\langle \mathcal{S}_{\omega}^{\prime \prime}(R_{\omega})u, u \rangle_{X^{*}, X} 
= \sum_{n \in \mathbb{Z}}\left( 
\langle L_{\omega, +, n} u_{n}^{R}, u_{n}^{R} 
\rangle_{H_{x}^{-1}(\mathbb{R}), H_{x}^{1}(\mathbb{R})} + 
\langle L_{\omega, -, n} u_{n}^{I}, u_{n}^{I} 
\rangle_{H_{x}^{-1}(\mathbb{R}), H_{x}^{1}(\mathbb{R})} \right). 
\end{equation*}
From the Weyl's essential theorem, we see that 
$\sigma_{\text{ess}}(L_{\omega, +, 0}), \sigma_{\text{ess}}(L_{\omega, -, 0})\subset [\omega, \infty)$. 
Concerning the discrete spectrum of 
the operators $L_{\omega, +, 0}$ and $L_{\omega, -, 0}$, 
we obtain the following: 
\begin{lemma}
\label{lem:neg-eigenvalue}
Assume that $1 < p < \infty$. 
The operator $L_{\omega, +, 0}$ has only one negative eigenvalue $-\frac{\omega}{\omega_p}$ with corresponding eigenfunction $R_{\omega}^{\frac{p+1}{2}}$ and $0$ eigenvalue with corresponding eigenfunction $\partial_{x} R_{\omega}$. 
\end{lemma}
See e.g. \cite[Theorems 3.1 and 3.4]{Chang-Gustaf-Nakanishi-Tsai} 
for the proof of Lemma \ref{lem:neg-eigenvalue}. 
From Lemma \ref{lem:neg-eigenvalue}, we have 
\begin{equation}\label{lowerbound-L-}
\langle L_{\omega, +, 0} v, v 
\rangle_{H_{x}^{-1}(\mathbb{R}), H_{x}^{1}(\mathbb{R})} 
\geq - \frac{\omega}{\omega_{p}} \|v\|_{L_{x}^{2}(\mathbb{R})}^{2}. 
\end{equation}

Since $L_{\omega, -, 0}R_{\omega} = 0$ and $R_{\omega}(x) > 0$ 
for all $x \in \mathbb{R}$, we obtain the following: 
\begin{lemma} \label{spectrum-operator}
Assume that $1 < p < \infty$. 
the operator $L_{\omega, -, 0}$ is non-negative and the zero eigenvalue 
is simple with corresponding eigenfunction $R_{\omega}$. 
\end{lemma}
From Lemma \ref{spectrum-operator}, 
we see that 
there exists a constant $C_{1} > 0$ such that 
\begin{equation} \label{positive-}
\langle L_{\omega, -, 0} v, v \rangle_{H_{x}^{-1}(\R), 
H_{x}^{1}(\R)} 
\geq C_{1} \|v\|_{L_{x}^{2}(\R)}^{2} 
\end{equation}
for all $v \in H_{x}^{1}(\mathbb{R}, \mathbb{R})$ satisfying 
$(v, R_{\omega})_{L_{x}^{2}(\R)} = 0$. 
Moreover, applying the abstract theorem of 
Grillakis, Shatah and Strauss~\cite{GSS}, 
we obtain the following: 
\begin{proposition} \label{sufficient-ope}
Assume that $\omega >0$ and $1 < p < 5$. 
There exists a constant $C_{2} > 0$ such that 
\begin{equation} \label{positive+}
\langle L_{\omega, +, 0} v, v \rangle_{H_{x}^{-1}
(\R), H_{x}^{1}(\R)} 
\geq C_{2} \|v\|_{L_{x}^{2}(\R)}^{2}
\end{equation} 
for all $v \in H^{1}_{x}(\mathbb{R}, \mathbb{R})$ satisfying 
$(v, R_{\omega})_{L^{2}_{x}(\R)} 
= (v, \partial_{x} R_{\omega})_{L^{2}_{x}(\R)} = 0$. 
\end{proposition}
\begin{remark}
We note that when $p > 5$, 
the result of Proposition \ref{sufficient-ope} does not hold 
(see e.g. Grillakis, Shatah and Strauss~\cite[Theorem 4.1]{GSS}). 
\end{remark}

We are now in a position to prove the first part of Theorem \ref{thm-1}. 
\begin{proof}[Proof of Theorem \ref{thm-1} {\rm (i)}]
%Let $u \in X$. 
Note that 
\begin{equation*}
L_{\omega, +, n} = L_{\omega, +, 0} + |n|, \qquad 
L_{\omega, -, n} = L_{\omega, -, 0} + |n|. 
\end{equation*}
This together with \eqref{lowerbound-L-} and \eqref{positive-} 
yields 
that for $0 < \omega < \omega_p$, we obtain 
\begin{equation*} %\label{sta-lowerbound-1}
\langle L_{\omega, -, n} u_{n}^{I}, u_{n}^{I} 
\rangle_{H_{x}^{-1}(\mathbb{R}), H_{x}^{1}(\mathbb{R})} 
= \langle L_{\omega, -, 0} u_{n}^{I}, u_{n}^{I} 
\rangle_{H_{x}^{-1}(\mathbb{R}), H_{x}^{1}(\mathbb{R})}
+ |n| \|u_{n}^{I}\|_{L_{x}^{2}(\R)}^{2} 
\geq 
|n|\|u_{n}^{I}\|_{L_{x}^{2}(\R)}^{2}, 
\end{equation*}
\begin{equation}
\begin{split}
\langle L_{\omega, +, n} u_{n}^{R}, u_{n}^{R} \rangle_{H_{x}^{-1}(\mathbb{R}), H_{x}^{1}(\mathbb{R})} 
& = \langle L_{\omega, +, 0} u_{n}^{R}, u_{n}^{R} \rangle_{H_{x}^{-1}(\mathbb{R}), H_{x}^{1}(\mathbb{R})} 
+ |n|\|u_{n}\|_{L_{x}^{2}(\R)}^{2} \\ 
& \geq (1 - \frac{\omega}{\omega_{p}})\|u_{n}^{R}\|_{L_{x}^{2}(\R)}^{2} 
\label{sta-lowerbound-3}
\end{split}
\end{equation}
for all $n\in \mathbb{Z} \setminus \{0\}$. 
Moreover, one has 
\begin{equation*}
\begin{split}
& (u, R_{\omega})_{L_{x, y}^{2}(\mathbb{R} \times \mathbb{T})} 
= 2 \pi (u_{0}^{R}, R_{\omega})_{L_{x}^{2}(\mathbb{R})}, \qquad 
(u, iR_{\omega})_{L_{x}^{2}(\mathbb{R} \times \mathbb{T})} 
= 2\pi (u_{0}^{I}, R_{\omega})_{L_{x}^{2}(\mathbb{R})}, \qquad \\
&(u, \partial_{x} R_{\omega})_{L_{x, y}^{2}(\mathbb{R} \times \mathbb{T})} 
= 2 \pi (u_{0}^{R}, 
\partial_{x} R_{\omega})_{L_{x}^{2}(\mathbb{R})}. 
\end{split}
\end{equation*}
Therefore, we see from \eqref{orthogonal-cond1} that 
\begin{equation*}
(u_{0}^{R}, R_{\omega})_{L_{x}^{2}(\mathbb{R})} 
= (u_{0}^{I}, R_{\omega})_{L_{x}^{2}(\mathbb{R})}
= (u_{0}^{R}, \partial_{x} R_{\omega})_{L_{x}^{2}(\mathbb{R})} = 0. 
\end{equation*}
This together with 
\eqref{positive-} and 
\eqref{positive+} yields that 
\begin{equation} \label{sta-lowerbound-2}
\langle L_{\omega, -, 0} u_{0}^{I}, u_{0}^{I} 
\rangle_{H_{x}^{-1}(\mathbb{R}), H_{x}^{1}(\mathbb{R})} 
\geq C_{1} \|u_{0}^{I}\|_{L_{x}^{2}(\R)}^{2}, \quad
\langle L_{\omega, +, 0} u_{0}^{R}, u_{0}^{R} 
\rangle_{H_{x}^{-1}(\mathbb{R}), H_{x}^{1}(\mathbb{R})} 
\geq C_{2} \|u_{0}^{R}\|_{L_{x}^{2}(\R)}^{2}. 
\end{equation}
It follows from \eqref{sta-lowerbound-3} 
and \eqref{sta-lowerbound-2} that 
when $0 < \omega < \omega_{p}$, 
we have 
\begin{equation*}
\langle L_{\omega, +} \Re u, \Re u 
\rangle_{X^{*}, X} 
\geq C_2 \|\Re u\|_{L_{x, y}^{2}(\R\times\T)}^{2} \geq p \delta \int_{\mathbb{R} \times \mathbb{T}} R_{\omega}^{p-1} |\Re u|^{2} dxdy. 
\end{equation*}
for sufficiently small $\delta>0$. We infer that 
\begin{equation*}
\begin{split}
0 
& \leq \langle L_{\omega, +} \Re u, \Re u 
\rangle_{X^{*}, X} 
- p \delta \int_{\mathbb{R} \times \mathbb{T}} R_{\omega}^{p-1} |\Re u|^{2} dxdy \\
%& \leq \int_{\mathbb{R} \times \mathbb{T}} 
%\left(|\partial_{x} \Re u(x, y)|^{2} + \omega |\Re u(x, y)|^{2} 
%+ |D_{y}| \Re u(x, y) \overline{\Re u(x, y)}\right)dxdy \\
%& - p (1 + \delta)\int_{\mathbb{R} \times \mathbb{T}}|\Re u(x, y)|^{p+1} dxdy \\
& = (1+\delta) 
\langle L_{\omega, +} \Re u, \Re u 
\rangle_{X^{*}, X} \\
& - \delta \int_{\mathbb{R} \times \mathbb{T}} \left(|\partial_{x} \Re u(x, y)|^{2} + \omega |\Re u(x, y)|^{2} 
+ |D_{y}| \Re u(x, y) \overline{\Re u(x, y)}\right)dxdy. 
\end{split}
\end{equation*}
This implies that 
\begin{equation*}
\begin{split}
& \langle L_{\omega, +} \Re u, 
\Re u \rangle_{X^{*}, X} \\
& \geq \frac{\delta}{1+ \delta} 
\int_{\mathbb{R} \times \mathbb{T}} 
\Big(|\partial_{x} %& 
\Re u(x, y)|^{2} + \omega |\Re u(x, y)|^{2} 
%\\ &
+ |D_{y}| \Re u(x, y) \overline{\Re u(x, y)}\Big)dxdy. 
\end{split}
\end{equation*}
Similarly, we can obtain 
\begin{equation*}
\langle L_{\omega, -} \Im u, \Im u \rangle_{X^{*}, X} 
\geq \delta^{\prime}
\int_{\mathbb{R} \times \mathbb{T}} 
\left(|\partial_{x} \Im u(x, y)|^{2} + \omega |\Im u(x, y)|^{2} 
+ ||D_{y}|^{\frac{1}{2}}\Im u(x, y)|^{2} \right)dxdy. 
\end{equation*}
Therefore, we see that 
for $0 < \omega < \omega_{p}$, 
there exists a constant $C>0$ such that 
\begin{equation*}
\langle \mathcal{S}_{\omega}^{\prime \prime}(R_{\omega})u, u \rangle_{X^{*}, X} \geq C 
\|u\|^{2}_{X}
\end{equation*}
satisfying 
$(u, R_{\omega})_{L_{x, y}^{2}(\mathbb{R} \times \mathbb{T})} = 
(u, i R_{\omega})_{L_{x, y}^{2}(\mathbb{R} \times \mathbb{T})} = 
(u, \partial_{x} R_{\omega})_{L_{x, y}^{2}(\mathbb{R} \times \mathbb{T})} = 0$.
This completes the proof. 
\end{proof}

\section{Instability in case of $\omega > \omega_{p}$}
This section is devoted to the proof of Theorem \ref{thm-1} {\rm (ii)}. 
In the sequel, we assume that $\omega>\omega_p$. 
We identify $z \in \mathbb{C}$ with ${}^{t} (\Re z, \Im z) 
\in \mathbb{R}^{2}$ and put $z = {}^{t} (\Re z, \Im z)$. 
Putting $f(z) = |z|^{p-1}z$ for $z \in \mathbb{R}^{2}$, 
we see that 
$f \in C^{1}(\mathbb{R}^{2}, \mathbb{R}^{2})$ and 
\begin{equation} \label{eq-deri}
D f(z) = 
\begin{pmatrix}
(p-1) |z|^{p-3} (\Re z)^{2} + |z|^{p-1} 
& (p-1) |z|^{p-3} \Re z \Im z \\
(p-1) |z|^{p-3} \Re z \Im z 
& (p-1) |z|^{p-3} (\Im z)^{2} + |z|^{p-1} 
\end{pmatrix}, 
\end{equation}
where $D f(z)$ is the derivative of $f$ at $z \in \mathbb{R}^{2}$. 

Note that the Hamiltonian form of \eqref{WS} is given by
\begin{equation}
\label{WS-Hamil}
u_t=- J \mathcal{H}^{\prime}(u),
\end{equation} 
where $\mathcal{H}^{\prime}$ is the Fr\'{e}chet 
derivative of $\mathcal{H}$ and 
\begin{equation} \label{complex}
J:=\begin{pmatrix}
0 &-1\\
1 &0
\end{pmatrix}.
\end{equation}
We write
\begin{equation} \label{eq-error}
u(t) =e^{i\omega t}\left(R_{\omega} + v(t)\right).
\end{equation} 
Then, we can easily verify that $v:={}^{t}(\Re v, \Im v)$ is a solution of\footnote{Abuse of notation $R_\omega:={}^t(R_\omega,0)$.}
\begin{equation}
\label{WS-Ham-v}
v_t= -J (\mathcal{S}_{\omega}^{\prime \prime}(R_{\omega}) v + NL(v,R_\omega)) 
\end{equation}
with 
\begin{equation} \label{linearized-nonlinear}
\begin{split}
NL(v,R_\omega) 
& = \mathcal{H}^{\prime}(R_{\omega} + v) -
\mathcal{H}^{\prime}(R_{\omega}) - 
\mathcal{H}^{\prime \prime}(R_{\omega}) v \\
& = f(R_{\omega} + v) - f(R_{\omega}) - D f(R_{\omega}) v.
\end{split}
\end{equation}
It follows from \eqref{eq-deri} that 
\begin{equation} \label{eq-deri2}
NL(v,R_\omega) 
= \begin{pmatrix}
|v+R_{\omega}|^{p-1}(\Re v+R_{\omega}) -p|R_{\omega}|^{p-1}\Re v -|R_{\omega}|^{p-1}R_{\omega}\\
|v+R_{\omega}|^{p-1}\Im v -|R_{\omega}|^{p-1}\Im v
\end{pmatrix}.
\end{equation}

We shall show that for $\omega>\omega_p$, 
$-J \mathcal{S}_{\omega}^{\prime \prime}(R_{\omega})$ has at least one positive eigenvalue. 
Since 
\begin{equation} \label{eq-linearizedop}
-J \mathcal{S}_{\omega}^{\prime \prime}(R_{\omega}) 
= \bigoplus_{n \in \mathbb{Z}} -J S_{\omega}(n),
\end{equation}
where, for $a \in \R$,
\begin{equation} \label{linerized-ope-1}
S_{\omega}(a):= \begin{pmatrix}
L_{\omega, +, a}& 0 \\
0 & L_{\omega, -, a} 
\end{pmatrix}.
\end{equation}
We know that $-J \mathcal{S}_{\omega}^{\prime \prime}(R_{\omega})$ 
has an eigenvalue if and only if there exists $n \in \mathbb{Z}$ such that 
$- J S_\omega(n)$ has the same eigenvalue. 
Therefore, 
in order to obtain the linear instability of the line soliton $R_{\omega}$, 
it is sufficient to show that 
$-J S_\omega(1)$ has a positive eigenvalue. 
We obtain the following:
\begin{lemma}
\label{lem:posi-eigen}
For $\omega>\omega_p$, 
$-J S_{\omega}(1)$ has at least one positive eigenvalue. 
\end{lemma}
\begin{proof}
The proof follows Yamazaki's steps~\cite[Proposition 1]{Yamazaki3} who employed the 
argument of Rousset and Tzvetkov~\cite[Theorem 1]{R-T-1} and 
Rousset and Tzvetkov~\cite[Lemma 3.1]{R-T-2}.

Let 
\begin{equation*}
\phi_{\omega} := 
\begin{pmatrix}
R^{\frac{p+1}{2}}_{\omega} \\
0
\end{pmatrix}
\end{equation*} 
and
\begin{equation*}
U(v, a, \lambda):= S_{\omega}(a) 
(\phi_{\omega}+v) + \lambda J^{-1} (\phi_{\omega}+v),
\end{equation*}
with $v \in {\rm Span}(\phi_{\omega})^{\perp}$ 
and $\omega, \ \lambda >0$, 
where 
\begin{equation*}
{\rm Span}(\phi_{\omega})^{\perp} = 
\left\{v \in L_{x}^{2}(\mathbb{R}) \times L_{x}^{2}(\mathbb{R}) 
\colon \langle v, \phi_{\omega} \rangle = 0 \right\}. 
\end{equation*}
From Lemmas \ref{lem:neg-eigenvalue} and \ref{spectrum-operator}, 
we know that the kernel of $S_{\omega} (\nu_{\omega})$ 
is spanned by $\phi_{\omega}$, where 
$\nu_{\omega} = \omega/\omega_{p}$. 
Note that $\nu_{\omega} > 1$ for $\omega > \omega_{p}$. 
We see that $U$ is a $C^\infty$ 
map from 
$(\phi_{\omega})^{\perp} \times \mathbb{R} \times \mathbb{R}$ 
to $L_{x}^2(\mathbb{R})\times L_{x}^2(\mathbb{R})$ 
%\textcolor{red}{IT SEEMS TO ME THAT $U$ is %UNBOUNDED OPERATOR if the %target space is $L^2$, %Any clarification?} 
and 
\begin{equation*}
U(0, \nu_{\omega}, 0)=0. 
\end{equation*}
Differentiate $U(v, a, \lambda)$ with respect to 
$a, v$ and substituting $v = 0, a = \nu_{\omega}$ and $\lambda = 0$, 
we have 
\begin{equation*}
(D_{a}U(0, \nu_{\omega}, 0), D_{v}U(0, \nu_{\omega}, 0)) 
(\mu, w)
= (\mu \phi_{\omega}, S_{\omega} (\nu_{\omega}) w).
\end{equation*}
for any $\mu \in \mathbb{R}$ and 
$w \in L_{x}^2(\mathbb{R})\times L_{x}^2(\mathbb{R})$. 
Thus, $D_{a} U(0, \nu_{\omega}, 0) \in \mathcal{L}(\mathbb{R}, 
L_{x}^2(\mathbb{R})\times L_{x}^2(\mathbb{R}))$ 
%\textcolor{red}{It seems to me that the target space %should be $L^2$ and not %$L^2\times L^2$. Same remark %in below} 
is invertible 
since it is a linear one-to-one mapping from $\mathbb{R}$ to ${\rm Span}(\phi_{\omega})$. 
In addition, from the fact that the kernel of 
$S_{\omega} (\nu_{\omega})$ is spanned by $\phi_{\omega}$, 
we infer that $D_{v}U(0, \nu_{\omega}, 0) \in 
\mathcal{L}((\phi_{\omega})^{\perp}, 
L_{x}^2(\mathbb{R}) \times L_{x}^2(\mathbb{R}))$ is invertible. 
Then, by the implicit function theorem, there exist two $C^\infty$ functions 
$a(\lambda) \in \mathbb{R}$ and 
$v(\lambda) \in (\phi_{\omega})^{\perp}$ 
satisfying $a(0)= \nu_{\omega}$, $v(0)=0$ and 
$U(v(\lambda), a(\lambda), \lambda)=0$, that is, 
\begin{equation*}
0 = U(v(\lambda), a(\lambda), \lambda) = 
S_{\omega} (a(\lambda)) (\phi_{\omega} 
+v(\lambda)) + \lambda J^{-1} (\phi_{\omega}+v(\lambda)),
\end{equation*}
for sufficiently small $\lambda$. 
Differentiating the previous identity with respect to $\lambda$, 
we have
\begin{equation} \label{identity:U'}
\begin{split}
0 = \frac{d }{d \lambda} U(v(\lambda), a(\lambda), \lambda)
& = 
a^{\prime}(\lambda)(\phi_{\omega} + v(\lambda)) 
+S_{\omega}(a(\lambda)) v^{\prime}(\lambda) \\
& \quad + J^{-1}(\phi_{\omega}+v(\lambda))
+\lambda J^{-1} v^{\prime}(\lambda).
\end{split}
\end{equation}
Since $v(0) = 0$ and $a(0) = \nu_{\omega}$, 
for $\lambda=0$, we obtain
\begin{equation}
\label{identity:lambda=0}
0 = a^{\prime}(0) \phi_{\omega} + 
S_{\omega}(\nu_{\omega})v^{\prime}(0) + J^{-1}\phi_{\omega}.
\end{equation}
Taking the $L^2$-scalar product of identity \eqref{identity:lambda=0} 
with $\phi_{\omega}$ and using the fact that $\phi_{\omega} \in 
\mathop{\mathrm{Ker}} S_{\omega} (\nu_{\omega})$ and 
$J$ is a skew-symmetric operator, we get
\begin{equation*}
0 = a^{\prime}(0) \|\phi_{\omega}\|_{L_{x}^2(\R) \times L_{x}^{2}(\R)}^{2}.
\end{equation*}
Thus, $a^{\prime}(0)=0$. 
Substituting $a^{\prime}(0) = 0$ into \eqref{identity:lambda=0}, we see that
\begin{equation}
\label{identity;v'(0)}
S_{\omega}(\nu_{\omega}) v^{\prime}(0) =- J^{-1}\phi_{\omega}.
\end{equation}
Differentiating \eqref{identity:U'} with respect to $\lambda$ again, 
we have 
\begin{equation*}
\begin{split}
0 = \frac{d^2 U}{d \lambda^2} (v(\lambda),\omega(\lambda),\lambda) 
= & a^{\prime \prime}(\lambda)(\phi_{\omega} + v(\lambda)) 
+ 2 a^{\prime}(\lambda)v^{\prime}(\lambda) 
+ S_{\omega}(a(\lambda))v^{\prime \prime}(\lambda) \\
& + 2 J^{-1}v^{\prime}(\lambda) + \lambda J^{-1}v^{\prime \prime}(\lambda). 
%\label{identity:U^{\prime \prime}}
\end{split}
\end{equation*}
Since $v(0) = 0$, $a(0) = \nu_{\omega}$ and 
$a^{\prime}(0) = 0$, putting $\lambda=0$, we have
\begin{equation*}
0 = a^{\prime \prime}(0) \phi_{\omega} + 
S_{\omega}(\nu_{\omega})v^{\prime \prime}(0) 
+ 2 J^{-1}v^{\prime}(0). 
\end{equation*}
Multiplying the above identity by $\phi_{\omega}$ and 
integrating the resulting equation, we have 
\begin{equation*}
0 = a^{\prime \prime}(0) 
\|\phi_{\omega}\|_{L_{x}^{2}(\R) \times L_{x}^{2}(\R)}^{2} 
+ 2 \langle J^{-1} v^{\prime}(0), \phi_{\omega} \rangle. 
\end{equation*}
This together with \eqref{identity;v'(0)} yields that 
\begin{equation*}
a^{\prime \prime}(0) 
= -2 \frac{\langle J^{-1}v^{\prime}(0), \phi_{\omega} \rangle}
{\|\phi_{\omega}\|_{L_{x}^{2}(\R) \times L_{x}^{2}(\R)}^{2}} 
= 2 \frac{\langle J^{-1}v^{\prime}(0), JS_{\omega}(\nu_{\omega})v^{\prime}(0) \rangle}
{\|\phi_{\omega}\|_{L_{x}^{2}(\R) \times L_{x}^{2}(\R)}^{2}} 
= - 2 \frac{\langle v^{\prime}(0), S_{\omega}(\nu_{\omega})
v^{\prime}(0) \rangle}
{\|\phi_{\omega}\|_{L_{x}^{2}(\R) \times L_{x}^{2}(\R)}^{2}} <0.
\end{equation*}
In the last inequality, we have used the fact that 
$v^{\prime}(0) \in (\phi_{\omega})^{\perp}$, 
$\mathop{\mathrm{Ker}} S_{\omega}(\nu_{\omega})
= \mathop{\mathrm{Span}} (\phi_{\omega})$ and $S_{\omega}(\nu_{\omega})$ 
is non-negative operator. 
Therefore, there exists a sufficiently small $\varepsilon >0$ such that the function $a(\lambda)$ 
on $(0, \varepsilon)$ has the inverse function $\lambda(a)$
defined on a neighborhood of $a(0)=\nu_\omega$.
In other words, $-JS_{\omega}(a)$ has a positive eigenvalue 
$\lambda(a)$ for $a \in (\lambda(\varepsilon), \nu_{\omega})$.

Let 
\begin{equation*}
a_0 := \inf\{ b > 0 \colon 
\mbox{$-JS_{\omega}(b)$ has a positive eigenvalue} 
\}.
\end{equation*}

From the above argument, we see that $a_{0} < \nu_{\omega}$. We shall show that $a_0 = 0$. Suppose the contrary that $0 < a_0 < \nu_{\omega}$. We take $\{a_n\}_{n=1}^{\infty} \subset (a_0, \nu_{\omega})$ such that 
\[
\lim_{n \to \infty} a_n = a_0.
\]
By the perturbation theory, the sequence $\{a_{n}\}$ satisfies 
either 
\[
\lim_{n \to \infty} \lambda(a_n) = 0 \quad
{\rm or} \quad \lim_{n \to \infty}\lambda(a_n) = \infty.
\]
We claim that 
the second limit cannot hold. 
Indeed, taking $u_n\neq 0$ the associated eigenfunction corresponding to the eigenvalue $\lambda(a_n)$ with 
$\|u_{n}\|_{L_{x}^{2}(\R) \times L_{x}^{2}(\R)} = 1$, 
we see that
\begin{equation} \label{iden-JS}
\begin{split}
\langle -J S_{\omega}(a_n)u_n,u_n \rangle 
&= \langle (L_{\omega, -, a_{n}})\Re u_n, \Im u_n \rangle - 
\langle (L_{\omega, +, a_{n}}) \Im u_n, \Re u_n \rangle \\
& = -(p-1) \langle R^{p-1}_{\omega}\Im u_n, \Re u_n \rangle.
\end{split}
\end{equation}
Since $-J S_{\omega}(a_n)u_n = \lambda(a_{n}) u_{n}$ and $\|R_{\omega}\|_{L_{x}^\infty(\R)}\lesssim 1$, it follows from \eqref{iden-JS} that 
\begin{equation*}
\lambda(a_n) \|u_n\|_{L_{x}^2(\R) \times L_{x}^{2}(\R)}^{2} 
= |\langle -J S_{\omega}(a_n)u_n, u_n \rangle| \lesssim 
\|u_n\|_{L_{x}^2(\R) \times L_{x}^{2}(\R)}^{2},
\end{equation*}
which means that $\lambda(a_n)\lesssim 1$.
As a consequence, we obtain that 
\begin{equation*}
\lim_{n \to \infty} \lambda(a_n) = 0. 
\end{equation*}
Next, 
%let $u_n \in L_{x}^{2}(\mathbb{R}) \times L_{x}^{2}(\mathbb{R})$ be the eigenfunction of the operator $-JS_{\omega}(a_n)$ corresponding to $\lambda(a_{n})$ with $\|u_n\|_{L_{x}^{2}(\R) \times L_{x}^{2}(\R)}=1$. 
since $-JS_{\omega}(a_n)u_n=\lambda(a_n)u_n$, we have 
\begin{equation*}
S_{\omega}(a_0)u_n = 
S_{\omega}(a_{n}) u_{n} + (a_{0} - a_{n})u_{n} =
- J^{-1} \lambda (a_{n}) u_{n} + (a_{0} - a_{n}) u_{n}.
\end{equation*}
Since $S_{\omega}(a_0)$ is invertible 
and the inverse $(S_{\omega}(a_0))^{-1}$ is bounded, 
we obtain
\begin{equation*}
\begin{split}
1 = \|u_n\|_{L_{x}^{2}(\R) \times L_{x}^{2}(\R)}&= 
\|(S_{\omega}(a_0))^{-1} (- J^{-1} 
\lambda (a_{n}) u_{n} + (a_{0} - a_{n}) u_{n})\|_{L_{x}^{2}(\R)} \\
&\lesssim 
\lambda(a_{n}) + |a_{0} - a_{n}|
\to 0 
\end{split}
\end{equation*}
as $n \to \infty$, 
which is a contradiction. 
Thus, we conclude that $a_0 = 0$, which implies that 
for $0 < a < \nu_{\omega}$, $-JS_{\omega}(a)$ 
has a positive eigenvalue. 
Since $\nu_{\omega} > 1$ for $\omega > \omega_{p}$, 
we see that $-JS_{\omega}(1)$ has a positive 
eigenvalue. 
This completes the proof. 
\end{proof}

For each $k \in \mathbb{Z}$, 
define an orthogonal projection $P_{\leq k}$ as
\begin{equation} \label{projection}
P_{\leq k} u(x,y) = \sum_{n=-k}^k u_n(x) e^{iny}, \quad (x,y) \in \mathbb{R} \times \mathbb{T}, 
\end{equation}
where
\begin{equation*}
u(x,y) = \sum_{n=-\infty}^{\infty} u_n(x)e^{iny}, \qquad 
u_{n}(x) = \frac{1}{2\pi} \int_{\mathbb{T}} u(x, y) e^{- i n y} dy. 
\end{equation*}
From Lemma \ref{lem:posi-eigen}, 
we see that there exists a
positive eigenvalue of 
$- J \mathcal{S}_{\omega}^{\prime \prime}
(R_{\omega})$ for $\omega > \omega_{p}$. 
Thus, we can define 
\begin{equation} \label{spectral-radius}
\lambda_{0} := \max \left\{ \lambda >0 \colon\;
\lambda \in \sigma_{\text{disc}} 
(- J \mathcal{S}_{\omega}^{\prime \prime}(R_{\omega})) \right\}. 
\end{equation}
where 
$\sigma_{\text{disc}}(- J \mathcal{S}_{\omega}
^{\prime \prime}(R_{\omega}))$ 
is the set of the discreate spectrum of $- J 
\mathcal{S}_{\omega}^{\prime \prime}(R_{\omega})$. 
We can estimate the semigroup for the low frequency part. 
\begin{lemma}%[\cite{Yamazaki1}]
\label{lem-3-2}
For a positive integer $k$ and $\varepsilon > 0$, there exists $C_{k, \varepsilon}>0$ 
such that
\begin{equation} \label{semi-est}
\left\|e^{-tJS^{\prime \prime}_{\omega}(R_{\omega})}P_{\leq k}v
\right\|_{L_{x, y}^2(\R \times \T)} \leq 
C_{k, \varepsilon} e^{(\lambda_{0} + \varepsilon)t}\|v\|_{L^2_{x, y}(\R \times \T)}, 
\quad t>0, \ v \in 
L_{x, y}^2(\mathbb{R} \times \mathbb{T},\mathbb{C}). 
\end{equation}
\end{lemma}
\begin{proof}
A similar semigroup estimate has already been obtained for example in 
\cite[Lemma 3.2]{Yamazaki3} and \cite[Lemma 3.3]{Yamazaki1}. Following their argument, we can prove estimate \eqref{semi-est}, and for the sake of completeness, 
we shall give the proof here.

By definition of $S(a)$, we have 
\begin{equation*}
- J S(a) = 
\begin{pmatrix}
0 & - \partial_{x}^{2} + \omega + |a| - R_{\omega}^{p-1} \\
\partial_{x}^{2} - \omega - |a| + p R_{\omega}^{p-1} & 0
\end{pmatrix}
\end{equation*}
Applying the argument of Proposition 10 
in \cite{Georgiev-Ohta}, we have 
\begin{equation} \label{spectral-mapping}
\sigma(e^{- JS(a)}) = e^{\sigma(- JS(a))}. 
\end{equation}
This together with \eqref{eq-linearizedop} and 
\eqref{spectral-radius} implies that 
the spectral radius of $e^{-JS(n)}$ is 
less than or equal to $e^{\lambda_{0}}$. 
Therefore, by Lemmas 2 and 3 in \cite{Shatah-Strauss2}, 
for any $\varepsilon>0$ and each $n \in \mathbb{N}$, 
there exists $C_{n, \varepsilon}>0$ such that 
\begin{equation} \label{spectral-mapping-eq1}
\|e^{- t J S(n)}v\|_{L_{x}^{2}(\R)} \leq C_{n, \varepsilon} 
e^{(\lambda_{0} + \varepsilon) t} \|v\|_{L_{x}^{2}(\R)} \qquad 
\mbox{for all $v \in L_{x}^{2}(\mathbb{R})$}. 
\end{equation}
%Indeed, 
%let $\lambda (a)$ be the spectral radius of $- JS(a)$. 
%Then, we have $\lambda (a) \leq \lambda_{0}$ for all $a \in \mathbb{R}$. 
%By the definition of $\lambda(a)$ and \eqref{spectral-mapping}, 
%the spactral radius of $e^{- JS(a)}$ is
%\begin{equation}
%e^{\Re \lambda(a) } = |e^{\lambda(a) }| = 
%\lim_{m \to \infty} \|e^{m L}\|^{1/m}. 
%\end{equation}
%Therefore, for any $\varepsilon >0$, there exists a constant 
%$S_{\varepsilon}>0$ such that 
%\begin{equation} \label{eq1-spectral-radius}
%e^{\Re \lambda(a) - \varepsilon} < \|e^{m L}\|^{1/m} 
%< e^{\Re \lambda(a) + \varepsilon}
%\end{equation}
%for every integer $m \geq S_{\varepsilon} -1$. 
%Let $t > S_{\varepsilon}$ and $m = [t]$. 
%Then, from \eqref{eq1-spectral-radius}, we have 
%\begin{equation}
%\|e^{- JS(a) t}\| \leq K \|e^{- JS(a) m}\| 
%< K e^{(\Re \lambda (a) + \varepsilon) m} 
%\leq K e^{(\Re \lambda(a) + \varepsilon) t} 
%\leq K e^{(\lambda_{0} + \varepsilon) t}, 
%\end{equation} 
%where $K =\sup\left\{\|\exp(\theta L)\| \colon 0 \leq \theta \leq 1 \right\}$. 
%This yields \eqref{spectral-mapping-eq1}. 

\noindent
Hence, for $t >0$ and $v \in L_{x, y}^{2}(\mathbb{R} \times \T)$, we have 
\begin{equation*}
\|e^{- t J \mathcal{S}_{\omega}^{\prime \prime}(R_{\omega})} 
P_{\leq k} v\|_{L^{2}_{x, y}(\R \times \T)} 
\leq \left\|\sum_{n = -k}^{k} 
e^{- t J S(n)} v_{n} e^{i n y}\right\|_{L_{x, y}^{2}(\R \times \T)} 
\leq C_{k, \varepsilon} e^{(\lambda_{0} + \varepsilon) t} \|v\|_{L_{x, y}^{2}
(\R \times \T)}, 
\end{equation*}
where 
\begin{equation*}
v(x, y) = \sum_{n \in \mathbb{Z}} v_{n}(x) e^{i n y}. 
\end{equation*}
This completes the proof. 
\end{proof}
Next, we control the $L^2$ norm of the low frequencies 
for the nonlinear term by the energy norm.% for $1<p\leq 5$. %we recall the following lemma ( see the proof of Lemma 3.2 in \cite{Yamazaki1} for more details).
\begin{lemma}
\label{lem-3-1}
For $k>1$, there exists $C>0$ such that
\begin{equation}\|P_{\leq k} NL(v,R_\omega)\|_{L^2_{x,y}(\R \times \T)} 
\leq 
\begin{cases}
C\|v\|_{X}^{p} \qquad &\mbox{if $1<p\leq 2$},\\
C k^{\frac{1}{2}}(\|v\|_{X}^2+\|v\|_{X}^{p}) \qquad 
&\mbox{if $2<p\leq 5$},
\end{cases}
\end{equation}
$NL(v,R_\omega)$ is given by \eqref{linearized-nonlinear} and $v\in X$.
\end{lemma}

\begin{remark}
Note that it seems impossible to 
uniformly control the 
$L^2$ norm of the nonlinear term by the energy norm when $p<3$. 
This is the reason why we are restricted only 
to low-frequencies 
which is sufficient for our purpose. 
\end{remark}

\begin{proof}
We shall split the proof into two cases according to the value of $p$.

\textbf{(Case 1).} $1<p\leq 2$. 
%In this case, we need not take the projection $P_{\leq k}$ and can estimate 
%the nonlinear term as it is. 
From \eqref{linearized-nonlinear}, we obtain 
\begin{equation*}
\begin{split}
NL(v, R_\omega)
& = \int_{0}^{1} \frac{d}{d \theta}
(f(\theta v + R_{\omega})) d \theta 
- D f(R_{\omega}) v \\
& = \int_0^1 
(D f(\theta v + R_\omega) - D f(R_\omega))v(x,y)d\theta.
\end{split}
\end{equation*}
From the result of \cite[Lemma 2.4]{ginibrle-velo}, we have
\begin{equation}
\label{estim:p-1}
||a|^{p-1}-|b|^{p-1}|\leq |a-b|^{p-1},
\end{equation}
which together with \eqref{eq-deri2} yields
\begin{equation*}
|NL(v, R_\omega)| \leq |v(x,y)|^{p} 
\qquad \mbox{for $1 < p \leq 2$}.
\end{equation*}
Hence, we have
\begin{equation*}
\|NL(v,R_\omega)\|_{L^2_{x,y}(\R \times \T)} \leq \|v\|_{L^{2p}_{x,y}(\R \times \T)}^p.
\end{equation*}
On the other hand, since $X \hookrightarrow L^q_{x,y}
(\R \times \T)$ for $2\leq q\leq 6$, 
we obtain 
\begin{equation*}
\|NL(v,R_\omega)\|_{L^2_{x,y}(\R \times \T)} \lesssim \|v\|_{X}^p.
\end{equation*}
\par
\textbf{(Case 2).} $2 < p \leq 5$.
Note that the embedding $X \hookrightarrow 
L^{2p}_{x,y}(\mathbb{R} \times \mathbb{T})$ 
does not hold for $p >3$. 
Thus, we shall restrict ourselves to the low-frequencies only and 
use the following Bernstein estimate (see, for example, 
\cite[Lemma 11.4]{K T V}) 
\begin{equation} \label{bernstein}
\|P_{\leq k} NL(v(x, \cdot),R_\omega)\|_{L^2_y(\T)} 
\lesssim k^\frac{1}{2} \| NL(v(x, \cdot),R_\omega)\|_{L^1_y(\T)}.
\end{equation}
We put 
\begin{equation*}
F(s) := f(sv + R_{\omega}), 
\end{equation*}
where $f(z) = |z|^{p-1}z$. 
We note that $F$ is $C^{2}((0, \infty), \mathbb{R})$ 
is $C^{2}$ if $p > 2$. 
From \eqref{linearized-nonlinear}, we can write 
\begin{equation*}
NL(v,R_\omega)= \int^1_0 (1-s) F^{\prime \prime}(s) ds.
\end{equation*} 
By the definition of $F(s)$, we have
\begin{equation*}
\begin{split}
F^{\prime \prime}(s)
= & \frac{p^{2} -1}{4} |s v + R_{\omega}|^{p-3} (s\overline{v} 
+ R_{\omega})v^{2} 
+ \frac{p^{2} -1}{2} |sv + R_{\omega}|^{p-3}
(sv + R_{\omega}) |v|^{2} \\
& + \frac{(p-1)(p-3)}{4} |sv + R_{\omega}|^{p-5}(sv + R_{\omega})^{3}
(\overline{v})^{2}. 
\end{split}
\end{equation*}
Thus, since $p>2$, we obtain
\begin{equation*}
F^{\prime \prime}(s) \lesssim |v|^2 |sv+R_\omega|^{p-2} \lesssim |v|^2 \left(|v|^{p-2} + |R_\omega|^{p-2}\right).
\end{equation*}
As a consequence, from the fact that 
$\|R_\omega\|_{L^\infty_x(\R)} \lesssim 1$, we get
\begin{equation*}
|NL(v,R_\omega)| \lesssim |v|^2 + |v|^p.
\end{equation*} 
Then, we obtain 
\begin{equation*}
\| NL(v(x, \cdot),R_\omega)\|_{L^1_y} \lesssim \|v(x, \cdot)\|_{L^2_y}^2+\|v(x, \cdot)\|_{L^p_y}^p.
\end{equation*}
Next, taking the $L^2$-norm in $x$ variable for the previous estimate, 
we have
\begin{equation}
\label{estim-NL-L2-L1}
\| NL(v,R_\omega)\|_{L^2_xL^1_y(\R \times \T)} 
\lesssim \|v\|_{L^4_x L^2_y(\R \times \T)}^2+\|v\|_{L^{2p}_x L^p_y(\R \times \T)}^p.
\end{equation}
It follows from \cite[Theorem 10.2]{BON} that 
$X \hookrightarrow L^{q_{1}}_{x}L^{q_{2}}_{y}(\mathbb{R} \times \mathbb{T})$ if 
$q_{1}$ and $q_{2}$ satisfies 
\begin{equation*}
q_{1}, q_{2} \geq 2, \qquad
\frac{1}{4} \leq \frac{1}{2 q_{1}} + \frac{1}{q_{2}}. 
\end{equation*}
This yields that
\begin{equation}
\label{estim:L2p-Lp}
\|v\|_{L^{2p}_x L^p_y(\R \times \T)} 
\lesssim \|v\|_{X} \quad {\rm for } \quad 2 \leq p \leq 5.
\end{equation}
Combining \eqref{bernstein} with \eqref{estim-NL-L2-L1} and \eqref{estim:L2p-Lp} finishes the proof in this case, and therefore completes the proof of this lemma.
\end{proof}

Lemma \ref{lem:posi-eigen} guarantees that there exists 
a positive eigenvalue $\lambda_{0}$ 
of the operator $- JS_{\omega}(1)$ 
and the corresponding eigenfunction 
$\chi_{0} \in H_{x}^2(\mathbb{R},\mathbb{C})$
with $\|\chi_{0}\|_{L_{x, y}^2(\R \times \T)}=1$. 
Then, putting 
\begin{equation}
\label{def:chi_0}
\chi(x,y)=\chi_0(x) e^{i y}, 
\end{equation}
we see that $\lambda_{0}$ is also the eigenvalue of 
the operator 
$-J \mathcal{S}^{\prime \prime}_{\omega}(R_{\omega})$ 
and $\chi \in H^2_{x}H^{1}_{y}
(\mathbb{R} \times \mathbb{T},\mathbb{C})$ 
is the corresponding eigenfunction. 

For $\delta>0$, let $u_{\delta}(t)$ be the solution of 
(\ref{WS-Hamil}) 
with initial data $u_{\delta}|_{t=0} 
= R_{\omega}+\delta\chi$ and 
$v_{\delta}(t)$ be the solution of (\ref{WS-Ham-v}) 
having an initial data $v_{\delta}|_{t=0}=\delta \chi$.
Then, we infer that $u_{\delta}(t)=e^{i\omega t}(R_{\omega}+v_{\delta}(t))$. 
We have the following:
\begin{lemma}\label{lem-3-3}
There exist a positive integer $K_\omega$ and 
positive constant $C$ such that for $\delta>0$ and $t>0$, 
\begin{equation*}
\|v_{\delta}(t)\|_{X} \leq 
C\|P_{\leq K_\omega}v_{\delta}(t)\|
_{L^2_{x, y}(\R \times \T)} + o(\delta).
\end{equation*}
\end{lemma}
\begin{proof}
Using the Taylor expansion and the fact that 
$\mathcal{S}^{\prime}_{\omega}(R_{\omega})=0$, for $v \in X$ and any fixed $t>0,$ we write
\begin{equation} \label{eq1-lem3-3}
\begin{split}
\mathcal{S}_{\omega}(u_{\delta}(t))
& = \mathcal{S}_{\omega}(R_{\omega} + v_{\delta}(t)) \\
& = 
\mathcal{S}_{\omega}(R_{\omega}) + 
\langle \mathcal{S}_{\omega}^{\prime}(R_{\omega}), v_{\delta}(t)\rangle_{X^{*}, X} 
+ \frac{1}{2} \langle \mathcal{S}_{\omega}^{\prime \prime}(R_{\omega})v_{\delta}(t), 
v_{\delta}(t)\rangle_{X^{*}, X} \\ 
& \quad+ o(\|v_{\delta}(t)\|_{X}^{2}) \\
& = 
\mathcal{S}_{\omega}(R_{\omega}) 
+ \frac{1}{2} \langle \mathcal{S}_{\omega}^{\prime \prime}(R_{\omega})v_{\delta}(t), 
v_{\delta}(t)\rangle_{X^{*}, X}+o(\|v_{\delta}(t)\|_{X}^{2}), 
\end{split}
\end{equation}
and, at $t=0$,
\begin{equation*}
\mathcal{S}_{\omega}(R_{\omega}+\delta \chi)
= \mathcal{S}_{\omega}(R_{\omega}) + \frac{\delta^{2}}{2}
\langle \mathcal{S}_{\omega}^{\prime \prime}(R_{\omega})\chi,\chi\rangle_{X^{*}, X}+o(\delta^{2}). 
\end{equation*}
From the conservation laws of the mass and the Hamiltonian \eqref{conserv}, we have 
\begin{equation} \label{eq2-lem3-3}
\mathcal{S}_{\omega}(u_{\delta}(t)) 
= \mathcal{S}_{\omega}(R_{\omega}+v_{\delta}(t))
= \mathcal{S}_{\omega}(R_{\omega}+\delta \chi)
\end{equation}
for any $t\geq 0$.
Moreover, 
since $J$ is a skew-symmetric operator, 
we see that
\begin{equation} \label{eq3-lem3-3}
\langle \mathcal{S}_{\omega}^{\prime \prime}(R_{\omega})\chi,\chi\rangle_{X^{*}, X} = 
- \langle J\mathcal{S}_{\omega}^{\prime \prime}
(R_{\omega})\chi, J^{-1}\chi\rangle_{X^{*}, X} = 
(\lambda_0 \chi, J^{-1}\chi)_{L_{x, y}^2(\R \times \T)}=0. 
\end{equation}
Thus, from \eqref{eq1-lem3-3}--\eqref{eq3-lem3-3}, 
we infer that
\begin{equation}\label{eq4-lem3-3}
\langle \mathcal{S}_{\omega}^{\prime \prime}
(R_{\omega})v_{\delta}(t), 
v_{\delta}(t)\rangle_{X^{*}, X} 
= o(\|v_{\delta}(t)\|_{X}^2)+o(\delta^2).
\end{equation}
Let 
\[
K_{\omega} = \max\left\{ 
k \in \mathbb{Z} \colon k \leq 1 + \frac{\omega}{\omega_{p}}
\right\}.
\]
Since $S_\omega(n)$ is also positive 
for $|n| > 1 + \frac{\omega}{\omega_p}$, 
we see that 
$\mathcal{S}_{\omega}^{\prime \prime}(R_{\omega})P_{> K_{\omega}}$ is positive, 
that is, there exists a constant $C_{0} > 0$ such that 
\begin{equation}\label{eq6-lem3-3}
\langle \mathcal{S}_{\omega}^{\prime \prime}(R_{\omega})
P_{> K_{\omega}}v_{\delta}(t), P_{> K_{\omega}}v_{\delta}(t) \rangle _{X^{*}, X} 
\geq C_{0} \|P_{> K_{\omega}}v_{\delta}(t)\|_{X}^{2}. 
\end{equation}
By the definition of $S_\omega(n)$, 
we can take constants $C_{1} \in (0, C_{0})$ and 
$C_{2} > 0$ such that
\begin{equation} \label{eq5-lem3-3}
\langle S_\omega(n)P_{\leq K_{\omega}}v_{\delta}(t), 
P_{\leq K_{\omega}}v_{\delta}(t) 
\rangle_{X^{*}, X} \geq 
C_{1} \|P_{\leq K_{\omega}}v_{\delta}(t)\|_{H_{x}^1(\R)}^2 
-C_{2} \|P_{\leq K_{\omega}}v_{\delta}(t)\|_{L_{x}^2(\R)}^2.
\end{equation}
Thus, from \eqref{eq4-lem3-3}, \eqref{eq6-lem3-3} 
and \eqref{eq5-lem3-3}, 
we obtain 
\begin{equation*} 
\begin{split}
\|v_{\delta}(t)\|_{X}^2
& =
\|P_{\leq K_{\omega}}v_{\delta}(t)\|_{X}^2 + \|P_{> K_{\omega}}v_{\delta}(t)\|_{X}^2\\
& \leq \|P_{\leq K_{\omega}}v_{\delta}(t)
\|_{H_x^{1}L^2_y(\R \times \T)}^2 
+ K_{\omega} \|P_{\leq K_{\omega}}v_{\delta}(t)
\|_{L_{x,y}^{2}(\R \times \T)}^2 
+ \|P_{> K_{\omega}}v_{\delta}(t)\|_{X}^2\\
&\leq \frac{1}{C_{1}} \langle \mathcal{S}_{\omega}^{\prime \prime}(R_{\omega})
P_{\leq K_{\omega}}v_{\delta}(t),P_{\leq K_{\omega}}v_{\delta}(t)\rangle_{X^{*}, X} 
+ \left(K_{\omega} + \frac{C_{2}}{C_{1}}\right) \|P_{\leq K_{\omega}}v_{\delta}(t)\|_{L^{2}_{x,y}
(\R \times \T)}^2\\
& \quad + \frac{1}{C_{0}} \langle \mathcal{S}_{\omega}^{\prime \prime}
(R_{\omega})P_{> K_{\omega}}v_{\delta}(t),P_{> K_{\omega}}v_{\delta}(t)\rangle_{X^{*}, X} 
\\
& \leq 
\left(K_{\omega} + \frac{C_{2}}{C_{1}}\right) \|P_{\leq K_{\omega}}v_{\delta}(t)\|_{L^{2}_{x,y}
(\R \times \T)}^2 
+ \frac{1}{C_{1}} \langle \mathcal{S}_{\omega}^{\prime \prime}(R_{\omega})
v_{\delta}(t), v_{\delta}(t)\rangle_{X^{*}, X} \\
&\leq \left(K_{\omega} + \frac{C_{2}}{C_{1}}\right)\|P_{\leq K_{\omega}}v_{\delta}(t)\|_{L^2_{x,y}(\R \times \T)}^2 
+ o(\delta^2) + o(\|v_{\delta}(t)\|_{X}^2), 
\end{split}
\end{equation*} 
where we have used the fact that $0 < C_{1} < C_{0}$ 
in the third inequality and 
the orthogonality between $P_{\leq K_{\omega}}$ and $\mathcal{S}_{\omega}^{\prime \prime}(R_{\omega})P_{> K_{\omega}}=P_{> K_{\omega}}\mathcal{S}_{\omega}^{\prime \prime}(R_{\omega})$. This finishes the proof of Lemma \ref{lem-3-3}.
\end{proof}
Now, we give the proof of Theorem \ref{thm-1} (ii).

%Let $\varepsilon _0 =\min\{(p-1)\lambda_0/2,\lambda_0/2\}$.
\begin{proof}[Proof of Theorem \ref{thm-1} (ii)] 
From the Duhamel formula associated to \eqref{WS-Ham-v}, 
$v_{\delta}$ satisfies the following: 
\begin{equation*}
\begin{split}
v_{\delta}(t) 
& = e^{- J \mathcal{S}_{\omega}^{\prime \prime}(R_{\omega})t} \delta \chi 
- J \int_{0}^{t} e^{- (t-s) J \mathcal{S}_{\omega}^{\prime \prime}(R_{\omega})}
NL(v_{\delta}(s), R_{\omega}) ds \\
& = \delta e^{\lambda_{0}t} \chi 
- J \int_{0}^{t} e^{- (t-s) J \mathcal{S}_{\omega}^{\prime \prime}(R_{\omega})}
NL(v_{\delta}(s), R_{\omega}) ds. 
\end{split}
\end{equation*}

Using Lemmas \ref{lem-3-3}, \ref{lem-3-2} and \ref{lem-3-1} 
in this order, we estimate the $X$ norm of $v_{\delta}$ 
in the following way
\begin{equation*}
\begin{split}
\|v_{\delta}(t)\|_{X} 
& \lesssim \delta e^{\lambda_{0} t} \|\chi\|_{X} 
+\int_{0}^{t} \|e^{- (t-s) J \mathcal{S}_{\omega}^{\prime \prime}(R_{\omega})}
P_{\leq K_{\omega}} N
L(v_{\delta}(s), R_{\omega})\|_{L_{x, y}^{2}(\R \times \T)}ds 
+ o(\delta) \\
&\lesssim \delta e^{\lambda_0t}
\|\chi\|_{L_{x, y}^2(\R \times \T)} 
+ \int_0^t e^{(\lambda_{0} + \varepsilon)(t-s)}
\|P_{\leq K_{\omega}}
NL(v_{\delta}(s),R_\omega)\|_{L_{x, y}^2(\R \times \T)} ds
+ o(\delta) \\
& \lesssim \delta e^{\lambda_0 t} + \int_0^t e^{(1+\varepsilon _0)\lambda_0(t-s)}(\|v_{\delta}(s)\|_{X}^2+\|v_{\delta}(s)\|_{X}^p)ds,
\end{split}
\end{equation*}
where $\varepsilon_0=\varepsilon/\lambda_0$ and $\varepsilon>0$ given in Lemma \ref{lem-3-2}. Thus, there exist constants $C_{1} > 1$ and $C_{2} >0$ such that 
\begin{equation*}
\|v_{\delta}(t)\|_{X} 
\leq C_{1} \delta e^{\lambda_0 t} + 
C_{2} \int_0^t e^{(1+\varepsilon _0)\lambda_0(t-s)}
(\|v_{\delta}(s)\|_{X}^2+\|v_{\delta}(s)\|_{X}^p)ds.
\end{equation*} 
We shall show that for sufficiently small 
$\delta>0$ and $\varepsilon_0 >0$, we have 
\begin{equation} \label{up-bound}
\|v_{\delta}(t)\|_{X} \leq 2 C_{1} \delta e^{\lambda_0 t} \quad 
\mbox{for $t \in [0,T_{\varepsilon_0 ,\delta}]$},
\end{equation}
where
\begin{equation*}
T_{\varepsilon _0,\delta}=\frac{1}{\lambda_0} 
\log \left(\frac{\varepsilon_0}{\delta}\right). 
\end{equation*}
We put 
\begin{equation*}
T_{*} = \sup\{t > 0 \colon \|v_{\delta}(t)\|_{X} \leq 2 C_{1} \delta e^{\lambda_{0}t} 
\}. 
\end{equation*}
Since $\|v_{\delta}(0)\|_{X} = \delta (< 2C_{1} \delta)$, 
we see that $T_{*} > 0$. 
In order to prove \eqref{up-bound}, 
it is enough to show that $T_{\varepsilon, \delta} \leq T_{*}$. 
Suppose the contrary that $T_{*} < T_{\varepsilon, \delta}$. 
Then, we have 
\begin{equation} \label{nonlinear-est1}
\begin{split}
2 C_{1} \delta e^{\lambda_{0} T_{*}} 
%& 
= \|v_{\delta}(T_{*})\|_{X} %\\& 
\leq C_{1} \delta e^{\lambda_{0} T_{*}} 
+ \int_{0}^{T_{*}} e^{(1 + \varepsilon_{0}) \lambda_{0} (T_{*} - s)} 
(\|v_{\delta}(s)\|_{X}^{2} + \|v_{\delta}(s)\|_{X}^{p}) ds. 
\end{split}
\end{equation}
Note that for $0 < s < T_{*} (< T_{\varepsilon_{0}, \delta})$, 
we obtain 
\begin{equation} \label{nonlinear-est2}
2 C_{1} \delta e^{\lambda_{0} s} < 
2 C_{1} \delta e^{\lambda_{0} T_{\varepsilon_{0}, \delta}} = 2 C_{1} \varepsilon_{0} \ll 1. 
\end{equation}
It follows from \eqref{nonlinear-est1} 
and \eqref{nonlinear-est2} 
that
\begin{equation*}
\begin{split}
C_{1} \delta e^{\lambda_{0} T_{*}} 
& \leq 
\int_{0}^{T_{*}} e^{(1 + \varepsilon_{0}) \lambda_{0} (T_{*} - s)} 
\left[(2C_{1} \delta)^{2} e^{2 \lambda_{0} s} 
+ (2C_{1} \delta)^{p} e^{p \lambda_{0} s}\right]ds \\
& \leq 8 C_{1}^{2} \delta^{2} e^{(1 + \varepsilon_{0})\lambda_{0} T_{*}} 
\int_{0}^{T_{*}} e^{(1 - \varepsilon_{0}) \lambda_{0} s}ds \\
& = 8 C_{1}^{2} \delta^{2} e^{(1 + \varepsilon_{0})\lambda_{0} T_{*}} 
\frac{e^{(1-\varepsilon_{0}) \lambda_{0} T_{*}}}{(1-\varepsilon_{0}) \lambda_{0}} \\
& = \frac{8 C_{1}^{2} \delta^{2}}{(1-\varepsilon_{0}) \lambda_{0}} e^{2 \lambda_{0}T_{*}}. 
\end{split}
\end{equation*}
This together with $T_{*} < T_{\varepsilon, \delta}$ 
yields that 
\begin{equation*}
1 \leq 
\frac{8 C_{1}^{2} \delta}{(1 - \varepsilon_{0}) \lambda_{0}} 
e^{\lambda_{0} T_{\varepsilon, \delta}} 
= \frac{8 C_{1}^{2} \delta}{(1 - \varepsilon_{0}) \lambda_{0}} 
\times \frac{\varepsilon_{0}}{\delta} 
= \frac{8 C_{1}^{2}}{(1 - \varepsilon_{0}) \lambda_{0}} 
\varepsilon_{0} < 1
\end{equation*}
for sufficiently small $\varepsilon_{0} > 0$, 
which is a contradiction. 
Thus, \eqref{up-bound} holds. 

Since $(e^{J \mathcal{S}_{\omega}^{\prime \prime} (R_{\omega})})^{*} 
= e^{(J \mathcal{S}_{\omega}^{\prime \prime}(R_{\omega}))^{*}}$ and 
\begin{equation*}
(J \mathcal{S}_{\omega}^{\prime \prime}(R_{\omega}))^{*} J \chi 
= \mathcal{S}_{\omega}^{\prime \prime}(R_{\omega}) J^{*} J \chi
= - \mathcal{S}_{\omega}^{\prime \prime}(R_{\omega}) \chi
= J(J \mathcal{S}_{\omega}^{\prime \prime}(R_{\omega}) \chi) 
= - \lambda_{0} J \chi, 
\end{equation*}
we have, by \eqref{up-bound}, that 
\begin{equation} \label{eq1-proof-insta}
\begin{split}
& \quad |(\chi,v_{\delta}(T_{\varepsilon _0,\delta}))
_{L_{x, y}^2(\R \times \T)}| \\
&=\left| \delta e^{\lambda_0 T_{\varepsilon _0, \delta}} + 
\int_0^{T_{\varepsilon _0,\delta}}(\chi, - J e^{-(T_{\varepsilon _0,\delta}-s)
J \mathcal{S}_{\omega}^{\prime \prime}(R_{\omega})}
NL(v_{\delta}(s),R_\omega))_{L^2}ds \right|\\
&\geq \varepsilon _0 -C \int_0^{T_{\varepsilon _0,\delta}} e^{(T_{\varepsilon _0,\delta}-s)\lambda_0}(\|v_{\delta}(s)\|_{X}^2+\|v_{\delta}(s)\|_{X}^p)ds\\
& \geq \varepsilon _0 -C \int_0^{T_{\varepsilon _0,\delta}} e^{(T_{\varepsilon _0,\delta}-s)\lambda_0} 
(\delta^{2} e^{2 \lambda_{0} s} + \delta^{p} e^{p \lambda_{0} s})ds \\ 
& \geq \varepsilon_0 - C \varepsilon_0 ^{\min \{p,2\}} \\
& \geq \frac{\varepsilon_{0}}{2}.
\end{split}
\end{equation}
From $P_{\leq 0} R_{\omega}=R_{\omega}$ and \eqref{eq-error},
there exists $\varepsilon _0 >0$ such that 
for $\varepsilon _0 > \delta >0$ and $\theta \in \R$, 
we have 
\begin{equation*}
\begin{split}
\|u_{\delta}(T_{\varepsilon _0,\delta})-e^{i\theta}R_{\omega}\|_{L^2_{x, y}(\R \times \T)}
&\geq \|(I-P_{\leq 0}) (u_{\delta}(T_{\varepsilon _0,\delta})-e^{i\theta}R_{\omega})\|_{L^2_{x, y}(\R \times \T)}\\
&= \|(I-P_{\leq 0})e^{-i\omega T_{\varepsilon _0,\delta}}u_{\delta}(T_{\varepsilon _0,\delta})\|_{L^2_{x, y}(\R \times \T)}\\
&= \|(I-P_{\leq 0})(e^{-i\omega T_{\varepsilon _0,\delta}}u_{\delta}(T_{\varepsilon _0,\delta})-R_{\omega})
\|_{L^2_{x, y}(\R \times \T)} \\
& \geq \|(P_{\leq 1} - P_{\leq 0})(e^{-i\omega T_{\varepsilon _0, \delta}}u_{\delta}(T_{\varepsilon _0,\delta})-R_{\omega})
\|_{L^2_{x, y}(\R \times \T)}. 
%\\
%& = \|(P_{\leq k_{0}} - P_{\leq k_{0}-1})v_{\delta}(T_{\varepsilon_{0}, %\delta})\|_{L^2_{x, y}(\R \times \T)}.
\end{split}
\end{equation*}
Let $v_{1}$ be the first Fourier mode of 
$v_{\delta}(T_{\varepsilon_{0}, \delta})$ on the $y$ variable.
From \eqref{def:chi_0} and 
the Cauchy-Schwarz inequality, we have
\begin{equation}
(\chi,v)_{L^2_{x,y}(\R \times \T)} 
= (e^{i y}\chi_0,v)_{L^2_{x,y}(\R \times \T)}= 
(\chi_0,v_{1})_{L^2_{x}} 
\leq \|\chi_0\|_{L^2_{x}(\R)} \|v_{1}\|_{L^2_{x}(\R)}
\end{equation} 
for $v \in L^2(\R \times \T)$. Thus, we have 
\begin{equation}\label{eq2-proof-insta}
\|(P_{\leq 1}-P_{\leq 0})v\|_{L^2_{x,y}(\R \times \T)} = 
\|v_{1}\|_{L^2_{x}(\R)} \geq \|\chi_0\|_{L^2_{x}(\R)}^{-1} 
|(\chi,v)_{L^2_{x,y}(\R \times \T)}| 
\end{equation}
for $v \in L_{x, y}^2(\R \times \T)$. 
It follows from \eqref{eq1-proof-insta}--
\eqref{eq2-proof-insta} that 
\begin{equation*}
\begin{split}
\|u_{\delta}(T_{\varepsilon _0,\delta})
-e^{i\theta}R_{\omega}\|_{L^2_{x,y}(\R \times \T)} 
&\geq \|(P_{\leq 1}-P_{\leq 0})v_\delta 
(T_{\varepsilon _0,\delta})\|_{L^2_{x,y}
(\R \times \T)} \\
&\gtrsim |(\chi,v_{\delta}(T_{\varepsilon _0,\delta}))_{L^2_{x,y}(\R \times \T)}| 
\geq \frac{\varepsilon _0}{2}.
\end{split}
\end{equation*}
This implies that the standing wave $e^{i\omega t}R_{\omega}$ is unstable.
\end{proof}

\section{Characterization of ground states}
In this section, following Berestycki and Wei~\cite{BW} and 
Terracini,Tzvetkov and Visciglia~\cite{TTV}, 
we shall show Theorem \ref{thm-chara-ground}. 
\subsection{Existence of ground states}
\label{sec-exist-ground}
We first give a proof of Theorem \ref{thm-exist-ground}. 
Namely, we show the existence of ground states for any $\omega > 0$. 
To this end, we prepare several lemmas, which are needed later. First, recall that
\[
m_{\omega} = \inf \left\{\mathcal{S}_{\omega}(u) \colon 
u \in X \setminus \{0\}, \; \mathcal{N}_{\omega}(u) = 0 \right\}, 
\]
where 
\[
\mathcal{N}_{\omega} (u) = 
\langle \mathcal{S}_{\omega}^{\prime}(u), u \rangle 
= 
\|\nabla u\|_{L_{x, y}^{2}(\R \times \T)}^{2} 
+ \||D_{y}|^{\frac{1}{2}} u\|_{L_{x, y}^{2}(\R \times \T)}^{2}
+ \omega \|u\|_{L_{x, y}^{2}(\R \times \T)}^{2} 
- \|u\|_{L_{x, y}^{p+1}(\R \times \T)}^{p+1}. 
\]
The following lemma provide an alternative form of $m_{\omega}$.
\begin{lemma} \label{existence-ground-thm2}
For any $\omega>0$, we have 
\begin{equation*}
m_{\omega} = \inf \left\{
\mathcal{I}_{\omega}(u) \colon \mathcal{N}_{\omega}(u) \leq 0 \right\}, 
\end{equation*}
where 
\begin{equation} \label{eq-0}
\mathcal{I}_{\omega}(u) := 
\left(\frac{1}{2} - \frac{1}{p+1}\right) 
\left\{\|\partial_{x} u\|_{L_{x, y}^{2}(\R \times \T)}^{2} 
+ \||D|_{y}^{\frac{1}{2}} u\|_{L_{x, y}^{2}(\R \times \T)}^{2} 
+ \omega \|u\|_{L_{x, y}^{2}(\R \times \T)}^{2} 
\right\}. 
\end{equation}
\end{lemma}
\begin{proof}
Let $u \in X$ satisfy $\mathcal{N}_{\omega}(u) = 0$. 
Then, from \eqref{action}, we obtain 
\begin{equation*}
\begin{split}
\mathcal{S}_{\omega}(u) 
& = \mathcal{S}_{\omega}(u) - \frac{1}{p+1} 
\mathcal{N}_{\omega}(u) \\
& = 
\left(\frac{1}{2} - \frac{1}{p+1}\right)
\left\{\|\partial_{x} u\|_{L_{x, y}^{2}(\R \times \T)}^{2} 
+ \||D|_{y}^{\frac{1}{2}} u\|_{L_{x, y}^{2}(\R \times \T)}^{2} 
+ \omega \|u\|_{L_{x, y}^{2}(\R \times \T)}^{2} \right\}
= \mathcal{I}_{\omega}(u). 
\end{split}
\end{equation*}
Therefore, we have 
\begin{equation} \label{eq-1-m-omega}
m_{\omega} = \inf\left\{\mathcal{I}_{\omega}(u) \colon \mathcal{N}_{\omega}(u) 
= 0 \right\}. 
\end{equation}
Suppose that $\mathcal{N}_{\omega}(u) < 0$. 
Then, there exists $t_{0} \in (0, 1)$ such that 
$\mathcal{N}_{\omega}(t_{0} u) = 0$. 
it follows from \eqref{eq-0} and 
\eqref{eq-1-m-omega} that 
\begin{equation*}
m_{\omega} \leq \mathcal{I}_{\omega}(t_{0}u) < 
\mathcal{I}_{\omega}(u). 
\end{equation*}
Taking an infimum on $u \in X$, we obtain \eqref{eq-0}. 
This completes the proof
\end{proof}
From the Sobolev inequality, 
we provide the following lemma:
\begin{lemma} \label{existence-ground-thm3}
Let $1 \leq p \leq 5$ and $\omega > 0$. We have 
\begin{equation*}
\|u\|_{L_{x, y}^{p+1}(\R \times \T)}^{p+1} \leq C 
\sup_{k \in \mathbb{Z}} 
\left(\|u\|_{L_{x, y}^{p+1}(Q_{k})} 
\right)^{p-1} \|u\|_{X}^{2}, 
\end{equation*}
where 
$
Q_{k} = [k, k+1) \times \mathbb{T}. 
$

\end{lemma}
\begin{proof}
It follows that 
\begin{equation*}
%\begin{split}
%& \quad 
\|Q\|_{L_{x, y}^{p+1}(\R \times \T)}^{p+1} 
%= \int_{\R \times \T} 
%|Q(x, y)|^{p+1} dxdy 
%= \sum_{k \in \mathbb{Z}} 
%\int_{Q_{k}} |Q(x, y)|^{p+1} dxdy 
%\\& 
= \sum_{k \in \mathbb{Z}} 
\|Q\|_{L^{p+1}_{x, y}(Q_{k})}^{p+1} 
\leq (\sup_{k \in \mathbb{Z}} 
\|Q\|_{L^{p+1}_{x, y}(Q_{k})})^{p-1} 
\sum_{k \in \mathbb{Z}} 
\|Q\|_{L^{p+1}_{x, y}(Q_{k})}^{2}. 
%\end{split}
\end{equation*}
By the Sobolev inequality, we have 
\begin{equation*} %\label{eq4-2-1}
\|Q\|_{L^{p+1}_{x, y}(Q_{k})} \leq C 
\|Q\|_{X(Q_{k})}, 
\end{equation*}
where $C > 0$ is the constant, 
which does not depend on $k \in 
\mathbb{Z}$ and 
\[
\|Q\|_{X(Q_{k})} 
= \left\{
\int_{Q_{k}} \left(
|\partial_{x} u(x, y)|^{2} + |D_{y}| u(x, y) 
\overline{u(x, y)} + |u(x, y)|^{2}
\right) dxdy
\right\}^{\frac{1}{2}}. 
\]
This yields that 
\[
\|Q\|_{L_{x, y}^{p+1}(\R \times \T)}^{p+1}
\leq C(\sup_{k \in \mathbb{Z}} 
\|Q\|_{L^{p+1}_{x, y}(Q_{k})})^{p-1} 
\sum_{k \in \mathbb{Z}} 
\|Q\|_{X(Q_{k})}^{2} 
= C(\sup_{k \in \mathbb{Z}} 
\|Q\|_{L^{p+1}_{x, y}(Q_{k})})^{p-1} 
\|Q\|_{X}^{2}.
\]
This completes the proof.
\end{proof}
We are now in position to prove Theorem \ref{thm-exist-ground}. 
\begin{proof}[Proof of Theorem \ref{thm-exist-ground}]
Let $\{u_{n}\}$ be a minimizing sequence for $m_{\omega}$. 
Then, for sufficiently large $n \in \mathbb{N}$, we have 
\begin{equation*}
\mathcal{I}_{\omega}(u_{n}) \leq m_{\omega} + 1. 
\end{equation*}
Thus, we see that $\{u_{n}\}$ is bounded in $X$. 
Then, it follows from the 
Sobolev inequality 
and $\mathcal{N}_{\omega}(u_{n}) = 0$ 
that 
\begin{equation*}
\|u_{n}\|_{L_{x, y}^{p+1}(\R \times \T)}^{2} 
\lesssim \|u_{n}\|_{X}^{2} \lesssim
\|u_{n}\|_{L_{x, y}^{p+1}(\R \times \T)}^{p+1}. 
\end{equation*}
This yields that there exists a constant $C_{1} >0$ such that 
$C_{1} \leq \|u_{n}\|_{L_{x, y}^{p+1}(\R \times \T)}$ 
for all $n \in \mathbb{N}$. 
It follows from Lemma \ref{existence-ground-thm3} that
\begin{equation*}
C_{1} \leq \sup_{k \in \mathbb{Z}} 
\left(\|u_{n}\|_{L_{x, y}^{2}(Q_{k})} 
\right)^{p-1} \|u_{n}\|_{X}^{2} 
\lesssim \sup_{k \in \mathbb{Z}} 
\left(\|u_{n}\|_{L_{x, y}^{2}(Q_{k})} 
\right)^{p-1}. 
\end{equation*}
Thus, there exists $\{\tau_{n}\} \subset \mathbb{Z}$, 
subsequence of $u_{n}$ (we still denote it by the same letter)
and $Q_\omega \neq 0$ 
such that $u_{n}(\cdot + \tau_{n}, \cdot)$ 
converges weakly to $Q_{\omega} \neq 0$ 
in $X$ as $n$ tends to infinity. 
We put $v_{n}(\cdot) = u_{n}(\cdot + \tau_{n}, \cdot)$. 
Then, by Brezis and Lieb lemma~\cite{Brezis-Lieb}, 
we have 
\begin{align}
& \lim_{n \to \infty}\left\{ 
\mathcal{I}_{\omega}(v_{n}) - \mathcal{I}_{\omega}(v_{n} - Q_{\omega}) 
- \mathcal{I}_{\omega}(Q_{\omega})
\right\} = 0, \label{eq-3-1}\\
& \lim_{n \to \infty}\left\{ 
\mathcal{N}_{\omega}(v_{n}) - \mathcal{N}_{\omega}(v_{n} - Q_{\omega}) 
- \mathcal{N}_{\omega}(Q_{\omega})
\right\} = 0. \label{eq-4-1}
\end{align}
Suppose that $\mathcal{N}_{\omega}(Q_{\omega}) >0$. 
Since $\mathcal{N}_{\omega}(v_{n}) 
= \mathcal{N}_{\omega}(u_{n})= 0$ for all $n \in \mathbb{N}$, 
we see from \eqref{eq-4-1} that $\mathcal{N}_{\omega}(v_{n} - Q_{\omega}) < 0$ 
for sufficiently large $n \in \mathbb{N}$. 
This yields that 
there exists $t_{n} \in (0, 1)$ such that 
$\mathcal{N}_{\omega}(t_{n}(v_{n} - Q_{\omega})) = 0$. 
It follows that 
\begin{equation*}
m_{\omega} \leq \mathcal{I}_{\omega}(t_{n}(v_{n} - Q_{\omega})) 
< \mathcal{I}_{\omega}(v_{n} - Q_{\omega}). 
\end{equation*} 
On the other hand, since 
$\lim_{n \to \infty} \mathcal{I}_{\omega}(v_{n}) 
=\lim_{n \to \infty} \mathcal{I}_{\omega}(u_{n})= m_{\omega}$, 
we obtain $\mathcal{I}_{\omega}(Q_{\omega}) < 0$, which is a contradiction. 
Therefore, we have $\mathcal{N}_{\omega}(Q_{\omega}) \leq 0$. 
From the lower semicontinuity of $X$ norm, we have 
\begin{equation*}
m_{\omega} \leq \mathcal{I}_{\omega}(Q_{\omega}) 
\leq \liminf_{n \to \infty} 
\mathcal{I}_{\omega}(v_{n}) = m_{\omega}. 
\end{equation*}
This implies that $m_{\omega} = \mathcal{I}_{\omega}(Q_{\omega}) 
= \lim_{n \to \infty} \mathcal{I}_{\omega}(v_{n})$. 
Then, by \eqref{eq-3-1}, we have 
\begin{equation*}
\lim_{n \to \infty} \mathcal{I}_{\omega}(v_{n} - Q_{\omega}) 
= \lim_{n \to \infty} \left\{ 
\mathcal{I}_{\omega}(v_{n} ) -
\mathcal{I}_{\omega}(Q_{\omega}) 
\right\} = 0. 
\end{equation*}
Therefore, $\{v_{n}\}$ converges to $Q_{\omega}$ strongly 
in $X$ as $n$ goes to infinity. 
It follows that $\mathcal{N}_{\omega}(Q_{\omega}) 
= \lim_{n \to \infty} \mathcal{N}_{\omega}(v_{n}) = 0$ and 
$Q_{\omega}$ is a minimizer for $m_{\omega}$. 

Next, we shall show that $Q_{\omega}$ is a ground state to \eqref{sp}. 
Since $Q_{\omega}$ is a minimizer for $m_{\omega}$, 
there exists a Lagrange multiplier $\lambda \in \mathbb{R}$ such that 
$\mathcal{S}_{\omega}^{\prime}(Q_{\omega}) 
= \lambda \mathcal{N}_{\omega}^{\prime}(Q_{\omega})$. 
This yields that 
\begin{equation*}
0 = \mathcal{N}_{\omega}(Q_{\omega}) 
= \langle \mathcal{S}_{\omega}^{\prime}(Q_{\omega}), Q_{\omega} \rangle
= \lambda \langle \mathcal{N}_{\omega}^{\prime}(Q_{\omega}), Q_{\omega} \rangle
= - \lambda (p-1)\|Q_{\omega}\|_{X}^{2}. 
\end{equation*}
Since $Q_{\omega} \neq 0$, we have $\lambda = 0$. 
This yields that $S_{\omega}^{\prime}(Q_{\omega}) = 0$, 
that is, $Q_{\omega}$ is a solution to \eqref{sp}. 
Moreover, we see that any non-trivial solution $u \in X$ 
to \eqref{sp} satisfies 
$\mathcal{N}_{\omega}(u) = 0$. 
It follows from the definition of $m_{\omega}$ 
that $\mathcal{S}_{\omega}(Q_{\omega}) 
\leq \mathcal{S}_{\omega}(u)$. 
Thus, we find that $Q_{\omega}$ is a 
ground state to \eqref{sp}. 

\end{proof}
\subsection{The high frequency case $\omega \gg 1$}
In this subsection, we shall show that $Q_{\omega} \neq R_{\omega}$ for sufficiently large $\omega > 0$. First, we recall that
\begin{equation} \label{R-mini}
m_{\omega, \mathbb{R}} 
=\inf\left\{\mathcal{S}_{\omega, \mathbb{R}}(u) \colon 
u \in H^{1}(\mathbb{R}), \; 
\mathcal{N}_{\omega, \mathbb{R}}(u) = 0 \right\}, 
\end{equation}
with 
\begin{align}
& \mathcal{S}_{\omega, \mathbb{R}}(u) := 
\frac{1}{2} \int_{\mathbb{R}}
\left\{|\partial_{x} u(x)|^{2} + \omega |u(x)|^{2} \right\} dx - 
\frac{1}{p+1} \int_{\mathbb{R}} |u(x)|^{p+1} dx, \label{R-action}\\
&\mathcal{N}_{\omega, \mathbb{R}}(u) := 
\int_{\mathbb{R}}
\left\{|\partial_{x} u(x)|^{2} + \omega |u(x, y)|^{2} \right\} dx - 
\int_{\mathbb{R}} |u(x)|^{p+1} dx, \label{R-nehari}
\end{align}
we have 
\begin{equation} \label{R-vari-chara}
\begin{split}
m_{\omega, \mathbb{R}} 
= \mathcal{S}_{\omega, \mathbb{R}}(R_{\omega})
= \inf\left\{\mathcal{I}_{\mathbb{R}}(u) \colon 
u \in H^{1}(\mathbb{R}), 
\mathcal{N}_{\omega, \mathbb{R}}(u) \leq 0 \right\},
\end{split}
\end{equation}
where 
\begin{equation} \label{I-R}
\mathcal{I}_{\omega, \mathbb{R}}(u) := 
\frac{p-1}{2(p+1)}\left\{\|\partial_{x} u\|_{L_{x}^{2}(\R)}^{2} 
+ \omega \|u\|_{L_{x}^{2}(\R)}^{2} \right\}.
\end{equation}
For simplicity, we denote 
$\mathcal{S}_{1, \mathbb{R}}, 
\mathcal{N}_{1, \mathbb{R}}, 
m_{1, \mathbb{R}}$ and $R_{1}$ 
by $\mathcal{S}_{\mathbb{R}}, 
\mathcal{N}_{\mathbb{R}}, m_{\mathbb{R}}$ and $R$. 
Next, we rescale the ground state by 
\begin{equation} \label{rescale-GroundState}
\widetilde{Q}_{\omega}(x, y) = \omega^{- \frac{1}{p-1}} 
Q_{\omega}(\omega^{-\frac{1}{2}}x, y). 
\end{equation}
Then, we see that $\widetilde{Q}_{\omega}$ satisfies the following: 
\begin{equation} \label{rescale-sp}
- \partial_{xx} \widetilde{Q}_{\omega} + \omega^{-1}
|D_{y}| \widetilde{Q}_{\omega} + \widetilde{Q}_{\omega} - 
|\widetilde{Q}_{\omega}|^{p-1}\widetilde{Q}_{\omega} = 0 \qquad 
\mbox{in $\mathbb{R} \times \mathbb{T}$}. 
\end{equation}
We put 
\begin{equation} \label{rescale-mini}
\widetilde{m}_{\omega} := 
\inf\left\{\widetilde{\mathcal{S}}_{\omega}(u) \colon 
u \in X \setminus \{0\}, \; \widetilde{\mathcal{N}}_{\omega}(u) = 0 \right\}, 
\end{equation}
where 
\begin{equation}
\begin{split}
& \widetilde{\mathcal{S}}_{\omega}(u) 
:= 
\frac{1}{2} \int_{\mathbb{R} \times \mathbb{T}}
\left\{|\partial_{x} u(x, y)|^{2} + \omega^{-1}|D_{y}|u(x, y) \overline{u(x, y)} 
+ |u(x, y)|^{2} \right\} dxdy 
\label{rescale-action}\\
& \hspace{1.4cm} - 
\frac{1}{p+1} \int_{\mathbb{R} \times \mathbb{T}} |u(x, y)|^{p+1} dxdy, 
\end{split}
\end{equation}
\begin{equation}
\begin{split}
& \widetilde{\mathcal{N}}_{\omega}(u) := 
\int_{\mathbb{R} \times \mathbb{T}}
\left\{|\partial_{x} u(x, y)|^{2} + \omega^{-1}|D_{y}|u(x, y) \overline{u(x, y)} 
+ |u(x, y)|^{2} \right\} dxdy
\label{rescale-nehari} \\
& \hspace{1.4cm}- \int_{\mathbb{R} \times \mathbb{T}} |u(x, y)|^{p+1} dxdy. 
\end{split} 
\end{equation}
Then, by a similar argument in Lemma \ref{existence-ground-thm2}, 
we have 
\begin{equation} \label{eq3-chara-ground}
\widetilde{m}_{\omega} = 
\inf\left\{\widetilde{\mathcal{I}}_{\omega}(u) \colon 
u \in X \setminus \{0\}, \; \widetilde{\mathcal{N}}_{\omega}(u) \leq 0 \right\}, 
\end{equation}
where
\begin{equation} \label{positive-function}
\begin{split}
%& 
\widetilde{\mathcal{I}}_{\omega}(u) := 
\frac{p-1}{2(p+1)} \int_{\mathbb{R} \times \mathbb{T}}
\big\{|\partial_{x} u(x, y)|^{2} + \omega^{-1}|D_{y}|u(x, y) \overline{u(x, y)} 
%\\ & \hspace{1.4cm} 
+ |u(x, y)|^{2} \big\} dxdy, 
\end{split}
\end{equation}
In addition, we see that 
\begin{equation} \label{eq-differe-func}
\int_{\mathbb{T}} \mathcal{I}_{\omega, \mathbb{R}}(u) dy \leq 
\widetilde{\mathcal{I}}_{\omega}(u), \qquad 
\int_{\mathbb{T}} \mathcal{N}_{\omega, \mathbb{R}}(u) dy \leq 
\widetilde{\mathcal{N}}_{\omega}(u) \qquad 
\mbox{for all $u \in X$}. 
\end{equation}

From the definition of $\widetilde{S}_{\omega}$, 
$\widetilde{\mathcal{S}}_{\omega}(R) = 2 \pi m_{\mathbb{R}}$. 
Since 
$R_{\omega}(x) = \omega^{\frac{1}{p-1}} R_{1}(\sqrt{\omega} x)$, 
we have 
\begin{equation*} \label{eq-vari-chara-R}
m_{\omega, \mathbb{R}} = \omega^{\frac{p+1}{p-1} - \frac{1}{2}} m_{\mathbb{R}}.
\end{equation*}

We shall show the following: 
\begin{proposition} \label{vari-omega-large}
There exists a sufficiently large $\omega_{L}>0$ such that 
$Q_{\omega} \neq R_{\omega}$ for all $\omega > \omega_{L}$
\end{proposition}
\begin{proof}[Proof of Proposition \ref{vari-omega-large}]
Since $m_{\omega} 
= \omega^{\frac{p+1}{p-1} - \frac{1}{2}} 
\widetilde{m}_{\omega}$, 
$m_{\omega, \mathbb{R}} = \omega^{\frac{p+1}{p-1} - \frac{1}{2}} m_{\mathbb{R}}$ and 
$\widetilde{S}_{\omega}(\widetilde{Q}_{\omega}) 
= m_{\omega}$ and 
$\widetilde{S}_{\omega} (R) = 2 \pi m_{\mathbb{R}}$, 
it is enough to show that 
$m_{\omega} < 2 \pi m_{\mathbb{R}}$ for sufficiently large $\omega >0$. The idea is to rescale $R$ with a $y$-dependent scaling that normalizes the $\dot{H}^1_xL^2_y$ and $L^{p+1}_{x,y}$ norms to $2\pi$ times the $\dot{H}^1_x$ and $L^{p+1}_{x}$ norms. Let $\rho \in C_{0}^{\infty}(\mathbb{T})$ 
be a positive function such that 
\begin{equation} \label{eq0-chara-ground}
\int_{\mathbb{T}} \rho^{\frac{p+3}{2(p-1)}}(y) dy = 2 \pi, \qquad 
\rho \not\equiv 1. 
\end{equation}
We put 
\begin{equation*}
\psi(x, y) := \rho^{\frac{1}{p-1}}(y) R
(\sqrt{\rho(y)} x). 
\end{equation*}
Then, we have 
\begin{equation} \label{eq1-chara-ground}
\begin{split}
& \|\partial_{x} \psi(\cdot, y)\|_{L_{x}^{2}(\mathbb{R})}^{2} 
= 
\rho(y)^{\frac{p+3}{2(p-1)}} 
\|\partial_{x} R\|_{L_{x}^{2}(\mathbb{R})}^{2}, \qquad 
\|\psi(\cdot, y)\|_{L_{x}^{2}(\mathbb{R})}^{2} = 
\rho(y)^{\frac{5-p}{2(p-1)}} 
\|R\|_{L_{x}^{2}(\mathbb{R})}^{2}, \qquad \\
& \|\psi(\cdot, y)\|_{L_{x}^{p+1}(\mathbb{R})}^{p+1} = 
\rho(y)^{\frac{p+3}{2(p-1)}} 
\|R\|_{L_{x}^{p+1}(\mathbb{R})}^{p+1}. 
\end{split}
\end{equation}
By the H\"{o}lder inequality and \eqref{eq0-chara-ground}, we have 
\begin{equation*}
\int_{\mathbb{T}} \rho^{\frac{5-p}{2(p-1)}}(y) dy 
< \left(\int_{\mathbb{T}} \rho^{\frac{p+3}{2(p-1)}}(y) dy\right)^{\frac{5-p}{p+3}}
\left(\int_{\mathbb{T}} 1 dy \right)^{\frac{2p-2}{p+3}} 
= 2\pi. 
\end{equation*}
We note that the above inequality is strict since 
$\rho \not\equiv 1$. 
We put 
\begin{equation} \label{eq3-1-chara-ground}
\delta = 2 \pi - \int_{\mathbb{T}} \rho^{\frac{5-p}{2(p-1)}}(y) dy \; 
(>0). 
\end{equation}
By \eqref{eq1-chara-ground}, 
\eqref{eq0-chara-ground}, 
\eqref{eq3-1-chara-ground} and 
$\mathcal{N}_{\mathbb{R}}(R) = 0$, we have
\begin{equation} \label{eq2-chara-ground}
\begin{split}
\widetilde{\mathcal{N}}_{\omega}(\psi)
& = 2\pi \|\partial_{x} R\|_{L_{x}^{2}(\mathbb{R})}^{2} 
+ \int_{\T} \rho(y)^{\frac{5-p}{2(p-1)}} dy \|R\|_{L_{x}^{2}(\mathbb{R})}^{2} 
+ \omega^{-1} \||D_{y}|^{\frac{1}{2}} \psi\|_{L_{x, y}^{2}(\mathbb{R} 
\times \mathbb{T})}^{2} \\
& \quad -2 \pi \|R\|_{L_{x}^{p+1}(\mathbb{R})}^{p+1} \\
& = 2 \pi \|\partial_{x} R\|_{L_{x}^{2}(\mathbb{R})}^{2} 
+ 2 \pi \|R\|_{L_{x}^{2}(\mathbb{R})}^{2} 
- 2\pi \|R\|_{L_{x}^{p+1}(\mathbb{R})}^{p+1}\\ 
& \quad 
- \delta \|R\|_{L_{x}^{2}(\mathbb{R})}^{2}
+ \omega^{-1} \||D_{y}^{\frac{1}{2}}| 
\psi\|_{L_{x, y}^{2}(\mathbb{R} \times 
\mathbb{T})}^{2} \\
& = - \delta \|R\|_{L_{x}^{2}(\mathbb{R})}^{2}
+ \omega^{-1} \||D_{y}^{\frac{1}{2}}| \psi
\|_{L_{x, y}^{2}(\mathbb{R} 
\times \mathbb{T})}^{2}. 
\end{split}
\end{equation} 
We take $\omega >0$ sufficiently large so that 
\begin{equation}\label{eq4-chara-ground}
- \delta \|R\|_{L_{x}^{2}(\R)}^{2}
+ \omega^{-1} \||D_{y}^{\frac{1}{2}}| \psi\|_{L_{x, y}^{2}
(\R \times \T)}^{2}
< - \frac{\delta}{2} \|R\|_{L_{x}^{2}(\R)}^{2} \leq 0. 
\end{equation}
This together with \eqref{eq2-chara-ground} yields 
that $\widetilde{\mathcal{N}}_{\omega}(\psi) < 0$. 
Then, from \eqref{eq3-chara-ground}, 
\eqref{R-vari-chara}, \eqref{I-R} and 
\eqref{eq4-chara-ground}, we have 
\begin{equation*} %\label{eq1-vari-chara}
\begin{split}
m_{\omega}
& \leq \widetilde{\mathcal{I}}_{\omega}(\psi) \\
& = \frac{p-1}{2(p+1)} \left\{ 
2\pi\|\partial_{x} R\|_{L_{x}^{2}(\R)}^{2} + 
2\pi\|R\|_{L_{x}^{2}(\R)}^{2} 
- \delta \|R\|_{L_{x}^{2}(\R)}^{2} + 
\omega^{-1}\||D_{y}|
^{\frac{1}{2}}\psi\|_{L_{x, y}^{2}(\R \times \T)}^{2} 
\right\} \\
& = 2 \pi m_{\mathbb{R}}
+ \frac{p-1}{2(p+1)} 
\left\{ 
- \delta \|R\|_{L_{x}^{2}(\R)}^{2} + 
\omega^{-1}\||D_{y}|^{\frac{1}{2}} 
\psi \|_{L_{x, y}^{2}(\R \times \T)}^{2} 
\right\} \\
& < 2 \pi m_{\mathbb{R}}. 
\end{split}
\end{equation*}
This completes the proof. 
\end{proof}
\subsection{The low frequency case $0 < \omega \ll 1$}
In this subsection, we shall show the following:
\begin{proposition} \label{vari-omega-small}
There exists a sufficiently small $\omega_{S}>0$ such that 
$Q_{\omega} = R_{\omega}$ for all $0 < \omega < \omega_{S}$. 
\end{proposition}
To prove Proposition \ref{vari-omega-small}, 
we prepare several preparation. We will use the same notations an in Section 4.2. We first show the following: 
\begin{lemma} 
\label{thm-2-10}
The minimum value of \eqref{mini-pro} 
and its corresponding ground state 
satisfy 
\begin{align}
& \lim_{\omega \to 0} 
\widetilde{m}_{\omega} = 
2 \pi m_{\mathbb{R}}, \label{eq-01} \\
& \lim_{\omega \to 0} \frac{1}{\omega} \int_{\mathbb{R} \times 
\mathbb{T}} |D_{y}| \widetilde{Q}_{\omega}(x, y) 
\overline{\widetilde{Q}_{\omega}(x, y)} dxdy = 0. 
\label{eq-02}
\end{align}
\end{lemma}
\begin{proof}
We take $\{\omega_{j}\} \subset \mathbb{R}_{+}$ with 
$\lim_{j \to \infty} \omega_{j} = 0$. 
Since $\widetilde{\mathcal{N}}_{\omega}(R) = 0$, it follows from 
the definition of $m_{\omega}$ that 
\begin{equation} \label{eq-2}
\widetilde{m}_{\omega_{j}} \leq \widetilde{S}_{\omega_{j}}(R) = 2 \pi 
m_{\mathbb{R}} \qquad 
\mbox{for all $j \in \mathbb{N}$}. 
\end{equation}
We shall show that 
\begin{equation} \label{eq4-1-chara-ground}
\lim_{j \to \infty} 
\int_{\mathbb{R} \times 
\mathbb{T}} |D_{y}| \widetilde{Q}_{\omega_{j}}(x, y) 
\overline{\widetilde{Q}_{\omega_{j}}(x, y)}dxdy = 0. 
\end{equation}
Arguing by contradiction and 
suppose that there exist a constant 
$\varepsilon_{0}>0$ and a subsequence of 
$\{\omega_{j}\}$ 
(we still denote it by the same letter) such that 
\begin{equation*}
\int_{\mathbb{R} \times 
\mathbb{T}} |D_{y}| \widetilde{Q}_{\omega_{j}}(x, y) 
\overline{\widetilde{Q}}_{\omega_{j}}(x, y)dxdy \geq \varepsilon_{0} >0
\qquad \mbox{for all $j \in \mathbb{N}$}. 
\end{equation*}
This implies that 
\begin{equation*}
\lim_{j \to \infty} 
\frac{1}{\omega_{j}}
\int_{\mathbb{R} \times 
\mathbb{T}} |D_{y}| \widetilde{Q}_{\omega_{j}}(x, y) 
\overline{\widetilde{Q}_{\omega_{j}}(x, y)} dxdy = \infty. 
\end{equation*}
On the other hand, it follows from \eqref{eq-2} that 
for sufficiently large $j \in \mathbb{N}$, we have 
\begin{equation} \label{eq2-vari-chara}
\begin{split}
& 2 \pi m_{\mathbb{R}} +1 \\
& \geq \widetilde{m}_{\omega_{j}} + 1\\
& \geq \widetilde{S}_{\omega_{j}}(\widetilde{Q}_{\omega_{j}}) 
- \frac{1}{p+1} \widetilde{N}_{\omega_{j}} (\widetilde{Q}_{\omega_{j}}) \\
& = \left(\frac{1}{2} - \frac{1}{p+1}\right) 
\left(
\int_{\mathbb{R} \times \mathbb{T}} 
\left\{|\partial_{x} \widetilde{Q}_{\omega_{j}}(x, y)|^{2} 
+ \frac{1}{\omega_{j}} |D_{y}| \widetilde{Q}_{\omega_{j}}(x, y) 
\overline{\widetilde{Q}_{\omega_{j}}(x, y)}
+ |\widetilde{Q}_{\omega_{j}}|^{2}
\right\} dxdy \right) \\
& \to \infty \qquad \mbox{as $j \to \infty$}, 
\end{split}
\end{equation}
which is a contradiction. 
Therefore, \eqref{eq4-1-chara-ground} holds. 
%In addition, we see from 
%\eqref{eq2-vari-chara} that 
%$\sup_{j \in \mathbb{N}} 
%\frac{1}{\omega_{j}} 
%\int_{\R \times \T} |D_{y}| \widetilde{Q}_{\omega_{j}}(x, y) 
%\overline{\widetilde{Q}_{\omega_{j}}(x, y)} 
%dxdy
%< \infty$. 

We first consider the case of $1 < p \leq 3$. 
By the H\"{o}lder and Gagliardo-Nirenberg inequalities 
and \eqref{eq4-1-chara-ground}, 
we have 
\begin{equation} \label{eq-3}
\begin{split}
& \quad 
\int_{\mathbb{R} \times \mathbb{T}} 
|\widetilde{Q}_{\omega_{j}}|^{p} ||D_{y}|^{\frac{1}{2}} \widetilde{Q}_{\omega_{j}}| dxdy \\
& \leq \|\widetilde{Q}_{\omega_{j}}\|_{L_{x, y}^{2p}(\R \times \T)}^{p} 
\||D_{y}|^{\frac{1}{2}} \widetilde{Q}_{\omega_{j}}\|_{L_{x, y}^{2}
(\R \times \T)} \\
& \lesssim \|\widetilde{Q}_{\omega_{j}}
\|_{L_{x, y}^{2}(\R \times \T)}^{\frac{3-p}{2}}
\||D_{y}|^{\frac{1}{2}} \widetilde{Q}_{\omega_{j}}
\|_{L_{x, y}^{2}(\R \times \T)}^{p}
\|\partial_{x} \widetilde{Q}_{\omega_{j}}
\|_{L_{x, y}^{2}(\R \times \T)}^{\frac{p-1}{2}} \\
& \to 0 \qquad \mbox{as $j \to \infty$}. 
\end{split}
\end{equation}
Next, we consider the case of $3 < p < 5$. 
It follows from \eqref{eq2-vari-chara} that 
$\sup_{j \in \mathbb{N}}\|\widetilde{Q}_{\omega_{j}}\|_{X} \leq M$ 
for some $M >0$. 
Then, by Theorem \ref{elliptic-regularity} below, 
there exists $C = C(M) > 0$ such that 
$\sup_{j \in \mathbb{N}}
\|\widetilde{Q}_{\omega_{j}}\|_{L_{x, y}^{\infty}} \leq C$. 
This yields that 
\begin{equation} \label{eq-3-30}
\begin{split}
& \quad \int_{\mathbb{R} \times \mathbb{T}} 
|\widetilde{Q}_{\omega_{j}}|^{p} ||D_{y}|^{\frac{1}{2}} \widetilde{Q}_{\omega_{j}}| dxdy \\
& \lesssim 
\|\widetilde{Q}_{\omega_{j}}\|_{L_{x, y}^{6}(\R \times \T)}^{3} 
\||D_{y}|^{\frac{1}{2}} \widetilde{Q}_{\omega_{j}}
\|_{L_{x, y}^{2}(\R \times \T)} \\
& \lesssim 
\|\widetilde{Q}_{\omega_{j}}\|_{X}^{3} 
\||D_{y}|^{\frac{1}{2}} \widetilde{Q}_{\omega_{j}}
\|_{L_{x, y}^{2}(\R \times \T)}
\to 0 \qquad \mbox{as $j \to \infty$}. 
\end{split}
\end{equation}
Multiplying \eqref{rescale-sp} by 
$|D_{y}|^{\frac{1}{2}}\widetilde{Q}_{\omega_{j}}$ and 
integrating the resulting equation, we have, by \eqref{eq-3} and \eqref{eq-3-30}, that 
\begin{equation*} 
\int_{\mathbb{R} \times \mathbb{T}} 
\left\{
||D_{y}|^{\frac{1}{4}} \partial_{x} \widetilde{Q}_{\omega_{j}}|^{2} 
+ ||D_{y}|^{\frac{1}{4}} \widetilde{Q}_{\omega_{j}}|^{2} 
+ \omega_{j}^{-1}||D_{y}|^{\frac{3}{4}}| \widetilde{Q}_{\omega_{j}}|^{2}
\right\} dx \to 0 \qquad \mbox{as $j \to \infty$}. 
\end{equation*}

This implies that
\begin{equation}
\label{eq-4}
\omega_{j}^{-1} \int_{\mathbb{R} \times \mathbb{T}} \left| | D_{y}|^{\frac{3}{4}} \widetilde{Q}_{\omega_{j}}\right|^{2} d x \rightarrow 0 \quad \text { as } j \rightarrow \infty.
\end{equation}
Thus,
\begin{align}
\label{eq-4.1}
\omega_{j}^{-1} \int_{\mathbb{R} \times \mathbb{T}} | D_{y}| \widetilde{Q}_{\omega_{j}} \overline{\widetilde{Q}}_{\omega_{j}} d x d y &=\omega_{j}^{-1} \int_{\mathbb{R} \times \mathbb{T}} \left| | D_{y}|^{\frac{1}{2}} \widetilde{Q}_{\omega_{j}}\right|^{2} d x d y \nonumber\\
& \leq \omega_{j}^{-1} \int_{\mathbb{R} \times \mathbb{T}} \left| | D_{y}|^{\frac{3}{4}} \widetilde{Q}_{\omega_{j}}\right|^{2} d x d y \rightarrow 0 \quad \text { as } j \rightarrow \infty.
\end{align}
This finishes the proof of \eqref{eq-02}.

In order to prove \eqref{eq-01}, 
multiplying \eqref{rescale-sp} by 
$\widetilde{Q}_{\omega_{j}}$ and 
integrating the resulting equation on $\mathbb{R}$, we have 
\begin{equation}
\label{def:Vj}
\mathcal{N}_{\mathbb{R}}\left(\widetilde{Q}_{\omega_{j}}(\cdot, y)\right)=V_{j}(y)
\end{equation}
where
\[
V_{j}(y):= - \omega_{j}^{-1} \int_{\mathbb{R}} | D_{y}| \widetilde{Q}_{\omega_{j}} \overline{\widetilde{Q}}_{\omega_{j}} d x \qquad (y \in \mathbb{T})
\]
Then, we see from Remark 5.1 that $V_{j} \in C_{y}^{1}(\mathbb{T}) .$ Therefore, from \eqref{eq-4.1}, we see that 
\[
\lim _{j \rightarrow \infty} V_{j}(y) 
= 0 \qquad 
\mbox{for all $y \in \mathbb{T}$}.
\]
This implies that 
$$ \lim _{j \rightarrow \infty} \mathcal{N}_{\mathbb{R}}\left(\widetilde{Q}_{\omega_{j}}(\cdot, y)\right)=0 . $$ 
Then, there exists $t_{j}(y)>0$ with $\lim _{j \rightarrow \infty} t_{j}(y)=1$ such that $\mathcal{N}_{\mathbb{R}}\left(t_{j}(y) \widetilde{Q}_{\omega_{j}}(\cdot, y)\right)=0 .$ From \eqref{R-vari-chara}, we have
\[
m_{\mathbb{R}} \leq \mathcal{I}_{\mathbb{R}}\left(t_{j}(y) \widetilde{Q}_{\omega_{j}}(\cdot, y)\right)=t_{j}(y)^{2} \mathcal{I}_{\mathbb{R}}\left(\widetilde{Q}_{\omega_{j}}(\cdot, y)\right)
\]
%Next, we recall the following result
%\begin{lemma}[\cite{kpv}]
%We have the following estimate:
%$$\||D_x|^s (fg) - g|D_x|^s f - f |D_x|^s g\|_{L^p}\lesssim \||D_x|^{s_1} f\|_{L^q} \||D_x|^{s_2} f\|_{L^r},$$
%where
%$$\frac 1 p=\frac 1q+ \frac 1r, \hbox{ } 1<p,q,r<\infty, \hbox{ } 1>s=s_1+s_2>0, s_i\geq 0.$$ 
%Moreover, if $s_2=0$ then $r=\infty$ is allowed.
%\end{lemma}
%Combining the Cauchy-Schwarz estimate with this lemma, we get
%\begin{align*}
%\int_{\R \times \T} \left( |D_{y}|^{\frac{1}{4}} \left( \left||D_{y}|^{\frac{1}{2}}| \widetilde{Q}_{\omega_{j}}\right|^{2} \right) - 2 \left||D_{y}|^{\frac{3}{4}}| \widetilde{Q}_{\omega_{j}}\right|^{2} \right) dydx
%\end{align*}

Integrating the above inequality on $y \in \mathbb{T}$, 
by \eqref{eq-differe-func} and \eqref{eq-2}, we obtain 
\begin{equation*} %\label{eq-5}
\begin{split}
& 2 \pi m_{\mathbb{R}} \leq 
\liminf_{j \to \infty} 
\int_{\mathbb{T}} t_{j}(y)^{2} \mathcal{I}_{\mathbb{R}}
(\widetilde{Q}_{\omega_{j}}(\cdot, y)) dy 
\leq \liminf_{j \to \infty} 
\int_{\mathbb{T}} \mathcal{I}_{\mathbb{R}}
(\widetilde{Q}_{\omega_{j}}(\cdot, y)) dy \\
& \leq \limsup_{j \to \infty}\widetilde{\mathcal{I}}_{\omega_{j}}
(\widetilde{Q}_{\omega_{j}}) 
= \limsup_{j \to \infty} \widetilde{m}_{\omega_{j}}
\leq 2 \pi m_{\mathbb{R}}. 
\end{split}
\end{equation*}
Thus, we find that \eqref{eq-01} holds. 

%\par 
%Finally, 
%From \eqref{eq-5}, we see that 
%\begin{equation*}
%\lim_{j \to \infty} 
%\int_{\mathbb{T}} \mathcal{I}_{\mathbb{R}}
%(\widetilde{Q}_{\omega_{j}}(\cdot, y)) dy 
%= \lim_{j \to \infty} \widetilde{m}_{\omega_{j}}
%= 2\pi m_{\R}. 
%\end{equation*}
%It follows that 
%\begin{equation*}
%\lim_{j \to \infty} 
%\left\{\widetilde{I}_{\omega_{j}}
%(\widetilde{Q}_{\omega_{j}}) - \int_{\mathbb{T}} \mathcal{I}_{\mathbb{R}}(\widetilde{Q}_{\omega_{j}}(\cdot, y)) dy\right\}
%= \lim_{j \to \infty} 
%\left\{ m_{\omega_{j}} - \int_{\mathbb{T}} \mathcal{I}_{\mathbb{R}}(\widetilde{Q}_{\omega_{j}}(\cdot, y)) dy\right\}
%= 0. 
%\end{equation*}
%Therefore, we have 
%\begin{equation*}
%\lim_{j \to \infty} \omega_{j}^{-1} 
%\int_{\mathbb{R} \times \mathbb{T}} |D_{y}|Q_{\omega_{j}}(x, y)
%\overline{Q_{\omega_{j}}(x, y)} dxdy 
%= \lim_{j \to \infty}
%\left\{\widetilde{I}_{\omega_{j}}
%(\widetilde{Q}_{\omega_{j}}) - \int_{\mathbb{T}} \mathcal{I}_{\mathbb{R}}(\widetilde{Q}_{\omega_{j}}(\cdot, y)) dy \right\} 
%= 0. 
%\end{equation*}
This completes the proof. 
\end{proof}
Using Lemma \ref{thm-2-10}, we shall show the 
following convergence:
\begin{proposition} \label{thm-3}
Let $1 < p < 5$ and $\{\omega_{j}\} \subset \R_{+}$ 
with $\lim_{j \to \infty} \omega_{j} = 0$. 
There exists a sequence 
$\{\tau_{j}\} \subset \mathbb{R}$ such that 
\begin{equation*}
\widetilde{Q}_{\omega_{j}}(\cdot + \tau_{j}, \cdot) \to R \qquad 
\mbox{in $X$ as $j \to 0$}
\end{equation*}
\end{proposition}
\begin{proof}
We take $\{\omega_{j}\} \subset \mathbb{R}_{+}$ 
with $\lim_{j \to \infty} \omega_{j} = 0$. 
From \eqref{eq-2}, we have 
\begin{equation*}
\begin{split}
2 \pi m_{\mathbb{R}}
& \geq \widetilde{\mathcal{S}}_{\omega_{j}}(\widetilde{Q}_{\omega_{j}}) 
- \frac{1}{p+1} \widetilde{\mathcal{N}}_{\omega_{j}}(\widetilde{Q}_{\omega_{j}}) \\
& = \left(\frac{1}{2} - \frac{1}{p+1}\right)
\int_{\mathbb{R} \times \mathbb{T}} 
\Big\{ |\partial_{x} \widetilde{Q}_{\omega_{j}}(x, y)|^{2}
+ \omega_{j}^{-1}|D_{y}|\widetilde{Q}_{\omega_{j}}(x, y) \overline{\widetilde{Q}_{\omega_{j}}(x, y)}
\\
& \qquad \qquad \qquad \qquad \qquad + |\widetilde{Q}_{\omega_{j}}(x, y)|^{2}
\Big\} dxdy. 
\end{split}
\end{equation*}
Thus, $\{\widetilde{Q}_{\omega_{j}}\}$ is bounded in $X$. 

Then, by the similar argument in the proof of Theorem \ref{thm-exist-ground}, 
we see that 
there exist a subsequence of $\{\widetilde{Q}_{\omega_{j}}\}$ 
(we still denote it by the same letter), $\{\tau_{j}\} \subset \mathbb{R}$ and 
$w \in X \setminus \{0\}$ such that 
$\lim_{j \to \infty} \widetilde{Q}_{\omega_{j}}(\cdot + \tau_{j}, \cdot) 
= w$ weakly in $X$. 
From the weak lower semicontinuity 
and \eqref{eq4-1-chara-ground}, we have 
\begin{equation*}
\begin{split}
\int_{\mathbb{R} \times \mathbb{T}}
|D_{y}|w(x, y) \overline{w(x, y)} dxdy 
\leq
\lim_{j \to \infty} 
\int_{\mathbb{R} \times \mathbb{T}}
|D_{y}|\widetilde{Q}_{\omega_{j}}(x, y) \overline{\widetilde{Q}_{\omega_{j}}(x, y)}dxdy 
= 0. 
\end{split}
\end{equation*}
Especially, we have $|D_{y}|^{1/2} w = 0$ in 
$L_{x, y}^{2}(\mathbb{R} \times 
\mathbb{T})$. 
Passing to the weak limit, we see that 
for all $y \in \mathbb{T}$, $w$ satisfies 
\begin{equation*}
- \partial_{xx} w + w - |w|^{p-1}w = 0, \qquad x \in \mathbb{R}
\end{equation*}
in the distribution sense.
Since the nontrivial solution of the above equation is unique, we see that 
$w = R(\cdot + \tau_{1})$. 
Choosing a subsequence, we may assume 
that 
$\lim_{j \to \infty}\widetilde{Q}_{\omega_{j}}
(\cdot + \tau_{j}, \cdot) = R$ weakly in 
$X$. 
Then, it follows %from the Brezis-Lieb Lemma~\cite{Brezis-Lieb} 
that 
\begin{align}
%& 
\widetilde{\mathcal{I}}_{\omega_{j}}(\widetilde{Q}_{\omega_{j}}(\cdot + \tau_{j}, \cdot))
- \widetilde{\mathcal{I}}_{\omega_{j}}
(\widetilde{Q}_{\omega_{j}}(\cdot + \tau_{j}, \cdot) - R)
- \widetilde{\mathcal{I}}_{\omega_{j}}(R) \to 0 \qquad \mbox{as $j \to \infty$}, 
\label{eq-6} 
%\\
%& \widetilde{\mathcal{N}}_{\omega_{j}}(\widetilde{Q}_{\omega_{j}}(\cdot + \tau_{j}, \cdot))
%- \widetilde{\mathcal{N}}_{\omega_{j}}
%(\widetilde{Q}_{\omega_{j}}(\cdot + \tau_{j}, \cdot) - R)
%- \widetilde{\mathcal{N}}_{\omega_{j}}(R) \to 0 \qquad \mbox{as $j \to \infty$}. 
%\label{eq-7} 
\end{align}
From \eqref{eq-01}, one has 
\[
\lim_{j \to \infty}\widetilde{\mathcal{I}}_{\omega_{j}}
(\widetilde{Q}_{\omega_{j}}(\cdot + \tau_{j}, \cdot)) = 
\lim_{j \to \infty} 
\widetilde{m}_{\omega_{j}} = 2 \pi m_{\mathbb{R}}. 
\]
This together with $\widetilde{\mathcal{I}}_{\omega_{j}}(R) = 2\pi m_{\mathbb{R}}$ and \eqref{eq-6}
yields that $\lim_{j \to \infty} \widetilde{\mathcal{I}}_{\omega_{j}}
(\widetilde{Q}_{\omega_{j}}(\cdot + \tau_{j}, \cdot) - R) = 0$. 
Thus, we see that 
$\widetilde{Q}_{\omega_{j}}(\cdot + \tau_{j}, \cdot)$ converges $R$
strongly in $H^{1}(\mathbb{R})$ as $j$ goes to infinity. 
This completes the proof. 
\end{proof}

In what follows, we put 
$\widehat{Q}_{\omega_{j}}(\cdot, \cdot) 
= \widetilde{Q}_{\omega_{j}} (\cdot + \tau_{j}, 
\cdot)$, 
where $\tau_{j} \in \R$
is given by Lemma \ref{thm-3}. 
Since $\widehat{Q}_{\omega_{j}} 
\in C^{2}(\R \times \T)$ 
(see Remark \ref{reg-Q} below) and 
$R$ is positive, 
it follows from Proposition \ref{thm-3} 
that $\widehat{Q}_{\omega_{j}}(x, y)$ 
is positive for all $(x, y) \in 
\R \times \T$. 

We now prove the following fractional calculus estimate: 
%(see \cite[Proposition 3.1]{Chirist-Weinstein} 
%and also \cite[Lemma 7.7]{Linares-Ponce}). 
\begin{lemma}
\label{A1}
Let $0 < \alpha < 1$ and $1 < p < 5$. Then, there exists a constant $C>0$ such that
%, for any $u\in L^{\infty}(\mathbb{R} \times \T) \cap L^2\dot{H}^{\alpha}\left(\R\times\mathbb{T}\right)$, we have
\begin{equation}
\begin{split}
& \quad 
\||D_{y}|^{\alpha} \left( |\widehat{Q}_{\omega_{j}} |^{p-1}\widehat{Q}_{\omega_{j}} -pR^{p-1}\widehat{Q}_{\omega_{j}} \right)
\|_{L_{x, y}^{2}(\mathbb{R} \times \T)} 
\\
& 
\leq 
\begin{cases}
C 
\|\widehat{Q}_{\omega_{j}} -R
\|^{p-1}_{L_{x, y}^\infty(\mathbb{R} \times \T)} 
\||D_{y}|^{\alpha} \widehat{Q}_{\omega_{j}} 
\|_{L_{x, y}^{2}(\mathbb{R} \times \T)} 
& %\qquad \\ & \hspace{-5.5cm}
\mbox{for $1 < p \leq 2$}, \\
C (\|\widehat{Q}_{\omega_{j}}
\|_{L_{x, y}^\infty(\mathbb{R} \times \T)}^{p-2} 
+ \|R
\|_{L_{x}^\infty(\mathbb{R} \times \T)}^{p-2})
\|\widehat{Q}_{\omega_{j}} -R\|
_{L_{x, y}^\infty(\mathbb{R} \times \T)} 
\||D_{y}|^{\alpha} \widehat{Q}_{\omega_{j}} 
\|_{L_{x, y}^{2}(\mathbb{R} \times \T)}
& %\\ & \hspace{-5.5cm} 
\mbox{for $2 < p < 5$}. 
\end{cases}
\end{split}
\end{equation}
\end{lemma}
\begin{proof}
We recall the following representation of the $\dot{H}^s(\T)$ norm 
(see \cite[Proposition 1.2]{BenOh} for more details)%\footnote{\begin{remark}
%In \cite[Proposition 1.4]{BenOh}, the lower and upper bounds are given by showing that
%$$B(n)=|n|^{-2 \alpha} \int_{-\pi}^\pi \frac{\left|e^{i y n}-1\right|^{2}}{|y|^{1+2 \alpha }} d y.$$
%is bounded both from above and below uniformly in $n \in \mathbb{Z}^{d} \backslash\{0\} $ which obviously do not depend on $x$.
%\end{remark}}
\begin{align*}
% &\frac{1}{\widehat{C}} \int_{\mathbb{T}} \int_{-\pi}^\pi \frac{|\widehat{Q}_{\omega_{j}} (x, y+z)-\widehat{Q}_{\omega_{j}} (x,z)|^{2}}{|y|^{1+2 \alpha}} d y d z \\
\|\widehat{Q}_{\omega_{j}}(x, \cdot) 
\|_{\dot{H}^{\alpha}\left(\mathbb{T}\right)}^{2} \sim \int_{\mathbb{T}} \int_{- 
\pi}^{\pi} \frac{|\widehat{Q}_{\omega_{j}} (x,y+z)-\widehat{Q}_{\omega_{j}} (x,z)|^{2}}{|y|^{1+2 \alpha}} d y d z,
\end{align*}
Here, the implicit constants depend only on 
$\alpha$. 
Then, we have 
\begin{equation} \label{eq-diff1}
\begin{split}
&\quad \left\| |D_{y}|^{\alpha} \left( |\widehat{Q}_{\omega_{j}}|^{p-1}\widehat{Q}_{\omega_{j}}(x)-pR^{p-1}(x)\widehat{Q}_{\omega_{j}}(x)\right)\right\|_{L_{y}^{2}( \T)}^{2} \\
& \sim \int_{\mathbb{T}} 
\int_{-\pi}^{\pi} \frac{\left|\left( |\widehat{Q}_{\omega_{j}}|^{p-1}\widehat{Q}_{\omega_{j}}-pR^{p-1}(x)\widehat{Q}_{\omega_{j}}\right)(x,y+z)-\left( |\widehat{Q}_{\omega_{j}}|^{p-1}\widehat{Q}_{\omega_{j}} - 
pR^{p-1}(x) 
\widehat{Q}_{\omega_{j}}\right)(x,z)\right|^{2}}{|y|^{1+2 \alpha}} d y d z\\
& = \int_{\mathbb{T}} \int_{-\pi}^{\pi} 
\frac{\left| (|\widehat{Q}_{\omega_{j}}|^{p-1}\widehat{Q}_{\omega_{j}})(x,y+z)-(|\widehat{Q}_{\omega_{j}}|^{p-1}\widehat{Q}_{\omega_{j}})(x,z) -pR^{p-1}(x)\left( \widehat{Q}_{\omega_{j}}(x,y+z)-\widehat{Q}_{\omega_{j}}(x,z)\right)\right|^{2}}{|y|^{1+2 \alpha}} d y d z\\
& = \int_{\mathbb{T}} \int_{- 
\pi}^{\pi} \frac{\left| \left( \int^1_0 p(s\widehat{Q}_{\omega_{j}}(x,y+z)+(1-s) 
\widehat{Q}_{\omega_{j}}(x,z))^{p-1}ds -pR^{p-1}(x)\right)\left( \widehat{Q}_{\omega_{j}}(x,y+z)-\widehat{Q}_{\omega_{j}}(x,z)\right)\right|^{2}}{|y|^{1+2 \alpha}} d y d z
\end{split}
\end{equation}
We consider the case of $1 < p \leq 2$. 
Using \eqref{estim:p-1}, we obtain 
\begin{equation} \label{eq-diff2}
\begin{split}
&\quad \left\| |D_{y}|^{\alpha} \left( |\widehat{Q}_{\omega_{j}}|^{p-1}\widehat{Q}_{\omega_{j}}(x)-pR^{p-1}(x)\widehat{Q}_{\omega_{j}}(x)\right)\right\|_{L_{y}^{2}( \T)}^{2} \\
&\lesssim \int_{\mathbb{T}} 
\int_{- 
\pi}^{\pi} 
\frac{ \int^1_0 p^{2}|s\widehat{Q}_{\omega_{j}}(x,y+z)+(1-s)\widehat{Q}_{\omega_{j}}(x,z)-R(x)|^{2p-2}ds
\left| \widehat{Q}_{\omega_{j}}(x,y+z)-\widehat{Q}_{\omega_{j}}(x,z)\right|^{2}}{|y|^{1+2 \alpha}} dy dz\\
& \lesssim \left\| \widehat{Q}_{\omega_{j}}
(x, \cdot) -R(x)\right\|^{2p-2}_{L_{y}^\infty(\T)} \int_{\mathbb{T}} \int_{- 
\pi}^{\pi} 
\frac{\left| \widehat{Q}_{\omega_{j}}(x,y+z)-\widehat{Q}_{\omega_{j}}(x,z)\right|^{2}}{|y|^{1+2 \alpha}} d y d z\\
&\lesssim \left\| \widehat{Q}_{\omega_{j}}(x, 
\cdot) -R(x)\right\|^{2p-2}_{L_{y}^\infty(\T)} \left\| |D_{y}|^{\alpha}\widehat{Q}_{\omega_{j}}(x, \cdot) \right\|_{L_{y}^{2}(\T)}^{2}.
\end{split}
\end{equation}
Next, we consider the case of $2 < p < 5$. 
It follows that 
\begin{equation*}
\begin{split}
& \quad 
\biggl|
\int^1_0 |s\widehat{Q}_{\omega_{j}}(x,y+z)+(1-s)\widehat{Q}_{\omega_{j}}(x,z)|^{p-1}ds -R^{p-1}(x)
\biggl| \\
& \leq 
\biggl|
\int_{0}^{1} \int_{0}^{1} 
\frac{d}{d t}
\left(t (s\widehat{Q}_{\omega_{j}}(x,y+z)+(1-s)\widehat{Q}_{\omega_{j}}(x,z))
+ (1-t) R\right)^{p-1} dtds
\biggl| \\
& \leq (p-1) 
\int_{0}^{1} \int_{0}^{1} 
\biggl|t (s\widehat{Q}_{\omega_{j}}(x,y+z)+(1-s)\widehat{Q}_{\omega_{j}}(x,z))
+ (1-t) R\biggl|^{p-2} \times \\
& \times 
\biggl|(s\widehat{Q}_{\omega_{j}}(x,y+z)+(1-s)\widehat{Q}_{\omega_{j}}(x,z)) - R\biggl|
dtds \\
& \lesssim (p-1)
(\|\widehat{Q}_{\omega_{j}}(x, \cdot)
\|_{L_{y}^{\infty}(\mathbb{T})} 
+ |R(x)|^{p-2}) 
\|\widehat{Q}_{\omega_{j}}(x, \cdot) 
- R(x)\|_{L_{y}^{\infty}(\mathbb{T})}. 
\end{split}
\end{equation*}
This together with \eqref{eq-diff1} 
yields that 
\begin{equation} \label{eq-diff3}
\begin{split}
&\left\| |D_{y}|^{\alpha} \left( |\widehat{Q}_{\omega_{j}}|^{p-1}\widehat{Q}_{\omega_{j}}(x)-pR^{p-1}(x)\widehat{Q}_{\omega_{j}}(x)\right)\right\|_{L_{y}^{2}( \T)}^{2} \\
& \lesssim 
(\|\widehat{Q}_{\omega_{j}}(x, \cdot)
\|_{L_{y}^{\infty}(\mathbb{T})} 
+ |R(x)|^{p-2}) 
\|\widehat{Q}_{\omega_{j}}(x, \cdot) 
- R(x)\|_{L_{y}^{\infty}(\mathbb{T})}
\int_{\mathbb{T}} 
\int_{- \pi}^{\pi} 
\frac{\left| 
\left( \widehat{Q}_{\omega_{j}}(x,y+z)-\widehat{Q}_{\omega_{j}}(x,z)\right)\right|^{2}}{|y|^{1+2 \alpha}} d y d z
\\
& \leq 
(\|\widehat{Q}_{\omega_{j}}(x, \cdot)
\|_{L_{y}^{\infty}(\mathbb{T})} 
+ |R(x)|^{p-2}) 
\|\widehat{Q}_{\omega_{j}}(x, \cdot) 
- R(x)\|_{L_{y}^{\infty}(\mathbb{T})}
\left\| |D_{y}|^{\alpha}\widehat{Q}_{\omega_{j}}(x, \cdot) \right\|_{L_{y}^{2}(\T)}^{2}. 
\end{split}
\end{equation}
Integrating \eqref{eq-diff2} and 
\eqref{eq-diff3} with respect to $x$ in $\R$, 
we conclude the proof of this lemma.
\end{proof}

Using Proposition \ref{thm-3} and Lemma \ref{A1}, 
we can show the following: 
\begin{lemma} \label{thm-4}
Let $1 < p < 5$. 
There exists a sufficiently small 
$\omega_{S} \in (0, \infty)$ 
such that 
\begin{equation*}
|D_{y}|^{\frac{1}{2}} \widetilde{Q}_{\omega} = 0 
\qquad \mbox{for all $0 < \omega < \omega_{S}$}. 
\end{equation*}
\end{lemma}
\begin{proof}
We take $\{\omega_{j}\} \subset \mathbb{R}_{+}$ 
with $\lim_{j \to \infty} \omega_{j} = 0$. 
We put 
$U_{j} = |D_{y}|^{1/2} \widehat{Q}_{\omega_{j}}$. 
From Remark \ref{reg-Q} below, we see that 
$U_{j} \in X$. 
Then, $U_{j}$ satisfies 
\begin{equation*}
- \partial_{xx}^{2} U_{j} + \omega_{j}^{-1} |D_{y}|U_{j} 
+ U_{j} - |D_{y}|^{1/2} 
(|\widehat{Q}_{\omega_{j}}|^{p-1} 
\widehat{Q}_{\omega_{j}}) = 0. 
\end{equation*}
Multiplying the above equation by $\overline{U}_{j}$ and integrating 
the resulting equation, we have 
\begin{equation} \label{eq-7}
\begin{split}
0 
& = \int_{\mathbb{R} \times \mathbb{T}} 
|\partial_{x} U_{j}(x, y)|^{2} dxdy 
+ \omega_{j}^{-1}\int_{\mathbb{R} \times \mathbb{T}} 
|D_{y}| U_{j}(x, y)\overline{U_{j}(x, y)} dxdy 
\\
& \quad+ \int_{\mathbb{R} \times \mathbb{T}}
|U_{j}(x, y)|^{2} dxdy - \int_{\mathbb{R} \times \mathbb{T}} 
|D_{y}|^{\frac{1}{2}} (|\widehat{Q}_{\omega_{j}}|^{p-1} \widehat{Q}_{\omega_{j}})
\overline{U_{j}(x, y)} dxdy \\
& = \int_{\mathbb{R} \times \mathbb{T}} 
|\partial_{x} U_{j}(x, y)|^{2} dxdy 
+ \omega_{j}^{-1}\int_{\mathbb{R} \times \mathbb{T}} 
|D_{y}| U_{j}(x, y)\overline{U_{j}(x, y)} dxdy \\
& \quad + \int_{\mathbb{R} \times \mathbb{T}}
|U_{j}(x, y)|^{2} dxdy - \int_{\mathbb{R} \times \mathbb{T}} 
|D_{y}|^{\frac{1}{2}} (p |R(x)|^{p-1} 
\widehat{Q}_{\omega_{j}}(x, y)) 
\overline{U_{j}(x, y)} dxdy \\
& \quad + \int_{\mathbb{R} \times \mathbb{T}}
|D_{y}|^{\frac{1}{2}} ( p|R(x)|^{p-1} 
\widehat{Q}_{\omega_{j}}(x, y) 
- |\widehat{Q}_{\omega_{j}}(x, y)|^{p-1} 
\widehat{Q}_{\omega_{j}}(x, y))
\overline{U_{j}(x, y)} dxdy. 
\end{split}
\end{equation}
Putting 
\begin{equation*}
U_{j}(x, y) = \sum_{k \in \mathbb{Z}} e^{i k y}\widehat{U}_{j}(x, k), \qquad 
\widehat{U}_{j}(x, k) = \frac{1}{2\pi} 
\int_{-\pi}^{\pi} U_{j}(x, y) e^{- i k y} dy, 
\end{equation*}
we have $\widehat{U}_{j}(x, 0) = 0$ 
for all $j \in \mathbb{N}$ because $U_{j} = |D_{y}|^{1/2} \widehat{Q}_{\omega_{j}}$. 
Then, it follows 
from the Parseval identity that 
\begin{equation} \label{eq-8}
\begin{split}
& \int_{\mathbb{R} \times \mathbb{T}} 
|\partial_{x} U_{j}(x, y)|^{2} dxdy 
+ \omega_{j}^{-1}\int_{\mathbb{R} \times \mathbb{T}} 
|D_{y}| U_{j}(x, y)\overline{U_{j}(x, y)} dxdy \\
& + \int_{\mathbb{R} \times \mathbb{T}}
|U_{j}(x, y)|^{2} dxdy - \int_{\mathbb{R} \times \mathbb{T}} |D_{y}|^{\frac{1}{2}} (p |R(x)|^{p-1} 
\widehat{Q}_{\omega_{j}}(x, y)) 
\overline{U_{j}(x, y)} dxdy \\
& = \int_{\mathbb{R} \times \mathbb{T}} 
|\partial_{x} U_{j}(x, y)|^{2} dxdy 
+ \omega_{j}^{-1}\int_{\mathbb{R} \times \mathbb{T}} 
|D_{y}| U_{j}(x, y)\overline{U_{j}(x, y)} dxdy \\
& + \int_{\mathbb{R} \times \mathbb{T}}
|U_{j}(x, y)|^{2} dxdy - p \int_{\mathbb{R} \times \mathbb{T}} |R(x)|^{p-1} 
|U_{j}(x, y)|^{2} dxdy \\
& \geq 
\int_{\mathbb{R} \times \mathbb{T}} 
|\partial_{x} U_{j}(x, y)|^{2} dxdy 
+ \omega_{j}^{-1}
\sum_{k \in \mathbb{Z} \setminus \{0\}}
\int_{\mathbb{R}} |k| |\widehat{U}_{j}(x, k)|^{2} dx \\
& + \int_{\mathbb{R} \times \mathbb{T}}
|U_{j}(x, y)|^{2} dxdy - C \sum_{k \in \mathbb{Z} \setminus \{0\}}
\int_{\mathbb{R}} |\widehat{U}_{j}(x, k)|^{2} dx \\
& \geq \int_{\mathbb{R} \times \mathbb{T}} 
|\partial_{x} U_{j}(x, y)|^{2} dxdy + \int_{\mathbb{R} \times \mathbb{T}}
|D_{y}| U_{j}(x, y)\overline{U_{j}(x, y)} dxdy 
+ \int_{\mathbb{R} \times \mathbb{T}} 
|U_{j}(x, y)|^{2} dxdy \\
& = \|U_{j}\|_{X}^{2} 
\end{split}
\end{equation}
for sufficiently large $j \in \mathbb{N}$. 
%It follows that 
%\begin{equation}
%I + II \geq \frac{1}{2} \|U_{j}\|_{X}^{2} \qquad 
%\mbox{for sufficiently large $j \in N$}. 
%\end{equation} 
\par
Next, we shall show that 
\begin{equation} \label{eq3-vari-chara}
\biggl|\int_{\mathbb{R} \times \mathbb{T}}
|D_{y}|^{\frac{1}{2}} ( p|R(x)|^{p-1} \widehat{Q}_{\omega_{j}}(x, y) 
- |\widehat{Q}_{\omega_{j}}(x, y)|^{p-1} \widehat{Q}_{\omega_{j}}(x, y))
\overline{U_{j}(x, y)} dxdy \biggl| \leq \frac{1}{2} \|U_{j}\|_{X}^{2} 
\end{equation}
for sufficiently large $j \in \mathbb{N}$.
We first consider the case of $1 < p \leq 2$. 
From Proposition \ref{thm-3}, we see that 
$\lim_{j \to \infty} \widehat{Q}_{\omega_{j}} 
= R$ strongly in the energy space $X$.
%$L^{q}_{x, y}
%(\mathbb{R} \times \mathbb{T})$ for 
%$2 \leq q \leq 6$. 
This together with Remark \ref{reg-Q} below 
yields that 
$\lim_{j \to \infty} \widehat{Q}_{\omega_{j}} 
= R$ strongly in $L^{\infty}_{x, y}(\mathbb{R} \times \mathbb{T})$. 
%for $2 \leq q < \infty$. 
%Then, it follows that $\lim_{j \to \infty} \widehat{Q}_{\omega_{j}}^{p-1} = R^{p-1}$ strongly in 
%$L^{\frac{q}{p-1}}_{x, y}
%(\mathbb{R} \times \mathbb{T})$ for 
%$2 \leq q < \infty$
%(see e.g. Willem~\cite[Appendix A]{Willem}). 
Thus, applying Lemma \ref{A1} with 
$\alpha = 1/2$, by the H\"{o}lder inequality, we have
\begin{equation} \label{eq-8-00}
\begin{split}
& \int_{\mathbb{R} \times \mathbb{T}}
|D_{y}|^{\frac{1}{2}}\left( p R(x)^{p-1}\widehat{Q}_{\omega_{j}}(x, y) 
- |\widehat{Q}_{\omega_{j}}(x, y)|^{p-1} 
\widehat{Q}_{\omega_{j}}(x, y)\right) \overline{U_{j}(x, y)} dxdy \\
& \leq 
\left(\int_{\mathbb{R} \times \mathbb{T}} 
\left||D_{y}|^{\frac{1}{2}} 
(p R^{p-1}(x) \widehat{Q}_{\omega_{j}}(x, y) 
- |\widehat{Q}_{\omega_{j}}(x, y)|^{p-1}(x, y) \widehat{Q}_{\omega_{j}}(x, y))\right|^{2} 
dxdy\right)^{\frac{1}{2}} \times \\ 
& \quad \times 
\|U_{j}\|_{L_{x, y}^{2}(\R \times \T)} \\
& \lesssim \|R - \widehat{Q}_{\omega_{j}}
\|^{p-1}_{L_{x, y}^{\infty}(\R \times \T)} \||D_{y}|^{\frac{1}{2}} \widehat{Q}_{\omega_{j}}
\|_{L_{x, y}^{2}(\R 
\times \T)} \|U_{j}\|_{L_{x, y}^{2}(\R \times \T)} \\
& < \frac{1}{2} \|U_{j}\|_{X} 
\end{split}
\end{equation}
for sufficiently large $j \in \mathbb{N}$. 
Thus, \eqref{eq3-vari-chara} holds 
for all $1 < p \leq 2$. 
For the case of $2 < p < 5$, 
we can prove \eqref{eq3-vari-chara} 
by the similar argument above. 

\par
It follows from 
\eqref{eq-7}, \eqref{eq-8} and \eqref{eq3-vari-chara} 
that 
there exists $j_{0} \in \mathbb{N}$ such that 
\begin{equation*}
0 \geq \frac{1}{2} \|U_{j}\|_{X}^{2} \qquad 
\mbox{for all $j \geq j_{0}$}. 
\end{equation*}
Thus, we find that $U_{j} = 0$ for all $j \geq j_{0}$. 
This completes the proof. 

\end{proof}
We are now in a position to prove 
Proposition \ref{vari-omega-small}. 
\begin{proof}[Proof of Proposition \ref{vari-omega-small}]
It follows from Lemma \ref{thm-4} that 
for each $y \in \mathbb{T}$, 
$\widetilde{Q}_{\omega}$ satisfies 
\begin{equation*}
- \partial_{xx} \widetilde{Q}_{\omega} + \widetilde{Q}_{\omega} 
- |\widetilde{Q}_{\omega}|^{p-1} \widetilde{Q}_{\omega} = 0 \qquad 
\mbox{in $\mathbb{R}$} 
\end{equation*}
for $0 < \omega < \omega_{S}$.
It is known that the solution of the above equation which 
decays at infinity is unique. 
Therefore, we have 
$\widetilde{Q}_{\omega}(x, y) = R(x)$. 
By scaling, we obtain $Q_{\omega}(x, y) = R_{\omega}(x)$. 
This completes the proof. 
\end{proof}

\subsection{Proof of Theorem \ref{thm-chara-ground}}
In this subsection, we give a proof of 
Theorem \ref{thm-chara-ground}. 
To this end, we prepare several lemmas, 
which are needed later. 
\begin{lemma}\label{variational-chara-thm4}
Suppose that there exists $\omega_{0} > 0$ such that 
$m_{\omega_{0}} < 2\pi m_{\omega_{0}, \mathbb{R}}$. 
Then, for any $\omega > \omega_{0}$, 
we have $m_{\omega} < 2\pi m_{\omega, \mathbb{R}}$. 
\end{lemma}
\begin{proof}
Since $m_{\omega_{0}} < 2\pi m_{\omega_{0}, \mathbb{R}}$, 
there exists a ground state $Q_{\omega_{0}} \neq R_{\omega_{0}}$ 
to \eqref{sp}. 
We put 
\begin{equation*}
V_{\omega}(x, y) = \lambda_{\omega}^{\frac{1}{p-1}}Q_{\omega_{0}}
(\sqrt{\lambda_{\omega}}x, y), \qquad
\lambda_{\omega} = \frac{\omega}{\omega_{0}} 
\ (> 1). 
\end{equation*}
Then, $V_{\omega}$ satisfies 
\begin{equation*}
- \partial_{xx} V_{\omega} + \lambda_{\omega} |D_{y}| V_{\omega} 
+ \omega V_{\omega} - V_{\omega}^{p} = 0 \qquad 
\mbox{in $\mathbb{R} \times \mathbb{T}$}. 
\end{equation*}
It follows that 
\begin{equation} \label{N-neg}
\mathcal{N}_{\omega}(V_{\omega}) 
= - (\lambda_{\omega} -1) 
\||D_{y}|^{\frac{1}{2}} V_{\omega}\|_{L^{2}(\R \times \T)}^{2} < 0. 
\end{equation}
Recall that
\[
m_{\omega, \mathbb{R}} = \omega^{\frac{p+1}{p-1} - \frac{1}{2}} m_{\mathbb{R}}.
\]
This together with Lemma \ref{existence-ground-thm2} 
and \eqref{N-neg} gives
\begin{equation*}
\begin{split}
%& 
m_{\omega} 
\leq \mathcal{I}_{\omega}(V_{\omega})
\leq \lambda_{\omega}^{\frac{p+1}{p-1} -\frac{1}{2}} 
\mathcal{I}_{\omega_{0}}(Q_{\omega_{0}})
= \lambda_{\omega}^{\frac{p+1}{p-1} 
- \frac{1}{2}} m_{\omega_{0}} % \\& 
< 2\pi 
\lambda_{\omega}^{\frac{p+1}{p-1} - \frac{1}{2}} m_{\omega_{0}, \mathbb{R}} 
= 2\pi m_{\omega, \mathbb{R}}. 
\end{split}
\end{equation*}
Thus, we obtain $m_{\omega} < 2\pi 
m_{\omega, \mathbb{R}}$. 
\end{proof}
\begin{lemma}\label{variational-chara-thm5}
Let $Q_{\omega}$ be the ground state to \eqref{sp}. 
Then, $Q_{\omega}$ becomes a minimizer for the 
following minimization problem: 
\begin{equation*}
M_{\omega} = 
\inf_{u \in X \setminus \{0\}}
\dfrac{\|\partial_{x} u\|_{L_{x, y}^{2}(\R \times \T)}^{2} + \||D_{y}|^{\frac{1}{2}} u\|_{L_{x, y}^{2}(\R \times \T)}^{2} 
+ \omega \|u\|_{L_{x, y}^{2}(\R \times \T)}^{2}}
{\left(\|u\|_{L_{x, y}^{p+1}(\R \times \T)}^{p+1} 
\right)^{\frac{2}{p+1}}}. 
\end{equation*}
\end{lemma}
\begin{proof}
We put 
\begin{equation*}
m_{1, \omega}
:= 
\frac{p-1}{2(p+1)} 
\inf_{u \in X \setminus \{0\}}
\dfrac{\left(\|\partial_{x} 
u\|_{L_{x, y}^{2}(\R \times \T)}^{2} 
+ \||D_{y}|^{\frac{1}{2}} u\|_{L_{x, y}^{2}(\R \times \T)}^{2} 
+ \omega \|u\|_{L_{x, y}^{2}
(\R \times \T)}^{2} \right)^{\frac{p+1}{p-1}}}
{\left(\|u\|_{L_{x, y}^{p+1}(\R \times \T)}^{p+1} 
\right)^{\frac{2}{p-1}}}.
\end{equation*}
Since $Q_{\omega}$ is a minimizer for $m_{\omega}$, 
in order to prove Lemma \ref{variational-chara-thm5}, 
it is enough to show that $m_{\omega} = m_{1, \omega}$. 

\par 
For each $u \in X \setminus \{0\}$, 
we put 
\begin{equation*}
t(u) = \left(\dfrac{\|\partial_{x} 
u\|_{L_{x, y}^{2}(\R \times \T)}^{2} 
+ \||D_{y}|^{\frac{1}{2}} u\|_{L_{x, y}^{2}(\R \times \T)}^{2} 
+ \omega \|u\|_{L_{x, y}^{2}(\R \times \T)}^{2}}
{\|u\|_{L_{x, y}^{p+1}(\R \times \T)}^{p+1}} 
\right)^{\frac{2}{p-1}}. 
\end{equation*}
Then, we see that $\mathcal{N}_{\omega}(t(u) u) = 0$. 
It follows from the definition of $m_{\omega}$ that 
\begin{equation*}
m_{\omega} \leq \mathcal{S}_{\omega}(t(u) u) 
= \frac{p-1}{2(p+1)} 
\dfrac{\left(\|\partial_{x} 
u\|_{L_{x, y}^{2}(\R \times \T)}^{2} + 
\||D_{y}|^{\frac{1}{2}} u\|_{L_{x, y}^{2}(\R \times \T)}^{2} 
+ \omega \|u\|_{L_{x, y}^{2}(\R \times \T)}^{2} \right)^{\frac{p+1}{p-1}}}
{\left(\|u\|_{L_{x, y}^{p+1}(\R \times \T)}^{p+1} 
\right)^{\frac{2}{p-1}}}. 
\end{equation*}
Taking an infimum on $u \in X$, we have 
\begin{equation} \label{eq1-vari-thm5}
m_{\omega} \leq m_{1, \omega}.
\end{equation} 

Let $Q_{\omega}$ be the ground state to \eqref{sp}. 
Since $\mathcal{N}_{\omega}(Q_{\omega}) = 0$, 
we have 
\begin{equation*}
\begin{split}
m_{1, \omega} 
& \leq \frac{p-1}{2(p+1)} 
\dfrac{\left(\|\partial_{x} Q_{\omega} 
\|_{L_{x, y}^{2}(\R \times \T)}^{2} 
+ \||D_{y}|^{\frac{1}{2}} Q_{\omega}\|_{L_{x, y}^{2}(\R \times \T)}^{2} 
+ \omega \|Q_{\omega}\|_{L_{x, y}^{2}(\R \times \T)}^{2} \right)^{\frac{p+1}{p-1}}}
{\left(\|Q_{\omega}\|_{L_{x, y}^{p+1}
(\R \times \T)}^{p+1} \right)^{\frac{2}{p-1}}} \\
& = \frac{p-1}{2(p+1)} 
\left(\|\partial_{x} Q_{\omega} \|_{L_{x, y}^{2}(\R \times \T)}^{2} 
+ \||D_{y}|^{\frac{1}{2}} Q_{\omega}\|_{L_{x, y}^{2}(\R \times \T)}^{2} 
+ \omega \|Q_{\omega}\|_{L_{x, y}^{2}(\R \times \T)}^{2} \right) \\
& = \mathcal{I}_{\omega} (Q_{\omega}) = m_{\omega}. 
\end{split}
\end{equation*}
This together with \eqref{eq1-vari-thm5} 
yields that $m_{\omega} = m_{1, \omega}$. 
\end{proof}
\begin{lemma}\label{variational-chara-thm6}
Let $Q_{\omega}$ be the ground state to \eqref{sp} 
and $L_{\omega, \text{g}, +}$ be the operator defined by 
\begin{equation} \label{lop-g}
L_{\omega, \text{g}, +} = - \partial_{xx} + |D_{y}| 
+ \omega - p |Q_{\omega}|^{p-1}
\end{equation}
with the domain $D(L_{\omega, \text{g}, +}) 
= X_{2}$, where 
$X_2:=H^2_xL^2_y \cap L^2_xH^{1}_y(\R \times \T)$. 
Let $\Lambda_{2}(\omega)$ be the second eigenvalue of 
the operator $L_{\omega, \text{g}, +}$. 
Then, we have $\Lambda_{2}(\omega) \geq 0$. 
\end{lemma}
\begin{proof}
For each $\varphi \in X$, we put 
\begin{equation*}
f_{\varphi}(t) := \dfrac{\|\partial_{x} 
(u + t \varphi) \|_{L_{x, y}^{2}(\R \times \T)}^{2} 
+ \||D_{y}|^{\frac{1}{2}} 
(u + t \varphi)\|_{L_{x, y}^{2}(\R \times \T)}^{2} 
+ \omega \|(u + t \varphi)\|_{L_{x, y}^{2}(\R \times \T)}^{2}}
{\left(\|(u + t \varphi)\|_{L_{x, y}^{p+1}
(\R \times \T)}^{p+1} \right)^{\frac{2}{p+1}}}. 
\end{equation*}
It follows from Lemma \ref{variational-chara-thm6} that 
$Q_{\omega}$ is the minimizer for $M_{\omega}$. 
Therefore, we have 
$f_{\varphi}^{\prime}(0) = 0$ and $f_{\varphi}^{\prime \prime}(0) \geq 0$. 
This yields that 
\begin{equation*}
0 \leq f_{\varphi}^{\prime \prime}(0) 
= \langle L_{\omega, \text{g}, +} \varphi, \varphi \rangle 
- \frac{p+1}{p-1} \frac{\langle Q_{\omega}^{p}, \varphi \rangle}
{\|Q_{\omega}\|_{L_{x, y}^{p+1}(\R \times \T)}^{p+1}}. 
\end{equation*}
Thus, it follows that 
\begin{equation*}
0 \leq \inf_{\varphi \in X, \varphi \perp Q_{\omega}^{p}} 
\langle L_{\omega, \text{g}, +} \varphi, \varphi \rangle. 
\end{equation*}
Then, by the mini-max theorem, we have 
\begin{equation*}
\Lambda_{2}(\omega) = 
\max_{\dim V = 1, V \subset X} \inf_{\varphi \perp V} 
\langle L_{\omega, \text{g}, +} \varphi, \varphi \rangle
\geq \inf_{\varphi \in X, \varphi \perp Q_{\omega}^{p}} 
\langle L_{\omega, \text{g}, +} \varphi, \varphi \rangle \geq 0. 
\end{equation*}
This completes the proof. 
\end{proof}

We are now in a position to prove Theorem \ref{thm-chara-ground}. 
\begin{proof}[Proof of Theorem \ref{thm-chara-ground}]
Let 
\begin{equation*}
\begin{split}
\omega_{*} 
& = \sup\left\{ 
\omega > 0 \colon m_{\omega} = 2\pi m_{\omega, \mathbb{R}}
\right\} \\
&= \sup\left\{ 
\omega > 0 \colon 
\mathcal{S}_{\omega}(Q_{\omega}) = 
\mathcal{S}_{\omega} (R_{\omega}) 
\right\}. 
\end{split}
\end{equation*}
Then, by Propositions \ref{vari-omega-large} and 
\ref{vari-omega-small}, we have $0 < \omega_{*} < \infty$. 
It follows from the definition of $\omega_{*}$ 
and Lemma \ref{variational-chara-thm4} that 
\[
\mathcal{S}_{\omega} (Q_{\omega}) = 
m_{\omega} < 2\pi m_{\omega, \R} = \mathcal{S}_{\omega} 
(R_{\omega})
\]
for $\omega > \omega_{*}$, 
which implies that $Q_{\omega} \neq R_{\omega}$. 
\par 
Next, we shall show that for $0 < \omega \leq \omega_{*}$, 
$Q_{\omega} = R_{\omega}$. 
From the definition of $\omega_{*}$ 
and Lemma \ref{variational-chara-thm4}, 
we see that $m_{\omega} = 2\pi m_{\omega, \R}$ for 
$0 < \omega < \omega_{*}$. 
Furthermore, since $m_{\omega}$ and $m_{\omega, \R}$ 
are continuous (see Lemma \ref{thmc-1} below), we have 
$m_{\omega_{*}} = 2\pi m_{\omega_{*}, \R}$.
Multiplying \eqref{sp} by $\overline{Q_{\omega}(x, y)}$ 
and integrating the resulting equation over $\R$, 
we have 
\[
\mathcal{N}_{\omega, \R}(Q_{\omega}) = 
- \int_{\R} ||D_{y}|^{\frac{1}{2}}Q_{\omega}|^{2} dx \leq 0. 
\]
This together with \eqref{R-vari-chara} 
yields that 
\[
m_{\omega, R} \leq \mathcal{I}_{\omega, \R}(Q_{\omega}). 
\]
Integrating the above over $\T$, one has 
\[
2\pi m_{\omega, \R} \leq 
\int_{\T} \mathcal{I}_{\omega, \R}(Q_{\omega}) dy 
= \mathcal{I}_{\omega}(Q_{\omega}) 
- \||D_{y}|^{\frac{1}{2}} Q_{\omega}\|_{L_{x, y}^{2}
(\R \times \T)}^{2} \leq 
\mathcal{I}_{\omega}(Q_{\omega}) = m_{\omega} 
= 2\pi m_{\omega, \R}. 
\]
This implies that 
$\||D_{y}|^{\frac{1}{2}} Q_{\omega}\|_{L_{x, y}^{2}
(\R \times \T)}^{2} = 0$. 
Thus, we see that $|D_{y}|^{1/2} Q_{\omega} = 0$. 
Since $Q_{\omega}$ is smooth (see 
Remark \ref{reg-Q}), 
$Q_{\omega}(x, y)$ does not depend on 
$y$ and satisfies \eqref{line-sp}. 
Thus, from the uniqueness of the solution to \eqref{line-sp}, 
we have $Q_{\omega} = R_{\omega}$. 
\par
Finally, we shall show that $\omega_{*} \leq \omega_{p}$ 
by contradiction. 
Suppose the contrary that $\omega_{*} > \omega_{p}$. 
Then, for $\omega \in (\omega_{p}, \omega_{*})$, 
one has $Q_{\omega} = R_{\omega}$, 
which means that $L_{\omega, g, +} = L_{\omega, +}$. 
Let $\lambda_{2}(\omega)$ be the second eigenvalue of $L_{\omega, +}$.
It follows from Lemma \ref{lem:neg-eigenvalue} that 
for $\omega > \omega_{p}$, $\lambda_{2}(\omega) < 0$. 
On the other hand, 
from Lemma \ref{variational-chara-thm6}, 
we see that $\lambda_{2}(\omega) 
= \Lambda_{2}(\omega) \geq 0$, 
which is a contradiction. 
Thus, we find that $\omega_{*} \leq \omega_{p}$. 
This completes the proof. 
\end{proof}

\section{Regularity of the ground state}
In this section, we would like to show that $Q_{\omega} 
\in L_{x, y}^{\infty}(\R \times \T)$ with a uniform bound.
\begin{theorem}\label{elliptic-regularity}
Assume that $\omega > 0$. 
Let $Q_{\omega} \in X$ be a solution to \eqref{sp}. 
For any $M>0$, there exists $C = C(M) > 0$ such that $\|Q_{\omega}\|_{X} \leq M$ implies 
$\|Q_{\omega}\|_{L_{x, y}^{\infty}(\R \times \T)} \leq C$. 
\end{theorem}
\begin{proof} 
For simplicity, we consider the case where $\omega = 1$ only. 
Putting $Q:= Q_{1}$, we see that 
$Q \in X (= H^1_xL^2_y \cap L^2_xH^{\frac{1}{2}}_y
(\R \times \T))$ satisfies 
\begin{equation}
\label{eq1}
- \partial_{xx} Q+ |D_{y}| Q + Q
- |Q|^{p-1}Q = 0 \quad \mbox{in $\R \times \T$}. 
\end{equation}
We first consider the case of $p \in (1, 3]$. 
It follows from the Sobolev 
embedding $X \hookrightarrow 
L_{x, y}^{q} (\R \times \T)$ for $2 \leq q \leq 6$ 
that $Q \in L_{x, y}^{2p}(\R \times \T)$. 
Using equation \eqref{eq1}, we have
\begin{equation} \label{eq1-regu}
\left\|(- \partial_{xx} + |D_{y}| + 1) Q\right\|_{L_{x, y}^2
(\R \times \T)} = 
\left\| |Q|^{p-1}Q \right\|_{L_{x, y}^{2}(\R \times \T)} 
= \left\| Q \right\|_{L_{x, y}^{2p}(\R \times \T)}^p,
\end{equation}
which yields that 
$Q \in X_2$, where 
$X_{2} = H_{x}^{2}L_{y}^{2} \cap 
L_{x}^{2} H_{y}^{1}(\R \times \T)$. 
We recall that from \cite[Theorem 10.4]{BON}, we have
\begin{equation}
\label{embedding-i}
W^{2, q}_xL^{q}_y\bigcap L^{q}_xW^{1,q}_y
(\R \times \T) 
\hookrightarrow
L_{x, y}^{\infty}(\R \times \T) 
\end{equation}
with $q \geq 1$ satisfying 
\[
\frac{1}{q} \leq \frac{2}{3}. 
\] 
This yields that
$X_2 \hookrightarrow L^\infty_{x,y}(\R \times \T)$ so that $Q \in L^\infty_{x,y}(\R \times \T)$.
%Now, since $\|\varphi\|_{X_{2}}=\|\varphi\|_{L^{2 p}}^{ p}$ is a direct 
%consequence from \eqref{eq1},
\par
Thus, we may restrict ourselves to $p \in (3, 5)$. 
Let $\boF_{\R \times \T}$ and 
$\boF^{-1}_{\R \times \T}$ denote the Fourier and the Fourier inverse transforms on $\R \times \T$, respectively. 
Namely, 
\[
\boF_{\R \times \T}(f)(\xi, k) = \frac{1}{(2\pi)^{\frac{3}{2}}} 
\int_{\R \times \T} f(x, y) e^{- i (\xi x + k y)} dxdy 
\]
for $f \in L^{2}(\R \times \T)$ and 
\[
\boF^{-1}_{\R \times \T}(g)(x, y) = 
\sum_{k \in \mathbb{Z}} 
e^{i k x} \int_{\R} e^{i \xi x} g(\xi, k) d\xi
\]
for $g \in L^{2}(\mathbb{R}, \ell^{2})$. 
We know that \eqref{eq1} can be written as 
%It remains to study the following problem
\begin{equation}
\label{convolution}
Q = \boF_{\R \times \T}^{-1}\left( m 
\boF_{\R \times \T}\left(|Q|^{p-1}Q\right) \right),
\end{equation}
%$$$$
where the Fourier multiplier $m$ is defined by
\[
m(\xi, \eta):= \frac{1}{\xi^2+|\eta|+1 } \quad {\rm for} \ \ 
(\xi, \eta)\in \R^{2}. 
\]

%Using Mikhlin-Hormander theorem (see [Schlag book Chapter 8 and Cardona Theorem 1.2]),
We put
\begin{equation} \label{fourier-m}
m_1(\xi, \eta) := 
\frac{\xi^2}{\xi^2+|\eta|+1 }, \qquad 
m_2(\xi, \eta):= \frac{|\eta|}{\xi^2+|\eta|+1 } 
\qquad \mbox{for $(\xi, \eta) \in \R^{2}$} 
\end{equation}
and define operators by 
$\widetilde{T}_{1} := \boF_{\R \times \T}^{-1}
[m_{1}(\xi, k) \boF_{\R \times \T}]$ 
and 
$\widetilde{T}_{2} := \boF_{\R \times \T}^{-1}
[m_{2}(\xi, k) \boF_{\R \times \T}]$.
Note that from \eqref{convolution}, we have 
\begin{equation} \label{conv2}
- \partial_{xx} Q = \widetilde{T}_{1} (|Q|^{p-1}Q), \qquad 
|D_{y}| Q = \widetilde{T}_{2} (|Q|^{p-1}Q). 
\end{equation} 
It is enough to show that 
$T_{i} \in \mathcal{L}(L_{x, y}^q(\R\times\T))$ 
for any $q \in (1,\infty)$ 
and $i = 1, 2$, i.e
\begin{equation}
\label{multiplier}
\|\widetilde{T}_{i} f\|_{L_{x, y}^{q}(\R \times \T)} = 
\|\boF_{\R \times \T}^{-1}\left( m_i 
\boF_{\R \times \T}\left(f\right) 
\right)\|_{L_{x, y}^q(\R \times \T)} 
\lesssim \| f\|_{L_{x, y}^q(\R \times \T)} \quad {\rm for} \ i=1,2.
\end{equation}
We admit \eqref{multiplier} for a moment and 
continue to prove Theorem \ref{elliptic-regularity}. 
Taking $q=\frac{p+1}{p}$ in \eqref{multiplier},
we obtain through \eqref{conv2}
\begin{equation*}
\begin{split}
&\|\partial_{xx} Q\|_{L_{x, y}^{\frac{p+1}{p}}(\R \times \T)} 
+ \||D_y| Q\|_{L_{x, y}^{\frac{p+1}{p}}(\R \times \T)} \\
& = \|T_{1}(|Q|^{p-1}Q) \|_{L_{x, y}^{\frac{p+1}{p}}(\R \times \T)} + 
\|T_{2} (|Q|^{p-1}Q)\|_{L_{x, y}^{\frac{p+1}{p}}(\R \times \T)}\\
&\lesssim \| Q\|_{L_{x, y}^{p+1}(\R \times \T)}^p 
\lesssim \| Q\|_{X}^p, 
\end{split}
\end{equation*}
which implies that 
\begin{equation} \label{1-re}
Q \in 
W_{x}^{2, \frac{p+1}{p}} L_{y}^{\frac{p+1}{p}} 
\cap 
L_{x}^{\frac{p+1}{p}} W_{y}^{1, \frac{p+1}{p}}
(\R \times \T).
\end{equation} 
Here, we recall that from \cite[Theorem 10.2]{BON}, 
we have
\begin{equation}
\label{embedding}
W^{2, q}_xL^{q}_y\bigcap L^{q}_xW^{1,q}_y
(\R \times \T) 
\hookrightarrow
L_{x, y}^{r}(\R \times \T) 
\end{equation}
with $r, q \geq 1$ satisfying 
\[
\frac{1}{q} - \frac{1}{r} \leq \frac{2}{3} 
\]
From \eqref{embedding}, 
we can take 
$q =\frac{p+1}{p}$, 
so that $Q \in L_{x, y}^{r}(\R \times \T)$ 
for all $2 \leq r \leq r_{1}$,
where 
\[
r_{1} = \frac{3(p+1)}{p-2}. 
\]
%Putting 
%\[
%p_{1} = \frac{7 + \sqrt{73}}{4}, 
%\]
%we see that 
%$r_{1} \geq 2p$ for $3 < p \leq p_{1}$. 
%Thus, we have $Q \in L^{2p}(\R \times \T)$. 
%This together with \eqref{eq1-regu} yields that 
%$Q \in L^{\infty}(\R \times \T)$. 
Note that $|Q|^{p-1}Q \in L_{x, y}^{\frac{r_{1}}{p}} 
(\R \times \T)$. 
Thus, repeating the above procedure, we see that 
\begin{equation} \label{2-re}
Q\in 
W_{x}^{2, \frac{r_{1}}{p}} L_{y}^{\frac{r_{1}}{p}} 
\cap 
L_{x}^{\frac{r_{1}}{p}} W_{y}^{1, \frac{r_{1}}{p}}
(\R \times \T).
\end{equation} 
Putting $p_{1} := 2 + \sqrt{6}$. 
we can easily verify that 
$p = p_{1}$ is a solution to the following 
quadratic equation:
\[
p^{2} - 4p -2 = 0. 
\]
We see from \eqref{embedding-i} and \eqref{2-re} 
that $Q \in L^{\infty}(\R \times \T)$ 
for $3 \leq p \leq p_{1}$. 
Thus, it suffices to consider the case of $p \in 
(p_{1}, 5)$. 
Putting $q =\frac{r_{1}}{p}$ in \eqref{embedding}, 
we infer that $Q \in L_{x, y}^{r}(\R \times \T)$ 
for all $2 \leq r \leq r_{2}$,
where 
\[
r_{2} = \frac{3(p+1)}{p^{2}-4p-2}. 
\] 
for $p_{1} < p < 5$. 

\par 
Consider the sequence $\{r_{i}\} \subset \R 
\cup \{\infty\}$ defined by 
\begin{equation*}
r_{1} = \frac{3(p+1)}{p-2}, \qquad 
\frac{1}{r_{i+1}} = 
\begin{cases}
\frac{p}{r_{i}} - \frac{2}{3} & \qquad \mbox{if $r_{i} < \frac{3p}{2}$}, \\
0 & \qquad \mbox{if $r_{i} \geq \frac{3p}{2}$}. 
\end{cases}
\end{equation*}
Namely, 
\begin{equation*}
r_{i+1} = 
\begin{cases}
\left(p^{i+1}\frac{p-5}{3(p^{2} -1)} + \frac{2}{3(p-1)} \right)^{-1} 
& \qquad \mbox{if $r_{i} < \frac{3p}{2}$}, 
\\
\infty & \qquad \mbox{if $r_{i} \geq \frac{3p}{2}$}. 
\end{cases}
\end{equation*}
Repeating the above argument, we see that 
$Q \in L^{r_{i}}(\R \times \T)$ for any $i \in \mathbb{N}$. 
It follows that 
for any $p \in (3, 5)$, 
$r_{i} \geq \frac{3p}{2}$ for $i \geq i_{p}$, where 
\[
i_{p} := \frac{1}{\log p} \log \left[\frac{2(p+1)}{p(5-p)} 
\right].
\] 
This implies $Q \in L^{\infty}(\R \times \T)$ 
for any $p \in (3, \infty)$. 

\par
Thus, admitting that \eqref{multiplier} holds, 
we obtain the desired result. 
\end{proof}
\begin{remark} \label{reg-Q}
As a consequence of this result, 
we see that $Q \in C_{x, y}^{\infty}(\R\times\T)$. 
More precisely, we infer that 
$Q \in \bigcap_{k\geq 1}X_k,$ with $X_k:=H^k_xL^2_y \cap L^2_xH^{\frac{k}{2}}_y(\R \times \T)$. 
Indeed, taking $\partial_{x}^{k-1}, |D_{y}|^{\frac{k-1}{2}}$ 
derivatives of \eqref{eq1}, 
and using the H\"{o}lder inequality, one gets 
\[
\|Q\|_{X_{k+1}} \lesssim\|Q\|_{X_{2}}^{p-1}\|Q\|_{X_{k-1}}
\]
for $k \geq 2$. A simple induction as in 
\cite{DebSau} then implies that
\[
Q \in \bigcap_{k \geq 1} X_{k}.
\]
\end{remark}

Now, we come back to prove \eqref{multiplier}. 
First, we recall the Mikhlin-Hormander theorem in the 
$\R^{2}$ case (see e.g. \cite[Chapter 8]{Schlag} for the proof). 
We denotes by $\mathcal{S}_{x, y}(\R^{2})$ the set of 
Schwartz class functions on $\R^{2}$. 
We denote by $\boF_{\R^{2}}[f]$ and 
$\boF_{\R^{2}}^{-1}[f]$ the Fourier transform and 
the inverse Fourier transform of $f$, that is, 
\[
\boF_{\R^{2}}[f](\xi) = \frac{1}{2\pi} \int_{\R^{2}} 
e^{- i x \xi} f(x) dx, \qquad 
\boF_{\R^{2}}^{-1}[g](x) 
= \int_{\R^{2}} e^{i x \xi} g(\xi) d\xi 
\] 
\begin{theorem}
\label{Hormander-thm}
Let $m: \mathbb{R}^{2} 
\rightarrow \mathbb{C}$ satisfy
\begin{equation}
\label{Cond:m}
\left|\partial^{\gamma} m(\xi)\right| \leq B|\xi|^{-|\gamma|}
\end{equation}
for any multi-index $\gamma$ of length
$|\gamma| \leq 4$ and for all $\xi \in \R^{2} 
\setminus \{0\} .$ Let $T$ be the operator defined on $\mathcal{S}_{x, y}(\R^{2})$ by 
$T(f):=\boF_{\R^{2}}^{-1}\left( m \boF_{\R^{2}} 
\left(f\right) \right)$. 
Then, for any $1<q<\infty,$ there exits 
a constant $C=C(q)$ such that
\[
\left\|T(f)\right\|_{L_{x, y}^q(\R^{2})} \leq C B\|f\|_{L_{x, y}^q(\R^{2})}
\]
for all $f \in \mathcal{S}_{x, y}(\R^{2})$. 
\end{theorem}
\begin{remark}
\label{Remark-1}
We define operators by 
$T_{1} := \boF_{\R^{2}}^{-1}[m_{1} \boF_{\R^{2}}]$ 
and $T_{2} := \boF_{\R^{2}}^{-1}[m_{2} \boF_{\R^{2}}]$, 
where $m_{1}$ and $m_{2}$ is given by 
\eqref{fourier-m}. Theorem \ref{Hormander-thm} asserts that $T_1$ and $T_2$ are bounded operators 
from $ L_{x, y}^{q}\left(\R^{2} \right)$ to itself 
since both $m_1$ and $m_2$ satisfy \eqref{Cond:m}. 
\end{remark}
Next, we will show that 
we can extend 
the result of Theorem \ref{Hormander-thm} 
to the infinite cylinder case $\R \times \T$.
\begin{theorem}
\label{Thm1-period}
Let $1 < q < \infty$ and 
$m \colon \R^{2} \to \mathbb{C}$ 
satisfying \eqref{Cond:m}. 
Putting 
\[
\widetilde{T} f(x,y)=
\boF_{\R \times \T}^{-1} 
\left(m(\xi,k) \boF_{\R \times \T}\right)\left(f\right)
\]
for all $L_{x, y}^{q}(\R \times \T)$, we have 
$\widetilde{T} 
\in \mathcal{L}
(L_{x, y}^{q}\left(\R \times \T \right))$ 
with 
\[
\|\widetilde{T}\|_{\mathcal{L}
(L_{x, y}^{q}\left(\R \times \T \right))} 
\leq \frac{1}{2 \sqrt{\pi}}\|T\|_{\mathcal{L}(L_{x, y}^{q}(\R^{2}))}, 
\]
where $T = \boF_{\R^{2}}^{-1}( m \boF_{\R^{2}} 
)\left(f\right)$.
\end{theorem}

\begin{remark}
Note that \eqref{multiplier} is a direct consequence of Theorem \ref{Thm1-period}, Theorem \ref{Hormander-thm} and 
the statement of Remark \ref{Remark-1}. 
We can easily find that 
$\widetilde{T}_1$ and $\widetilde{T}_2$ are 
the periodized versions of $T_1$ and $T_2$, which are the operators defined in Remark \ref{Remark-1}, respectively. 
Since $T_{1}, T_{2} \in \mathcal{L}
(L^{q}\left(\R^{2}\right))$, we conclude that 
$\widetilde{T}_1, \widetilde{T}_2 
\in \mathcal{L}(L_{x, y}^{q}\left(\R \times \T \right))$.
\end{remark}

The proof is the same as the one of Theorem 3.8 
in \cite[Chapter VII]{Stein-Weiss}. 
We will give it here for the sake of completeness. 
We start with the following definition.
\begin{definition}
Let $N \in \N \cup \{0\}$. 
We say that $f_N$ 
is a \textit{trigonometric polynomial} on $y$ if 
it is given by
\[
f_N(x,y)= \sum_{|k|\leq N} f_k(x) e^{i k y} \quad {\rm for \ any} \ \ (x,y) \in \R \times \T,
\]
with
\[
f_k(x):=\frac{1}{\sqrt{2\pi}} 
\int_{-\pi}^{\pi} f(x, y) e^{- i k y} dy.
\]
\end{definition}

\begin{remark}
Note that the set of trigonometric polynomial on $y$ is dense in $ L^{q}_{x, y} \left(\R \times \T \right)$ and in $C_y^0\left( \T , L_x^1(\R) \right)$. We refer to the next section for the proof.
\end{remark}
To prove Theorem \ref{Thm1-period}, 
we need to prepare the following two lemmas: 
\begin{lemma}
\label{lem:limR2:RC}
Let $f \in C_y^0\left( \T , L_x^1(\R) \right)$, then 
\begin{equation}
\label{limR2:RC}
\lim _{\varepsilon \rightarrow 0} 
2 \pi \varepsilon^{\frac{1}{2}} 
\int_{\R^{2}} f(x,y) 
e^{-\varepsilon \pi y^{2}} d x dy = 
\int_{\R \times \boC} f(x,y) d xdy
\end{equation}
where $\boC:=[- \pi, \pi)$.
\end{lemma}

\begin{proof}
By density, it is sufficient to show \eqref{limR2:RC} for trigonometric polynomials on $y$. In addition, by linearity, we can restrict the proof to $f(x,y)=
f_k(x) e^{i k y}$ for 
$(x,y) \in \R \times \T$. 
Thus, we write
\begin{align*}
\varepsilon^\frac{1}{2} 
\int_{\R^{2}} f(x,y) e^{- \varepsilon \pi y^{2}} d x dy
= \varepsilon^\frac{1}{2} \int_{\R}e^ {i k y} 
e^{- \varepsilon \pi y^{2}} d y \int_{\R} f_k(x)dx 
= e^{- \frac{k^{2}}{4 \varepsilon \pi}} 
\int_{\R} f_k(x) d x 
\end{align*}
Hence,
\[
\lim _{\varepsilon \rightarrow 0} 
\varepsilon^\frac{1}{2} \int_{\R^{2}} f(x,y) 
e^{- \varepsilon \pi y^{2}} d x dy = \left\lbrace 
\begin{aligned}
&0 \quad &{\rm if} \ k\neq 0 \\
& \int_{\R} f_0(x) d x &{\rm if} \ k= 0 
\end{aligned} \right. 
= \frac{1}{2 \pi}\int_{\R \times \boC} f(x,y) d xdy.
\]
This finishes the proof of this lemma.
\end{proof}

\begin{lemma}
\label{lem:limTPQR2:TPQRC}
Let $\mathcal{P} $ and $\mathcal{Q}$ are trigonometric polynomials in $y$. 
For $m \colon \R^{2} \to \mathbb{C}$ 
satisfying \eqref{Cond:m}, 
we define 
\[
T = \boF_{\R^{2}}^{-1}\left[ m \boF_{\R^{2}}\right], \qquad 
\widetilde{T} = 
\boF_{\R \times \T}^{-1}\left[ m \boF_{\R \times \T}\right]
\]
For each $\delta>0$, we 
put $w_{\delta}(y)= e^{- \delta \pi |y|^{2}}$ 
for $y \in \R$. 
Then, one has 
\begin{equation}
\label{limTPQR2:TPQRC}
\lim _{\varepsilon \rightarrow 0} (\pi \varepsilon)^{\frac{1}{2}} 
\int_{\R^{2}} T(\mathcal{P} \left.w_{\varepsilon \alpha}\right)(x,y) \overline{\mathcal{Q}(x,y)} w_{\varepsilon\beta}(y) d x d y 
= \int_{\R \times \boC}
\left(\widetilde{T} \mathcal{P}\right)(x,y) 
\overline{\mathcal{Q}(x,y)} d x d y
\end{equation}
where $\alpha, \beta>0$ satisfying $\alpha+\beta=1$.
\end{lemma}

\begin{proof}
As in the proof of Lemma \ref{lem:limR2:RC}, since \eqref{limTPQR2:TPQRC} is linear in $P$ and $\mathcal{Q} $, 
it is sufficient to consider 
$\mathcal{P} (x,y)=\mathcal{P}_k(x) e^{i k y} $ and 
$\mathcal{Q} (x,y)=\mathcal{Q}_\ell(x) 
e^{i \ell y}$. 
Let $\mathcal{F}_{\R}[f](\xi)$ be 
the Fourier transform of $f \in L^{2}(\R)$, 
that is, 
\[
\mathcal{F}_{\R}[f](\xi) = 
\frac{1}{\sqrt{2\pi}} \int_{\R} e^{- i x \xi} f(x) dx 
\]
for $f \in L^{2}(\R)$.
Using Plancherel's theorem on $\R^{2}$, we write
\begin{equation}
\begin{split}
\label{Plancherel}
& \quad 
\int_{\R^{2}} T\left(\mathcal{P} 
w_{\varepsilon \alpha}\right)(x,y) \overline{\mathcal{Q} }(x,y) w_{\varepsilon\beta}(y) d x d y \\
&= \int_{\R^{2}} m(\xi,\eta) \boF_{\R}[\mathcal{P}_k](\xi) 
\varphi_k(\eta) \overline{
\boF_{\R}[\mathcal{Q}_\ell]
}(\xi) 
\psi_\ell(\eta) d \xi d \eta, 
\end{split}
\end{equation}
where 
$\varphi_k(\eta)$ and $\psi_\ell(\eta)$ are the Fourier transforms in $y$ for $e^ {i k y} w_{\varepsilon \alpha}(y)$ and 
$e^ {i \ell y} w_{\varepsilon \beta}(y)$, respectively. 
This means that
\begin{equation}
\begin{split}
\label{Fourier:phi-psi}
& \varphi_k(\eta) 
= \frac{1}{\sqrt{2\pi}} \int_{\R} e^{- i \eta y} 
e^ {i k y} w_{\varepsilon \alpha}(y) dy 
= (2 \alpha \varepsilon \pi)^{-\frac{1}{2}} 
e^{- \frac{(\eta-k)^{2} }{4 \alpha\varepsilon \pi}}, \\
& \psi_\ell(\eta) 
= \frac{1}{\sqrt{2\pi}} \int_{\R} e^{- i \eta y} 
e^ {i \ell y} w_{\varepsilon \alpha}(y) dy 
= (2 \beta\varepsilon \pi)^{-\frac{1}{2}} 
e^{- \frac{(\eta-\ell)^{2} }
{4 \beta \varepsilon \pi}} .
\end{split}
\end{equation}
\par
We consider the case $k\neq \ell$. 
Including \eqref{Fourier:phi-psi} and the fact that $m$ is bounded into \eqref{Plancherel}, we get
\begin{align*}
& \int_{\R^{2}} T\left(\mathcal{P} w_{\varepsilon \alpha}\right)(x,y) \overline{\mathcal{Q} }(x,y) w_{\varepsilon\beta}(y) d x d y \\
&\lesssim \int_{\R} |\boF_{\R}
[\mathcal{P}_k](\xi) 
\overline{\boF_{\R}[\mathcal{Q}_\ell]}
(\xi)| d \xi 
\int_{\R} (2 \alpha\varepsilon \pi)^{-\frac{1}{2}} 
e^{- \frac{(\eta-k)^{2} }{4 \alpha\varepsilon \pi}} 
(2 \beta \varepsilon \pi)^{-\frac{1}{2}} 
e^{- \frac{(\eta-\ell)^{2} }{4\beta \varepsilon \pi}}d\eta.
\end{align*}
Since $k\neq \ell$, 
we know that $|k- \ell| \geq 1$ and 
get
\begin{align*}
&\int_{\R} (2 \alpha\varepsilon \pi)^{-\frac{1}{2}} 
e^{- \frac{(\eta-k)^{2} }{4\alpha \varepsilon \pi}} 
(2 \beta \varepsilon \pi)^{-\frac{1}{2}} 
e^{- \frac{(\eta - \ell)^{2} }{4 \beta \varepsilon \pi}}d\eta \\
&\lesssim \varepsilon^{-1} 
\left( \int_{|\eta-k| \geq \frac{1}{2} } 
e^{- \frac{(\eta-k)^{2} }{4 \alpha \varepsilon \pi}} 
e^{- \frac{(\eta - \ell)^{2} }{4 \beta \varepsilon \pi}}d\eta 
+ \int_{|\eta - \ell| \geq \frac{1}{2} } 
e^{- \frac{(\eta-k)^{2} }{4 \alpha \varepsilon \pi}} 
e^{- \frac{(\eta- \ell)^{2} }{4 \beta \varepsilon \pi}} 
d\eta \right)\\
&\lesssim 
\varepsilon^{- 1} 
e^{- \frac{1}{16 \alpha\varepsilon \pi}} 
\int_{\R} e^{- \frac{\pi \eta^{2}}{4 \beta \varepsilon}} d\eta 
+ \varepsilon^{- \frac{1}{2}}
e^{- \frac{1}{16 \beta \varepsilon \pi}} 
\int_{\R} e^{- \frac{\pi \eta^{2}}{4 \alpha \varepsilon}}d\eta \\
& \lesssim 
\varepsilon^{-1} e^{- \frac{1}{16 \alpha\varepsilon \pi}} 
+ \varepsilon^{-1}
e^{- \frac{1}{16 \beta \varepsilon \pi}}.
\end{align*}
Hence,
\[
\lim _{\varepsilon \rightarrow 0} \varepsilon^{1 / 2} 
\int_{\R^{2}} T(\mathcal{P} \left.w_{\varepsilon \alpha}\right)(x,y) \overline{\mathcal{Q} }(x,y) w_{\varepsilon\beta}(y) d x d y =0.
\]
On the other side, 
using Plancherel's theorem, we have
\begin{equation}
\begin{split}
\label{Plancherel2}
\int_{\R \times \boC}\left(\widetilde{T} 
\mathcal{P} \right)(x,y) \overline{\mathcal{Q}}(x,y) d x d y 
& 
= \sum_{n \in \mathbb{Z}}
\int_{\R_{x}} m(\xi, n) \boF_{\R}
[\mathcal{P} _k](\xi) \delta_{k, n} 
\overline{\boF_{\R}[\mathcal{Q}_\ell]}(\xi) \delta_{\ell, n} d \xi\\
& 
= \sum_{n \in \mathbb{Z}}
\int_{\R} m(\xi, n) \boF_{\R}[\mathcal{P}_k](\xi)
\overline{\boF_{\R}[\mathcal{Q}_\ell]}(\xi) d \xi 
\delta_{k, n} \delta_{\ell, n} \\
& =0.
\end{split}
\end{equation}
Thus, we conclude \eqref{limTPQR2:TPQRC} for $k\neq \ell$. 
\par
Next, we consider the case when 
$k = \ell$. 
We infer from \eqref{Plancherel} and \eqref{Fourier:phi-psi} that
\begin{align*}
& \int_{\R^{2}} T\left(\mathcal{P} 
w_{\varepsilon \alpha}\right)(x,y) 
\overline{\mathcal{Q} }(x,y) w_{\varepsilon\beta}(y) d x d y \\
&= (4 \alpha \beta \pi)^{-\frac{1}{2}} \varepsilon^{-1} \int_{\R} \boF_{\R}[\mathcal{P}_k](\xi) 
\overline{
\boF_{\R}[\mathcal{Q}_k]}(\xi) \int_{\R} m(\xi,\eta) 
e^{- \frac{(\eta-k)^{2}}{4 \varepsilon \pi}(\frac{1}{\alpha} + 
\frac{1}{\beta})} d\eta d \xi .
\end{align*}
Since $\alpha+\beta=1$, we have 
$\frac{1}{\alpha} + \frac{1}{\beta} = \frac{1}{\alpha\beta}$. 
It follows that 
\begin{align*}
&\lim _{\varepsilon \rightarrow 0} \varepsilon^{\frac{1}{2}} \int_{\R^{2}} T\left(\mathcal{P} w_{\varepsilon \alpha}\right)(x,y) \overline{\mathcal{Q} }(x,y) w_{\varepsilon\beta}(y) d x d y \\
&= \lim _{\varepsilon \rightarrow 0} 
(4 \alpha \beta \varepsilon \pi)^{-\frac{1}{2}} 
\int_{\R} 
\boF_{\R}[\mathcal{P}_k](\xi) 
\overline{
\boF_{\R}[\mathcal{Q}_k]}(\xi) \int_{\R} m(\xi,\eta) 
e^{- \frac{ (\eta-k)^{2}}{4 \varepsilon \alpha \beta \pi}} 
d\eta d\xi \\
&= \int_{\R} \boF_{\R}[\mathcal{P}_k](\xi) \overline{
\boF_{\R}[\mathcal{Q}_k]}(\xi) 
\lim _{\varepsilon \rightarrow 0} 
(4 \alpha \beta \varepsilon \pi)^{-\frac{1}{2}} 
\int_{\R} m(\xi,\eta) e^{- \frac{ (\eta-k)^{2}}
{4 \varepsilon \alpha \beta \pi}} d\eta d \xi \\
&= 
\frac{1}{\sqrt{\pi}} \int_{\R} m(\xi,k) \boF_{\R}[\mathcal{P}_k](\xi) \overline{
\boF_{\R}[\mathcal{Q}_k]}(\xi) d \xi .
\end{align*}
On the other hand, similarly to \eqref{Plancherel2}, we infer that
\[
\int_{\R \times \boC}\left(\widetilde{T} \mathcal{P} \right)(x,y) \overline{\mathcal{Q} }(x,y) d x d y = 
\int_{\R} m(\xi,k) \boF_{\R}[\mathcal{P}_{k}](\xi) \overline{
\boF_{\R}[\mathcal{Q}_k]} (\xi) d \xi.
\]
This finishes the proof of Lemma \ref{lem:limTPQR2:TPQRC}.
\end{proof}
We are now in position to prove Theorem \ref{Thm1-period}.
\begin{proof}[Proof of Theorem \ref{Thm1-period}]
By density, it is sufficient to restrict the proof to trigonometric polynomials on $y$. 
Let $\mathcal{P} $ and $\mathcal{Q} $ be 
the trigonometric polynomials on $y$ and 
$w_{\delta}$ be the function defined in Lemma 
\ref{lem:limTPQR2:TPQRC}. In addition, let 
$r> 1$ be the conjugate exponent of $q$, that is, 
$1/r = 1 - 1/q$.
We put $\alpha=1/q$ and $\beta=1/r$.
Using the 
H\"{o}lder inequality and the continuity 
of the linear operator $T$ (see Theorem \ref{Hormander-thm}),
we write
\begin{align}
\int_{\R^{2}}T\left(\mathcal{P} 
w_{\frac{\varepsilon}{p}}\right)(x,y) 
\overline{\mathcal{Q}(x,y)} 
w_{\frac{\varepsilon}{q}}(y) d x d y 
& \leq \left\| T \mathcal{P} 
w_{\frac{\varepsilon}{p}} \right
\|_{L_{x, y}^q(\R^{2})} \left\| \mathcal{Q} 
w_{\frac{\varepsilon}{q}} 
\right\|_{L_{x, y}^r(\R^{2})}\nonumber \\
& \leq 
\|T\|_{\mathcal{L}(L_{x, y}^{q}(\R^{2}))}
\left\| \mathcal{P} w_{\frac{\varepsilon}{p}} 
\right\|_{L_{x, y}^q(\R^{2})} \left\| \mathcal{Q} 
w_{\frac{\varepsilon}{q}} \right\|_{L_{x, y}^r(\R^{2})} 
\label{TPQR2:estim} 
\end{align}
We multiply both sides in \eqref{TPQR2:estim} by $\varepsilon^{\frac{1}{2}}$ and 
take the limit in order to get
\begin{equation*}
\begin{split}
& \lim_{\varepsilon \to 0 } \varepsilon^{\frac{1}{2}} \int_{\R^{2}}T\left(\mathcal{P} 
w_{\frac{\varepsilon}{p}}\right)(x,y) 
\overline{\mathcal{Q}(x,y)} 
w_{\frac{\varepsilon}{q}}(y) d x d y \\
& \leq \|T\|_{\mathcal{L}
(L_{x, y}^{q}(\R^{2}))}\lim_{\varepsilon \to 0 } 
\varepsilon^{\frac{1}{2}} \left\| \mathcal{P} 
w_{\frac{\varepsilon}{p}} 
\right\|_{L_{x, y}^q(\R^{2})} 
\left\| \mathcal{Q} w_{\frac{\varepsilon}{q}} 
\right\|_{L_{x, y}^r(\R^{2})}.
\end{split}
\end{equation*}
From Lemma \ref{lem:limTPQR2:TPQRC}, we have
\begin{equation}
\begin{split} 
\int_{\R \times \boC}\left(\widetilde{T} 
\mathcal{P} \right)(x,y) 
\overline{\mathcal{Q}}(x,y) d x d y 
& = \lim_{\varepsilon \to 0 } (\pi \varepsilon)^{\frac{1}{2}} \int_{\R^{2}}T\left(\mathcal{P} 
w_{\frac{\varepsilon}{p}}\right)(x,y) 
\overline{\mathcal{Q}(x,y)} 
w_{\frac{\varepsilon}{q}}(y) d x d y\\
& \leq \|T\| \lim_{\varepsilon \to 0 } 
(\pi \varepsilon)^{\frac{1}{2}} 
\left\| \mathcal{P} w_{\frac{\varepsilon}{p}} \right\|_{L_{x, y}^q(\R^{2})} \left\| \mathcal{Q} w_{\frac{\varepsilon}{q}} 
\right\|_{L_{x, y}^r(\R^{2})}.\label{TPQRC:estim} 
\end{split}
\end{equation}
For the right-hand side limit in \eqref{TPQRC:estim}, 
we use Lemma \ref{lem:limR2:RC} in order to obtain
\begin{align*}
&\lim_{\varepsilon \to 0 } 
(\pi \varepsilon)^{\frac{1}{2}} 
\left\| \mathcal{P} w_{\frac{\varepsilon}{p}} \right\|_{L_{x, y}^q(\R^{2})} \left\| \mathcal{Q} w_{\frac{\varepsilon}{q}} 
\right\|_{L_{x, y}^r(\R^{2})}\\
&=\lim _{\varepsilon \rightarrow 0}
\sqrt{\pi} \left[\varepsilon^{\frac{1}{2}} 
\int_{\R^{2} }|\mathcal{P} (x,y)|^{p} 
e^{-\varepsilon \pi |y|^{2}} d x \, dy \right]^{1 / p}
\left[\varepsilon^{\frac{1}{2}} 
\int_{\R^{2} } |\mathcal{Q} (x,y)|^{r} 
e^{-\varepsilon \pi |y|^{2}} d x \, dy \right]^{1 / r} \\
&= \frac{1}{2 \sqrt{\pi}} 
\left\| \mathcal{P} \right\|_{L_{x, y}^q(\R \times \boC )} \left\| \mathcal{Q} \right\|_{L_{x, y}^r(\R \times \boC )}.
\end{align*}
Hence, including this into \eqref{TPQRC:estim}, 
we infer that
\[
\langle \widetilde{T} \mathcal{P}, Q\rangle_{\R \times \T}
= 
\int_{\R \times \boC}
\left(\widetilde{T} \mathcal{P} \right)(x,y) 
\overline{\mathcal{Q}}(x,y) d x d y 
\leq \frac{1}{2 \sqrt{\pi}}
\|T\| \left\| \mathcal{P} \right\|_{L_{x, y}^q(\R \times \boC )} 
\left\| \mathcal{Q} \right\|_{L_{x, y}^r(\R \times \boC )}.
\]
By the duality, we obtain 
$\|\widetilde{T}\|_{\mathcal{L}
(L_{x, y}^{q}\left(\R \times \T \right))} 
\leq \frac{1}{2 \sqrt{\pi}}\|T\|_{\mathcal{L}(L_{x, y}^{q}(\R^{2}))}$. 
This finishes the proof of the theorem.
\end{proof}

\appendix
\section{Density of $y$-trigonometric polynomials}
In this appendix, we shall show that 
the set of the trigonometric polynomial on 
$y$ is dense in $C_{0}^{y}(\T, L_{x}^{1}(\R))$ and 
$L_{x, y}^{q}(\R \times \T)$. 
We can prove this by a classical argument. 
However, 
we will give a proof here for the sake of completeness. 
Note that it suffices to prove 
the density in $C_y^0\left( \T , L_x^1(\R) \right)$ only 
because the density in $L_{x, y}^{q}\left(\R \times \T \right)$ 
follows exactly in the same way. 
We start with reviewing the following definition.
\begin{definition}
A family of functions $\left\{\varphi_{n} \in C^0(\mathbb{T}) 
\colon n \in \mathbb{N}\right\}$ is an approximate identity if:
\begin{eqnarray}
&&\varphi_{n}(y) \geq 0 \quad 
\mbox{for every $y \in \T$} \label{(a)}\\
& &\int_{\mathbb{T}} \varphi_{n}(y) d y=1 \quad 
\mbox{for every $n \in \mathbb{N}$} \label{(b)}\\
&&\lim _{n \rightarrow \infty} \int_{\delta \leq |y| \leq 
\pi} \varphi_{n}(y) d y=0 \quad 
\mbox{for every $\delta>0$}. 
\label{(c)} 
\end{eqnarray}
In \eqref{(c)}
we identify $\T$ with the interval 
$\boC =[- \pi, \pi)$. 
\end{definition}
We now provide the following approximation lemma:
\begin{lemma}
\label{approx-idenity}
Let $f \in C_y^0\left( \T , L_x^1(\R) \right)$ and $\left\{\varphi_{n} \in \boC^0(\T) 
\colon n \in \mathbb{N}\right\}$ be an approximate identity. Then, $\lim_{n \to \infty}\varphi_{n} *_y f = f$ 
in $C_y^0\left( \T , L_x^1(\R) \right)$, 
where 
\[
\varphi_{n} *_{y} f = \int_{-\pi}^{\pi} 
\varphi_{n}(y - t) f(t) dt. 
\]
\end{lemma}
\begin{proof}
From \eqref{(a)} and \eqref{(b)}, we write 
\begin{align*}
&\|\varphi_{n} *_y f(\cdot,y)\|_{L_{x}^1(\R)}-\| f(\cdot,y)\|_{L_{x}^1(\R)}\\
& \qquad\leq \int_\R \int_{-\pi}^{\pi} \varphi_{n} (t) 
|f(x, y-t)|dt\ dx - \int_\R|f(x,y)| \ dx\\
&\qquad\leq \int_{-\pi}^{\pi} \varphi_{n} (t) \left( \int_\R |f(x,y-t)| dx - \int_\R|f(x,y)| \ dx \right) \ dt\\
&\qquad \leq \int_{|t|\leq \delta} \varphi_{n} (t) \left( \| f(\cdot,y-t)\|_{L_{x}^1(\R)} - \| f(\cdot,y)\|_{L_{x}^1(\R)} \right) \ dt\\
& \qquad \qquad + \int_{\pi > |t|> \delta} \varphi_{n} (t) \left( \| f(\cdot,y-t)\|_{L_{x}^1(\R)} - \| f(\cdot,y)\|_{L_{x}^1(\R)} \right) \ dt. 
\end{align*}
On the other hand, since $f \in C_y^0\left( \T , L_x^1(\R) \right)$, we infer that the function $y\mapsto \| f(\cdot,y)\|_{L_{x}^1(\R)} $ is uniformly continuous. Combining this with \eqref{(c)}, we deduce that for any $\varepsilon>0$
\[
\sup_{y \in \T} 
\biggl|\|\varphi_{n} *_y f(\cdot,y)\|_{L_{x}^1(\R)}-\| f(\cdot,y)\|_{L_{x}^1(\R)}\biggl| 
\leq \varepsilon \qquad 
\mbox{for sufficiently large 
$n \in \mathbb{N}$}. 
\]
This finishes the proof of Lemma \ref{approx-idenity}.
\end{proof}

As a consequence, we obtain the density property.
\begin{lemma}
The set of trigonometric polynomials on $y$ are dense in $ C_y^0\left(\T, L_x^1(\R) \right)$.
\end{lemma}

\begin{proof}
For each $n\in \N$, we define a function 
$\varphi_{n}$ by 
\[
\varphi_n(y):=c_n\left(1+\cos y \right)^n,
\]
where 
\[
c_{n} = \left(\int_{-\pi}^{\pi} (1 + \cos y)^{n} dy \right)^{-1}.
\]
Clearly, 
the sequence 
$\{\varphi_{n}\}_{n \in \mathbb{N}}$ satisfies \eqref{(a)} and \eqref{(b)}. 
We claim that $\{\varphi_n\}_{\mathbb{N}}$ also satisfies \eqref{(c)}. 
Putting $t = \tan \frac{y}{2}$. 
Then, we have 
\[
\int_{\delta}^{\pi} (1 + \cos y)^{n} dy 
= \int_{\tan \frac{\delta}{2}}^{\infty} 
\left(\frac{2}{1 + t^{2}} \right)^{n+1} dt, \qquad 
\int_{0}^{\pi} (1 + \cos y)^{n} dy 
= \int_{0}^{\infty} 
\left(\frac{2}{1 + t^{2}} \right)^{n+1} dt. 
\]
We can easily verify that 
\begin{align}
& 
\int_{\tan \frac{\delta}{2}}^{\infty} 
\left(\frac{2}{1 + t^{2}} \right)^{n+1} dt 
\leq \left(\frac{2}{1 + (\tan \frac{\delta}{2})^{2}} \right)^{n} 
\int_{0}^{\infty} 
\frac{2}{1 + t^{2}} dt 
= \left(\frac{2}{1 + (\tan \frac{\delta}{2})^{2}} \right)^{n}
\pi, \label{eq1-lemB2}\\
& \int_{0}^{\infty} 
\left(\frac{2}{1 + t^{2}} \right)^{n+1} dt
\geq \int_{0}^{\frac{\delta}{4}}
\left(\frac{2}{1 + t^{2}} \right)^{n+1} dt 
\geq \left(\frac{2}{1 + \frac{\delta^{2}}{16}} \right)^{n} 
\frac{\delta}{2}. 
\label{eq2-lemB2} 
\end{align}
Since $\tan s \geq s$ for all $s>0$, 
we have, by \eqref{eq1-lemB2} and 
\eqref{eq2-lemB2}, that 
\begin{equation}
\left[\int_{0}^{\infty} 
\left(\frac{2}{1 + t^{2}} \right)^{n+1} dt \right]^{-1}
\int_{\tan \frac{\delta}{2}}^{\infty} 
\left(\frac{2}{1 + t^{2}} \right)^{n+1} dt 
\leq 
\frac{2\pi}{\delta}
\left(\dfrac{1 + \frac{\delta^{2}}{16}}
{1 + \frac{\delta^{2}}{4}} \right)^{n}. 
\label{eq3-lemB2}
\end{equation}
Note that $\varphi_{n}$ is an even function for each $n \in 
\mathbb{N}$. 
This together with \eqref{eq3-lemB2} yields that 
\[
\lim_{n \rightarrow \infty} \int_{\delta \leq |y| \leq 
\pi} \varphi_{n}(y) d y 
\leq \lim_{n \to \infty} 
\frac{4\pi}{\delta}
\left(\dfrac{1 + \frac{\delta^{2}}{16}}
{1 + \frac{\delta^{2}}{4}} \right)^{n} = 0. 
\] 
Therefore, (B.3) holds. 

\par
Hence, 
the sequence 
$\{\varphi_{n}\}_{n \in \mathbb{N}}$ is an approximate identity. 
Thus, from Lemma \ref{approx-idenity}, 
we infer that for any $f \in C_y^0\left( \T , L_x^1(\R) \right)$, $\varphi_{n} *_y f$ converges to $f$ in 
$C_y^0\left(\T, L_x^1(\R) \right)$ 
as $n \rightarrow \infty$. 

\par
It remains then to show that $\varphi_{n} *_y f$ 
is a trigonometric polynomials on $y$. 
We claim that $\varphi_{n}$ is a 
trigonometric polynomial. 
By the binomial theorem, we have 
\begin{equation}
\begin{split}
(1 + \cos y)^{n} 
& = 2^{n} \left(\cos \frac{y}{2} \right)^{2n} \\
& = 2^{n} \dfrac{\left(e^{i \frac{y}{2}} 
+ e^{- i \frac{y}{2}}\right)^{2n}}{2^{2n}} \\
& = 2^{-n} \sum_{k=0}^{2n} 
\begin{pmatrix}
2n \\
k
\end{pmatrix}
e^{\frac{ky}{2}}e^{- \frac{2n-k}{2} y} \\
& = 2^{-n} \sum_{k=0}^{2n} 
\begin{pmatrix}
2n \\
k
\end{pmatrix}
e^{i (k-n)y} 
= 2^{-n} \sum_{k=-n}^{n} 
\begin{pmatrix}
2n \\
k
\end{pmatrix}
e^{i k y}. 
\end{split}
\end{equation}
Thus, we can write
\[
\varphi_{n}(y)=\sum_{k=-n}^{n} a_{n k} e^{i k y}, 
\quad \mbox{where 
$a_{n k}=2^{-n} c_{n}\left(\begin{array}{c}
{2 n} \\
{n+k}
\end{array}
\right)
$
}
\]
Namely, $\varphi_{n}$ is a 
trigonometric polynomial.
This implies that 
\[
\begin{aligned}
\varphi_{n} *_y f(x,y) 
= \int_{-\pi}^{\pi} \sum_{k=-n}^{n} a_{n k} 
e^{i k(y-t)} f(x,t) d t 
&=\sum_{k=-n}^{n} a_{n k} e^{i k y} \int_{-\pi}^{\pi} e^{i k t} 
f(x,t) d t \\
&=\sum_{k=-n}^{n} b_{k}(x) e^{i k y} ,
\end{aligned}
\]
where 
\[
b_{k}(x) =a_{n k} e^{i k y} 
\int_{-\pi}^{\pi} e^{i k t} f(x,t) d t.
\]
This finishes the proof of this lemma.
\end{proof}

\section{Continuity of the minimization value $m_{\omega}$}

\begin{lemma}\label{thmc-1}
Let $p\in (1,5)$. There exists a constant $C( p)>0$ such that if $0< \omega_{1}<\omega_{2} < \infty$, 
and $Q_{\omega_{1}}$ and $Q_{\omega_{2}}$ are minimizers of the variational problems for $m_{\omega_{1}}$ and $m_{\omega_{2}}$, respectively, then, we have 
\begin{align}
\label{eqc-5}
m_{\omega_{1}}
&\le 
m_{\omega_{2}}
-
\frac{\mathcal{M}(Q_{\omega_{2}})}{p+1}(\omega_{2}-\omega_{1}) 
+
C(d, p)\frac{\mathcal{M}(Q_{\omega_{2}})^{2}}{m_{\omega_{2}}}|\omega_{2}-\omega_{1}|^{2},
\\[6pt]
\label{eqc-6}
m_{\omega_{2}}
&\le 
m_{\omega_{1}}
+
\frac{\mathcal{M}(Q_{\omega_{1}})}{p+1}(\omega_{2}-\omega_{1}) 
+
C(d, p)\frac{\mathcal{M}(Q_{\omega_{1}})^{2}}{m_{\omega_{1}}}|\omega_{2}-\omega_{1}|^{2}.
\end{align}
In particular, $m_{\omega}$ is continuous and strictly increasing on 
$(0, \infty)$. 
\end{lemma}
\begin{proof}
Let us begin with a proof of \eqref{eqc-5}. Put $Q_{\omega_{2},\lambda}(x, y) 
:= \lambda Q_{\omega_{2}}(x, y)$ for $\lambda>0$. 
Since $\mathcal{N}_{\omega_{2}}(Q_{\omega_{2}})=0$, 
we see that 
\begin{equation}\label{eqc-1}
\begin{split} 
\mathcal{N}_{\omega_{1}}(Q_{\omega_{2},\lambda})
&=
\lambda^{2}\|\partial_{x} Q_{\omega_{2}}\|_{L^{2}}^{2}
+ \lambda^{2} \||D|_{y}^{\frac{1}{2}} Q_{\omega_{2}}\|_{L^{2}}^{2} 
+ 
\lambda^{2}\omega_{1} \mathcal{M}(Q_{\omega_{2}})
-\lambda^{p+1} \|Q_{\omega_{2}}\|_{L^{p+1}}^{p+1} 
\\[6pt]
&=
\lambda^{2}
\Bigm\{ 
\|Q_{\omega_{2}}\|_{L^{p+1}}^{p+1}
+ (\omega_{1} -\omega_{2}) \mathcal{M}(Q_{\omega_{2}})
-\lambda^{p-1} \|Q_{\omega_{2}}\|_{L^{p+1}}^{p+1}
\Bigm\}. 
\end{split}
\end{equation}
We define $\lambda_{*} < 1$ by 
\begin{equation}\label{eqc-3}
\lambda_{*}
:=
\left(1 + (\omega_{1}-\omega_{2})
\frac{\mathcal{M}(Q_{\omega_{2}})}
{\|Q_{\omega_{2}}\|_{L^{p+1}}^{p+1}}\right)^{\frac{1}{p-1}}.
\end{equation}
so that $\mathcal{N}_{\omega_{1}}(Q_{\omega_{2},\lambda_{*}})=0$. Thus, 
\begin{equation}\label{eqc-2}
m_{\omega_{1}} 
\le 
\mathcal{I}_{\omega_{1}}(Q_{\omega_{2},\lambda_{*}})
=
\lambda_{*}^{2}
\mathcal{I}_{\omega_{1}}(Q_{\omega_{2}})
=
\lambda_{*}^{2}m_{\omega_{2}}.
\end{equation}
The Taylor expansion yields 
that there exists $\theta_{*} \in (0,1)$, 
depending on $\omega_{1}$ and $\omega_{2}$, 
satisfying 
\begin{equation}\label{eqc-4}
\begin{split}
\lambda_{*}^{2}
&=
1+ \frac{2}{p-1} (\omega_{1}-\omega_{2})
\frac{\mathcal{M}(Q_{\omega_{2}})}
{\|Q_{\omega_{2}}\|_{L^{p+1}}^{p+1}} 
\\[6pt]
&\quad 
+
\frac{2(3-p)}{(p-1)^{2}}
\Big\{
1+ \theta_{*} (\omega_{1}-\omega_{2})
\frac{\mathcal{M}(Q_{\omega_{2}})}
{\|Q_{\omega_{2}}\|_{L^{p+1}}^{p+1}}
\Big)^{\frac{2(2-p)}{p-1}}
\Big((\omega_{1}-\omega_{2})
\frac{\mathcal{M}(Q_{\omega_{2}})}
{\|Q_{\omega_{2}}\|_{L^{p+1}}^{p+1}} \Big)^{2}
\\[6pt]
&\le 
1 - \frac{\mathcal{M}(Q_{\omega_{2}})}{(p+1)m_{\omega_{2}}}(\omega_{2} -\omega_{1})
+
C(d, p)\frac{\mathcal{M}(Q_{\omega_{2}})^{2}}{m_{\omega_{2}}^{2}}|\omega_{2}-\omega_{1}|^{2},
\end{split}
\end{equation} 
where $C( p)>0$ is some constant depending only on 
and $p$. In the last estimate, we used 
\[
m_{\omega_{2}} = \mathcal{S}_{\omega_{2}}(Q_{\omega_{2}}) 
= \frac{p-1}{2(p+1)} \|Q_{\omega_{2}}\|_{L^{p+1}}^{p+1}
\] which is a consequence of \eqref{eq-0} since $\mathcal{N}_{\omega_{2}}(Q_{\omega_{2}}) = 0$. 
Combining \eqref{eqc-2} with \eqref{eqc-4}, we obtain the desired inequality 
\begin{equation*}
m_{\omega_{1}}
\le 
m_{\omega_{2}}
-
\frac{\mathcal{M}(Q_{\omega_{1}})}{p+1}(\omega_{2}-\omega_{1}) 
+
C(d, p)\frac{\mathcal{M}(Q_{\omega_{2}})^{2}}{m_{\omega_{2}}}|\omega_{2}-\omega_{1}|^{2} .
\end{equation*} 
Thus, \eqref{eqc-5} holds. 
We can obtain \eqref{eqc-6} similarly. 
This completes the proof. 
\end{proof}

\section{Table of notations}
\begin{center}
\begin{tabular}{|c|c|}
\hline
Symbols 
& Descriptions or equation numbers 
\\ 
\hline
$X$
&
\eqref{energyspace}
\\
$X_{2}$
& 
$X_2 =H^2_xL^2_y \cap L^2_xH^{1}_y(\R \times \T)$
\\
$X_{k}$
& 
$X_k =H^k_xL^2_y \cap L^2_xH^{\frac{k}{2}}_y(\R \times \T)$
\\
$\mathcal{M}$ 
& \eqref{mass}
\\
$\mathcal{H}$ 
& \eqref{hamiltonian}
\\
$\mathcal{S}_{\omega}, \widetilde{\mathcal{S}}_{\omega}, 
\mathcal{S}_{\omega, \mathbb{R}}$ 
& \eqref{action}, \eqref{rescale-action}, \eqref{R-action}
\\
$\mathcal{N}_{\omega}, \widetilde{\mathcal{N}}_{\omega}, 
\mathcal{N}_{\omega, \mathbb{R}}$ 
& \eqref{nehari}, \eqref{rescale-nehari}, \eqref{R-nehari} 
\\
$\mathcal{I}_{\omega}, \widetilde{\mathcal{I}}_{\omega}$ 
& \eqref{eq-0}, \eqref{positive-function}
\\
$R_{\omega}$ 
& 
ground state of \eqref{line-sp}, 
$(2\omega)^{\frac{1}{p-1}} \sech(\sqrt{\omega} x)$ 
\\
$Q_{\omega}, \widetilde{Q}_{\omega}$ 
& 
ground state of \eqref{sp}, \eqref{rescale-GroundState} 
\\
$\mathcal{S}_{\R}, \mathcal{N}_{\R}, m_{\R}, R, Q$ 
& 
$\mathcal{S}_{\R} = \mathcal{S}_{1, \R}, \
\mathcal{N}_{\R} = \mathcal{N}_{1, \R}, \
m_{\R} = m_{1, \R}, \
R = R_{1}, \
Q = Q_{1}
$
\\
$m_{\omega}, \widetilde{m}_{\omega}, m_{\mathbb{R}}$
& \eqref{mini-pro}, \eqref{rescale-mini}, \eqref{R-mini}
\\
$\omega_{p}$
&
$\frac{4}{(p-1)(p+3)}$
\\
$\omega_{*}$
&
given in Theorem \ref{thm-chara-ground}
\\
$\nu_{\omega}$
&
$\frac{\omega}{\omega_{p}}$
\\
$L_{\omega, +}, L_{\omega, -}$ 
& \eqref{1d-ope} 
\\
$L_{\omega, +, n}, L_{\omega, -, n}$
& \eqref{2d-ope}
\\
$L_{\omega, \text{g}, +}$
& \eqref{lop-g}
\\
$S_{\omega}(a)$ 
& \eqref{linerized-ope-1}
\\
$A_{n}$ 
& $L_{\omega_{p}, +, n}$
\\
$J$ 
& \eqref{complex} 
\\
$NL(v, R_{\omega})$
& \eqref{linearized-nonlinear}
\\
$f(z)$
& $f(z) = |z|^{p-1}z$
\\
$F(s)$
& $F(s) = f(sv + R_{\omega})$
\\
$P_{\leq k}$
& \eqref{projection}
\\
$\lambda_{0}$
& positive eigenvalue of 
$- J\mathcal{S}_{\omega}^{\prime \prime}$, 
\eqref{spectral-radius}
\\
$\chi$ 
& eigenfunction of 
$- J\mathcal{S}_{\omega}^{\prime \prime}$ corresponding to $\lambda_{0}$
\\
$\lambda_{2}(a)$
& second eigenvalue of $-\partial_{xx} + |D_{y}| + \omega(a) 
- p \varphi(a)^{p-1}$ 
\\
$\lambda (\omega_{p})$ 
& second eigenvalue of $L_{\omega_{p}, +}$
\\
$\mathcal{P} $, $\mathcal{Q} $
& trigonometric polynomials on $y$
\\
\hline
\end{tabular}
\end{center}

\begin{thank}
The authors would like to thank the anonymous reviewers for their useful comments.
This work was done while H.K. was visiting at University of Victoria. H.K. thanks all members of the Department of Mathematics and Statistics for their warm hospitality. 
Y.B was supported by PIMS grant and NSERC grant (371637-2014).
S.I was supported by NSERC grant (371637-2019).
H.K. was supported by JSPS KAKENHI Grant Number JP17K14223.

\end{thank}

\bibliographystyle{plain}

\vspace{24pt}

\end{document}